\documentclass[12pt,reqno]{amsart}

\usepackage{amsmath,amsthm,amsfonts,amscd,amssymb,euscript}
\numberwithin{equation}{section}
\newcommand{\dbar}{\ensuremath{\bar \partial}}
\newcommand{\dbarb}{\ensuremath{\bar \partial _b }}
\newcommand{\ad}{\ensuremath{\bar \partial _b ^{*}  }}
\newcommand{\adt}{\ensuremath{\bar \partial _{b,t} ^{*}  }}
\newcommand{\adp}{\ensuremath{\bar \partial _b ^{*,\, +}  }}
\newcommand{\adm}{\ensuremath{\bar \partial _b ^{*, \, -}  }}
\newcommand{\ado}{\ensuremath{\bar \partial _b ^{*, \, 0}  }}

\newcommand{\boxbt}{\ensuremath{\square_{b,t}}}
\newcommand{\half}{\ensuremath{\frac {1}{2} \,}}
\newcommand{\C}{\ensuremath{{\mathbb C}}}

\newcommand{\N}{\ensuremath{{\mathbb N}}}

\newcommand{\rb}[1]{\ensuremath{r_{\bar {#1}}}}
\newcommand{\eps}{\epsilon}

\newcommand{\smooth}{\ensuremath{C^{\infty}}}

\newcommand{\lbar}[1]{\ensuremath{\bar L _{#1}}}
\newcommand{\omi}[1]{\ensuremath{{\overline{ \omega}}_{#1}}}
\newcommand{\dzi}[1]{\ensuremath{\frac {\partial} {\partial z_{#1}}}}
\newcommand{\dwi}[1]{\ensuremath{\frac {\partial} {\partial w_{#1}}}}
\newcommand{\dzib}[1]{\ensuremath{\frac {\partial} {\partial \bar z_{#1}}}}
\newcommand{\dwib}[1]{\ensuremath{\frac {\partial} {\partial \bar w_{#1}}}}
\newcommand{\nul}{\ensuremath{\mathcal{N}}}
\newcommand{\nulperp}{\ensuremath{{\,}^\perp \hspace{-0.02em}\mathcal{N}}}
\newcommand{\ran}{\ensuremath{\mathcal{R}}}
\newcommand{\ranperp}{\ensuremath{{\,}^\perp \hspace{-0.02em}\mathcal{R}}}
\newcommand{\hol}{\ensuremath{{\mathfrak H}}}
\newcommand{\harmon}{\ensuremath{{\mathcal{H}}_t}}
\newcommand{\harmof}[1]{\ensuremath{{\mathcal{H}}_{#1}}}
\newcommand{\harperp}{\ensuremath{{\,}^\perp \hspace{-0.02em}{\mathcal{H}}_t}}
\newcommand{\ze}{\ensuremath{\zeta}}
\newcommand{\dze}{\ensuremath{\tilde{\zeta}}}

\newcommand{\ppt}{\ensuremath{{\Psi}_t^+}}
\newcommand{\pot}{\ensuremath{{\Psi}_t^0}}
\newcommand{\pmt}{\ensuremath{{\Psi}_t^-}}
\newcommand{\pptd}{\ensuremath{{\widetilde{\Psi}}_t^+}}
\newcommand{\potd}{\ensuremath{{\widetilde{\Psi}}_t^0}}

\newcommand{\pptg}[1]{\ensuremath{{\Psi}_{t, \, #1}^+}}
\newcommand{\potg}[1]{\ensuremath{{\Psi}_{t, \, #1}^0}}
\newcommand{\pmtg}[1]{\ensuremath{{\Psi}_{t, \, #1}^-}}
\newcommand{\pptdg}[1]{\ensuremath{{\widetilde{\Psi}}_{t, \, #1}^+}}
\newcommand{\potdg}[1]{\ensuremath{{\widetilde{\Psi}}_{t, \, #1}^0}}

\newcommand{\pmtdg}[1]{\ensuremath{{\widetilde{\Psi}}_{t, \, #1}^-}}
\newcommand{\tnorm}[2]{\ensuremath{{\langle| {#1},{#2}|\rangle}_t}}
\newcommand{\sqtnorm}[1]{\ensuremath{{\langle| {#1}|\rangle}^2_t}}

\newcommand{\cp}{\ensuremath{{\mathcal{C}}^+}}
\newcommand{\co}{\ensuremath{{\mathcal{C}}^0}}
\newcommand{\cm}{\ensuremath{{\mathcal{C}}^-}}
\newcommand{\cpd}{\ensuremath{{\widetilde{\mathcal{C}}}^+}}
\newcommand{\cod}{\ensuremath{{\widetilde{\mathcal{C}}}^0}}
\newcommand{\cmd}{\ensuremath{{\widetilde{\mathcal{C}}}^-}}
\newcommand{\qbtg}[1]{\ensuremath{Q_{b,t} ({#1},{#1})}}
\newcommand{\qbtgs}[2]{\ensuremath{Q_{b,t} ({#1},{#2})}}

\newcommand{\qbtp}[1]{\ensuremath{Q^{\, l}_{b,\, +,t} ({#1},{#1})}}
\newcommand{\qbtm}[1]{\ensuremath{Q^{\, l}_{b,\, -,t} ({#1},{#1})}}
\newcommand{\qbto}[1]{\ensuremath{Q^{\, l}_{b,\, 0,t} ({#1},{#1})}}
\newcommand{\qbtos}[2]{\ensuremath{Q^{\, l}_{b,\, 0,t} ({#1},{#2})}}
\newcommand{\qbtk}[2]{\ensuremath{Q^{\, k}_{b,t} ({#1},{#2})}}
\newcommand{\ord}[1]{\ensuremath{\langle {#1}\rangle}}

\begin{document}
\title[Global regularity for $\dbarb$]{Global regularity for $\dbarb$\\ on weakly pseudoconvex CR
manifolds}
\author{Andreea C. Nicoara}



\maketitle

\vspace{-0.37in}

\begin{center}
Harvard University, Department of Mathematics\\
1 Oxford Street, Cambridge, MA 02138 \\
anicoara@math.harvard.edu
\end{center}



\tableofcontents

\section{Introduction}

A CR manifold is essentially the generalization of a real
hypersurface of codimension one in a complex manifold, hence an
odd dimensional manifold whose tangent bundle splits into a
complex subbundle, which is the sum of holomorphic and
anti-holomorphic directions, and another bundle, which is the
totally real part. The tangential Cauchy-Riemann operator $\dbarb$
is then roughly the restriction of the $\dbar$ operator to the
complex subbundle. This paper will only be concerned with CR
manifolds of hypersurface type, namely whose totally real part of
the tangent bundle is a line bundle, these being the most useful
class in the study of functions of several complex variables. We
shall inquire whether $\dbarb$ has closed range on compact,
orientable, weakly pseudoconvex CR manifolds of dimension $2n-1$
embedded in $\C^N,$ $n \leq N.$ The closed range property is
related to existence and regularity results for $\dbarb$ and for
abstract CR manifolds also to questions of embedding. In 1973 in
\cite{Kohnglobalreg}, J. J. Kohn solved the $\dbar$-Neumann
problem on bounded weakly pseudoconvex domains in $\C^n$ using
weights. He then formalized these results, introduced a type of
microlocalization suited for the study of CR manifolds, and proved
 $\dbarb$ had closed range on boundaries of bounded weakly pseudoconvex domains in a
complex manifold in his 1986 paper \cite{Kohnrange}. Mei-Chi Shaw in \cite{Shaw} and Mei-Chi Shaw and Harold Boas
in \cite{BoasShaw} independently proved that the range of $\dbarb$ was closed on forms of degree $n-3$ and $n-2$
respectively, on boundaries of bounded weakly pseudoconvex domains in $\C^n.$ It is in the spirit of \cite{Kohnglobalreg} and
\cite{Kohnrange} that the present paper is written.

We shall extend the results in \cite{Kohnglobalreg} and
\cite{Kohnrange} to compact, orientable, weakly pseudoconvex CR
manifolds embedded in $\C^N$ and of codimension higher than one.
The proof is done using microlocalization by dividing the
cotangent space into two truncated cones and the rest with some
overlap. One of the truncated cones, which we shall call $\cp,$ is
taken around the positive part of the axis corresponding to the
totally real direction of the tangent space, often dubbed the "bad
direction." Naturally, the other truncated cone $\cm$ is taken
around the negative part of that axis. Finally, the rest of the
cotangent space, including a neighborhood of the origin, is often
called the elliptic portion and represents the part on which the
holomorphic directions control the "bad direction," hence we get
very good estimates; we shall denote it by $\co.$ To keep the main
result as general as possible, we make no curvature assumptions on
the manifold $M,$ thus this division of the cotangent space is
only valid locally. Global estimates are then obtained by piecing
together local ones.

The idea that really makes the proof work is a new concept: CR plurisubharmonicity. Its strength lies in the
fact that it is very useful microlocally, and that any plurisubharmonic function on $\C^N$ is also
CR plurisubharmonic. This guarantees that there exists a canonical CR plurisubharmonic function on $\C^N,$ namely the
sum of squares of $|z_i|,$ where $1 \leq i \leq N.$ It is significant to note here that the only part of the proof that
uses the fact the CR manifold is embedded in $\C^N$ is showing the existence of a CR plurisubharmonic function, so in
principle, these results can be stated more abstractly, as long as the existence of a CR plurisubharmonic
function is assumed.

The CR plurisubharmonic function $\lambda$ is then used as a weight function. We impose different weight functions on each
of $\cp$ and $\cm,$ namely $e^{-t \lambda}$ on $\cp$ and $e^{t \lambda}$ on $\cm,$ and we construct very carefully
a covering of the CR manifold so that it satisfies a number of technical conditions. Taking all these into consideration,
we define a microlocal norm ${\langle| \, \cdot \,|\rangle}_t$ using a partition of unity subordinate to that covering.
This norm is equivalent to the $L^2$ norm, but it is much better suited for our purpose. We then define the
energy form $\qbtg{\, \cdot \,}$ with respect to it, which allows us to prove a very powerful main estimate
in Proposition~\ref{global}. Exploiting this main estimate, we can prove the various parts of the following theorem:

\medskip
\newtheorem{maintheorem}{Main Theorem}[section]
\begin{maintheorem}
Let $M^{2n-1}$ be a compact, orientable, weakly pseudoconvex manifold of dimension at least $5,$ embedded in
$\C^N$ ($n \leq N$), of codimension one or above, and endowed with the induced CR structure. Then the following hold:
\begin{enumerate}
\item[(i)] The ranges of $\dbarb,$ its adjoint $\adt$ with respect to ${\langle| \, \cdot \,|\rangle}_t,$ and
the Kohn Laplacian defined as $\boxbt \, = \, \dbarb \, \adt + \adt \,
\dbarb$ are closed in $L^2$ and all $H^s$ spaces for $s>0$ if $t$
is large enough, for $(p,q)$ forms such that $0 \leq q \leq n-3,$ $2 \leq q \leq n-1,$ and $1 \leq q \leq n-2$ respectively;
\item[(ii)] Let the harmonic space $\harmon$ consist of forms in $L^2$ for which the energy form is zero.
Let $\alpha$ be a $(p,q)$ form, $1 \leq q \leq n-2,$ such that $\alpha \perp \harmon.$ If $\alpha \in H^s$
for some $s\geq 0,$ then there exists a positive constant $T_s$ and a unique $(p,q)$ form $\varphi_t \perp \harmon$
such that $$\qbtgs{ \varphi_t}{\phi} =\tnorm{ \alpha}{\phi},$$ $\forall \, \phi \in Dom(\dbarb) \cap Dom(\ad)$
and $\varphi_t \in H^s$  for all $t \geq T_s;$
\item[(iii)] If $\alpha$ is closed $(p,q)$ form, $1 \leq q \leq n-2,$ then there exists a $(p, q-1)$ form $\omega$
such that $$\dbarb \, \omega \, = \, \alpha$$ and if $\alpha \in H^s$ for $s \geq 0,$ then $\omega \in H^s;$
\item[(iv)] If $\alpha \in \smooth$ is closed $(p,q)$ form, $1 \leq q \leq n-2,$ then there exists a $(p, q-1)$
form $\omega \in \smooth$ such that $$\dbarb \, \omega \, = \, \alpha;$$
\item[(v)] Let $H_0^{p,q}(M,\dbarb),$ $H_s^{p,q}(M,\dbarb),$ and $H_\infty^{p,q}(M,\dbarb)$ be the $\dbarb$ cohomology
of $M$ with respect to $L^2,$ $H^s,$ and $\smooth$ coefficients respectively. Then $H_0^{p,q}(M,\dbarb),$
$H_s^{p,q}(M,\dbarb),$ and $H_\infty^{p,q}(M,\dbarb)$ are finite for $1 \leq q \leq n-2$ and
$$H_0^{p,q}(M,\dbarb) \cong H_s^{p,q}(M,\dbarb) \cong H_\infty^{p,q}(M,\dbarb) \quad \forall \: s>0.$$
\end{enumerate}
\end{maintheorem}

\medskip
Part (i) shows that $\dbarb$ and its related operators $\adt$ and
$\boxbt$ all have closed range. Part (ii) provides a weak solution
for $\boxbt$ for datum $\alpha$ that is orthogonal to the harmonic
space. Moreover, the solution has the same regularity as the datum.
Defining the Neumann operator $N_t$ in the standard way and
using this solution, one solves the $\dbarb$ problem in part
(iii), and once again the solution is just as regular as the datum.
A Mittag-Leffler type argument as in \cite{Kohnmethods} shows that
starting with smooth datum, one can find a smooth solution in part
(iv). Finally, in (v) it is a standard result for $(p,q)$ forms with $1 \leq q \leq n-2$ that the $\dbarb$
cohomology groups with $L^2$ and Sobolev $s>0$ coefficients are
finite and isomorphic to each other, given the usual Kohn setup
for the $\dbar$-Neumann problem. A fact not explicitly
written in the literature and which follows
also by a Mittag-Leffler type argument similar to the one used to
prove part (iv) is that the $\dbarb$ cohomology group with smooth
coefficients is isomorphic to the two cohomology groups above,
hence also finite. This is significant since there exists no good
characterization of the $\dbarb$ cohomology group with smooth
coefficients after Michael Christ proved in \cite{Christ} that for worm
domains the canonical solution was not always smooth.

The paper is organized as follows: In Section~\ref{defnot}, we
give all the necessary definitions and introduce all the relevant
notation. Section~\ref{CRpsh} we devote to defining CR
plurisubharmonicity globally, giving its local characterization,
and proving plurisubharmonicity implies CR plurisubharmonicity. In
order to keep the notation simple, in
Section~\ref{microlocsection} we establish the main estimate only
for $(0,1)$ forms, which we use in Section~\ref{regdbar} in order
to prove the main theorem also for $(0,1)$ forms. Then comes
Section~\ref{highforms}, which outlines the proof of the main
theorem of all $(p,q)$ forms with $1 \leq q \leq n-2.$ Finally,
Section~\ref{Gardsect} and Section~\ref{Compdetails} have more or
less the status of appendices. In Section~\ref{Gardsect} we state
the matrix form of the sharp G{\aa}rding inequality, a deep result
in the theory of pseudodifferential operators that is a key
ingredient in the microlocatization, whereas in
Section~\ref{Compdetails} we include details of some computations
from the derivation of the main estimate for the interested
reader.

This work represents the author's PhD thesis written under the direction of J.J. Kohn at Princeton University.

\bigskip
\section{Definitions and Notation}
\label{defnot}

\medskip
\medskip
\newtheorem{def CR}{Definition}[section]
\begin{def CR}
Let $M$ be a smooth manifold of real dimension $2n-1$. An
integrable CR structure is a complex subbundle of the complexified
tangent bundle $\C T(M)$ denoted by $T^{1,0}(M)$ satisfying:
\begin{enumerate}
\item[(i)] $dim_{\C} T_x^{1,0}(M) \, = \, n-1$, where
$T_x^{1,0}(M)$ is the fiber at each $x \in M$;
\item[(ii)] $T^{1,0}(M) \cap T^{0,1}(M) \, = \, \{ 0 \}$,
where $T^{0,1}(M)$ is the complex conjugate of $T^{1,0}(M)$;
\item[(iii)] For any two vector fields $L$ and $L'$ in $T^{1,0}(M)$ defined
on some open subset $U$ of $M$, their bracket given by $[L,L']\, =
\,LL' -L'L$ is a section of $T^{1,0}(M)$ (the integrability
condition).
\end{enumerate}
Such a manifold $M$ endowed with an integrable CR structure is
called a CR manifold. \label{defCR}
\end{def CR}

\smallskip
\noindent We now consider CR manifolds embedded in a complex space
$\C^N$. The complexified tangent bundle of such a CR manifold $M$,
namely $\C T(M),$ is a subbundle of $T(\C^N) \, = \, \C^N$.

\newtheorem{defCRCN}[def CR]{Definition}
\begin{defCRCN}
Let $M$ be a smooth CR manifold of real dimension $2n-1$ embedded
in $\C^N$, $N \geq n$. If the complex subbundle of the
complexified tangent bundle $\C T(M)$ given by $T^{1,0}(M) \, = \,
\C T(M) \cap T^{1,0} (\C^N)$ satisfies the condition (i) of
Definition~\ref{defCR}, then $T^{1,0} (M)$ is called the induced
CR structure on $M$.
\end{defCRCN}

\smallskip
\noindent Notice that in this case conditions (ii) and (iii) of
Definition~\ref{defCR} are automatically satisfied as follows:
$T^{0,1} (M) \, = \, \overline{T^{1,0} (M)} \, = \, \C T(M) \cap
T^{0,1} (\C^N)$ and $T^{1,0} (\C^N) \cap T^{0,1} (\C^N) \, = \, \{
0 \}$, hence (ii). Now take any two vector fields $L$ and $L'$ in
$T^{1,0}(M)$ defined on some open subset $U$ of $M$. By the
definition of $T^{1,0} (M)$, $L \, = \, \sum_{i=1}^N \alpha_i (z)
\dzi {i}$ and $L' \, = \, \sum_{i=1}^N \beta_i (z) \dzi {i}$, so
$[L,L'] \, = \, \sum_{i=1}^N \gamma_i (z) \dzi {i}$, which means
$[L,L']$ is a section of $T^{1,0}(M)$, yielding (iii). We will
consider only CR manifolds embedded in $\C^N$ which have an
induced CR structure on them.

Let $B^q(M)$ be the bundle of $(0,q)$ forms which consists of
skew-symmetric multi-linear maps of $(T^{0,1}(M))^q$ into $\C$.
$B^0 (M)$ is then the set of functions on $M$. Since the induced
CR structure is integrable, meaning it satisfies condition (iii)
of Definition~\ref{defCR}, there exists a natural restriction of
the de Rham exterior derivative to $B^q (M)$, which we will denote
by $\dbarb$. Thus, $\dbarb: \, B^q(M) \rightarrow B^{q+1}(M)$.

\smallskip
\noindent {\bf Note:} The main theorem was stated for $(p,q)$ forms,
but in fact the holomorphic part of the forms is irrelevant for the proof
of the result, hence it is entirely appropriate to discard it and
introduce notation only for $(0,q)$ forms, thus improving the readability
of the paper.

\medskip
\noindent Next, we need to endow $\C T(M)$ with an appropriate
Riemannian metric:

\newtheorem{def CRmetric}[def CR]{Definition}
\begin{def CRmetric}
Let $M$ be a smooth CR manifold of real dimension $2n-1$. A
Riemannian metric is compatible with the CR structure on $M$ if
for all $P \in M$ the spaces $T^{1,0}_P (M)$ and $T^{0,1}_P (M)$
are orthogonal under the Hermitian inner product induced by this
metric on $\C T(M)$.
\end{def CRmetric}

\smallskip
\noindent Since we are considering only CR manifolds embedded in
$\C^N$ which have an induced CR structure on them, it is only
natural to choose as a metric the restriction on $\C T(M)$ of the
usual Hermitian inner product on $\C^N$. Clearly, this metric is
compatible with the induced CR structure because $\langle \dzi{i},
\dzib{j} \rangle \, = \, 0$ for all $1 \leq i,j \leq N$. We then
define a Hermitian inner product on $B^q (M)$ by $$(\phi, \psi) =
\int {\langle \phi, \psi \rangle}_x \, dV,$$ where $dV$ denotes
the volume element and $ {\langle \phi, \psi \rangle}_x$ the inner
product induced on $B^q (M)$ by the metric on $\C T(M)$ at each $x
\in M$. Let $||\, \cdot \,||$ be the corresponding norm and $L_2^q
(M)$ the Hilbert space obtained by completing $B^q (M)$ under this
norm. All the subsequent discussion will only concern the $L^2$
closure of $\dbarb$, so in order to keep the notation simple, we
will also denote it by $\dbarb$. Let us first define its domain:

\smallskip
\newtheorem{def domdbarb}[def CR]{Definition}
\begin{def domdbarb}
$Dom(\dbarb)$ is the subset of $L_2^q (M)$ composed of all forms
$\phi$ for which there exists a sequence of $\{ \phi_\nu \}_\nu$
in $B^q (M)$ satisfying:
\begin{enumerate}
\item[(i)] $\phi \, = \, \lim_{\nu \rightarrow \infty} \phi_\nu$
in $L^2,$ where $\phi_\nu$ is smooth and
\item[(ii)] $\{ \dbarb \phi_\nu \}_\nu$ is a Cauchy sequence in
$L_2^{q+1} (M)$.
\end{enumerate}
\end{def domdbarb}

\noindent For all $\phi \in Dom(\dbarb)$, let $\lim_{\nu
\rightarrow \infty} \dbarb \phi_\nu \, = \, \dbarb \phi$ which is
thus well-defined. We need to define next $\ad$, the $L^2$ adjoint
of $\dbarb$. Again, we first define its domain:

\smallskip
\newtheorem{def domad}[def CR]{Definition}
\begin{def domad}
$Dom(\ad)$ is the subset of $L_2^q (M)$ composed of all forms
$\phi$ for which there exists a constant $C>0$ such that
$$|(\phi, \dbarb \psi)| \leq C ||\psi||$$ for all $\psi \in
Dom(\dbarb).$
\end{def domad}

\noindent For all $\phi \in Dom(\ad)$, we let $\ad \phi$ be the
unique form in $L_2^q (M)$ satisfying $$(\ad \phi, \psi) = (\phi,
\dbarb \psi),$$ for all $\psi \in Dom(\dbarb)$.

\smallskip
The induced CR structure has a local basis $L_1, \dots, L_{n-1}$
for $(1,0)$ vector fields in a neighborhood $U$ of each point $x
\in M$. Let $\omega_1, \dots, \omega_{n-1}$ be the dual basis of
$(1,0)$ forms which satisfies $\langle \omega_i, L_j \rangle \, =
\, \delta_{ij}$, where $\langle \cdot, \cdot \rangle$ is the
natural pairing of $(1,0)$ vector fields with $(1,0)$ forms. This
implies that $\lbar{1}, \dots, \lbar{n-1}$ is a local basis for
the $(0,1)$ vector fields in the neighborhood $U$, and $\omi{1},
\dots, \omi{n-1}$ is the dual basis for $(0,1)$ forms. Moreover,
the tangent space to $M$ in the neighborhood $U$, $T(U)$ is
spanned by $L_1, \dots, L_{n-1}$, $\lbar{1}, \dots, \lbar{n-1}$,
and one more vector $T$ taken to be purely imaginary, i.e.
$\overline T = -T$. Let us now restrict ourselves to orientable CR
manifolds only, and define the Levi form globally. We let $\gamma$
be a purely imaginary global $1$-form on $M$ which annihilates
$T^{1,0}(M)\oplus T^{0,1} (M).$ Notice that any two such forms are
proportional. We normalize $\gamma$ by choosing it in such a way
that $\langle \gamma, T \rangle \, = \, -1.$

\smallskip
\newtheorem{def Levi}[def CR]{Definition}
\begin{def Levi}
The Levi form at a point $x \in M$ is the Hermitian form given by
$\langle d\gamma_x, L \wedge {\overline{L}}' \rangle,$ where $L$
and $L'$ are two vector fields in $T_x^{1,0} (U),$ $U$ a
neighborhood of $x$ in $M.$ We call $M$ weakly pseudo-convex if
there exists a form $\gamma$ such that the Levi form is positive
semi-definite at all $x \in M$ and strongly pseudo-convex if there
exists a form $\gamma$ such that the Levi form is positive
definite again at all $x \in M.$
\end{def Levi}

\smallskip \noindent Notice that this is a stronger definition of
pseudoconvexity than it is usually found in the literature where
$\gamma$ is defined locally, not globally, and pseudoconvexity
means $\langle d\gamma_x, L \wedge {\overline{L}}' \rangle$ does
not change sign in a neighborhood about $x.$ It will become
obvious in the next section the reason orientability and this
stronger definition of pseudoconvexity it permits are necessary.

\medskip
\noindent For the rest of this section, we shall restrict the discussion to functions and $(0,1)$ forms
and postpone the definition of $\dbarb$ and its various adjoints for $(0,q)$ forms with $1 <q \leq n-2$
until Section~\ref{highforms} when it is actually needed.
In this setup, $\dbarb$ is given by the following:
If $u$ is a smooth function on $U$, then
$$\dbarb (u)= \sum_j \lbar {j} (u) \ \omi {j}$$ If $\varphi$ is a
$(0,1)$ form on $U$, $\varphi \, = \, \sum_i \varphi_i \ \omi{i}$,
we have that
$$ \dbarb \  \varphi = \sum_i \dbarb (\varphi_i \wedge \omi{i})
=  \sum_{i<j} (\lbar {i} \varphi_j - \lbar {j} \varphi_i + \sum_k
m_{ij}^k \varphi_k)  \ \omi {i} \wedge \omi {j} \ ,$$ where
$m_{ij}^k$ are in $\smooth (U)$. As for the $L^2$ adjoint, $\ad$,
$$\ad \varphi = \sum_i \lbar{i}^* (\varphi_i) = - \sum_i (L_i
(\varphi_i) + f_i \varphi_i) \ ,$$ where $f_i$ are also in
$\smooth (U)$ and $\lbar{i}^*$ is the $L^2$ adjoint of $\lbar{i}$
for each $i$.

\newtheorem{def dom}[def CR]{Definition}
\begin{def dom}
Let $P$ be a pseudodifferential operator of order zero, then
another pseudodifferential operator of order zero, $\tilde{P}$ is
said to dominate $P$, if the symbol of $\tilde{P}$ is identically
equal to $1$ on a neighborhood of the support of the symbol of $P$
and the support of the symbol of $\tilde{P}$ is slightly larger
than the support of the symbol of $P$.
\end{def dom}

\noindent In particular, the previous definition also applies to
cutoff functions which are pseudodifferential operators of order
zero.

\smallskip

\smallskip
\noindent Since we are working with a parameter $t$, the following
two definitions are in order:

\smallskip
\newtheorem{def tdep}[def CR]{Definition}
\begin{def tdep}
A pseudodifferential operator $P$ is called zero order $t$
dependent if it can be written as $P_1+t \, P_2$, where $P_1$ and
$P_2$ are pseudodifferential operators independent of $t$ and
$P_2$ has order zero.
\end{def tdep}

\smallskip
\newtheorem{def invtdep}[def CR]{Definition}
\begin{def invtdep}
A pseudodifferential operator $P$ of order zero is called inverse
zero order $t$ dependent if its symbol $\sigma(P)$ satisfies
$D^{\alpha}_{\xi} \sigma (P) \, = \, D^{\alpha}_{\xi} p(x, \xi) \,
= \, \frac {1}{t^{|\alpha|}} \, q(x,\xi)$ for $|\alpha| \geq 0$,
where $q(x,\xi)$ is bounded independently of $t$.
\end{def invtdep}

\medskip
We will work with an inner product on $L_2^q (M)$ which will be
the sum of three inner products: the inner product without any
weight defined above which we will denote from now on by ${( \
\cdot \ , \ \cdot \ )}_0$; the inner product with weight
$e^{-t\lambda}$ , ${( \ \cdot \ , \ \cdot \ )}_t \, = \,
{(e^{-t\lambda} \ \cdot \ , \ \cdot \ )}_0$; and the inner product
with weight $e^{t\lambda}$, ${( \ \cdot \ , \ \cdot \ )}_{-t} \, =
\, {(e^{t\lambda} \ \cdot \ , \ \cdot \ )}_0$. Notice that each of
these three inner products determines an $L^2$ adjoint, so we
denote by $\adp$ the $L^2$ adjoint on ${( \ \cdot \ , \ \cdot \
)}_t$, by $\adm$ the $L^2$ adjoint on ${( \ \cdot \ , \ \cdot \
)}_{-t}$, and by $\ado$ the $L^2$ adjoint on ${( \ \cdot \ , \
\cdot \ )}_0$, although since there is no weight function in this
case, $\ado$ equals precisely $\ad$. Now let us compute $\adp$ and
$\adm$. Let $f$ be a smooth function and $\varphi$ a smooth
$(0,1)$ form, then
$$(\dbarb f, \varphi)_t = (\dbarb f, e^{-t \lambda} \varphi)_0 =
(f, \ad (e^{-t \lambda} \varphi))_0 = (f,e^{-t \lambda}(\ad -t [\ad,
\lambda])\varphi)_0 = (f,( \ad -t [\ad,\lambda])\varphi)_t \,
. $$ Therefore, $\adp \varphi =(\ad -t [\ad,\lambda]) \varphi =
\sum_i \lbar {i}^{*, \, t} (\varphi_i) =  - \sum_i (L_i
(\varphi_i) + f_i \varphi_i - t L_i (\lambda)\varphi_i) ,$ where
$\lbar {i}^{*, \, t}$ is the $L^2$ adjoint of $\lbar{i}$ with
respect to ${( \ \cdot \ , \ \cdot \ )}_t$. Similarly,$$(\dbarb f,
\varphi)_{-t} = (\dbarb f, e^{t \lambda} \varphi)_0 = (f, \ad
(e^{t \lambda} \varphi))_0 = (f,e^{t \lambda}( \ad + t[\ad,
\lambda])\varphi)_0 = (f,( \ad + t[\ad,\lambda])\varphi)_{-t} \,
. $$ So then $\adm \varphi = (\ad + t[\ad,\lambda]) \varphi =
\sum_i \lbar {i}^{*, \, -t} (\varphi_i) =  - \sum_i (L_i
(\varphi_i) + f_i \varphi_i + t L_i (\lambda)\varphi_i), $ where
$\lbar {i}^{*, \, -t}$ is the $L^2$ adjoint of $\lbar{i}$ with
respect to ${( \ \cdot \ , \ \cdot \ )}_{-t}$. Notice that the two
adjoint operators $\adp$ and $\adm$ are zero order $t$ dependent,
according to the definition given above.

\smallskip\noindent {\bf Notational remark:} \label{triplenotation} The norms corresponding to the
inner products ${( \ \cdot \ , \ \cdot \ )}_0,$ ${( \ \cdot \ , \
\cdot \ )}_t,$ and ${( \ \cdot \ , \ \cdot \ )}_{-t}$ will be
denoted by $||\, \cdot \, ||_{\, 0},$ $|||\, \cdot \, |||_{\, t},$
and $|||\, \cdot \, |||_{\, -t}$ respectively. The extra bars in
the notation for the two weighted norms are there to distinguish
them from the Sobolev norm that will be introduced in Section
\ref{microlocsection}.

\smallskip We will devote the next section to describing what sort
of function $\lambda$ we need for the upcoming argument.

\bigskip
\section{CR Plurisubharmonicity}
\label{CRpsh}

\medskip
\medskip
We restrict our attention to compact, orientable, weakly
pseudoconvex CR manifolds $M$ of dimension $2n-1$ which are
embedded in $\C^N,$ where $N \geq n$, and have an induced CR
structure. In other words, we will deal not only with CR
hypersurfaces but also with CR manifolds of codimension bigger
than $1$. To accomplish this, we need to define a class of
functions which are plurisubharmonic in the sense that the complex
Hessian defined using the local basis of $(1,0)$ vector fields of
$M$ to which we add a positive multiple of the Levi form has to be
either positive definite or semi-definite. We shall call this
notion CR plurisubharmonicity, and we will prove that any strongly
plurisubharmonic function on $\C^N$ is also strongly CR
plurisubharmonic on an embedded CR manifold $M$, and similarly,
any plurisubharmonic function on $\C^N$ is CR plurisubharmonic.

\medskip
\noindent We start with the invariant definition of CR
plurisubharmonicity, and then we will eventually derive the local
expression for CR plurisubharmonicity stated above. Note that this
invariant definition assumes the Levi form is globally defined,
which, as pointed above, amounts to the CR manifold being
orientable.

\smallskip
\newtheorem{def CRpshglobal}{Definition}[section]
\begin{def CRpshglobal}
Let $M$ be a CR manifold. A $\smooth$ real-valued function
$\lambda$ defined in the \label{CRpshdefglobal} neighborhood of
$M$ is called strongly CR plurisubharmonic if $\, \exists$ $A_0>0$
such that $\langle  \frac{1}{2} \, (\partial_b \dbarb \lambda -
\dbarb \partial_b \lambda)+ A_0 \,d\gamma, L \wedge {\overline{L}}
\rangle$ is positive definite $\, \forall$ $L \in T^{1,0} (M),$
where $\langle d\gamma, L \wedge {\overline{L}} \rangle$ is the
invariant expression of the Levi form. $\lambda$ is called weakly
CR plurisubharmonic if $\, \langle \frac{1}{2} \, (\partial_b
\dbarb \lambda - \dbarb \partial_b \lambda)+ A_0 \, d\gamma, L
\wedge {\overline{L}} \rangle$ is just positive semi-definite.
\end{def CRpshglobal}

\smallskip
\noindent In order to exploit this definition, let us first prove
the existence of a special basis in a neighborhood of each point
of the CR manifold.

\smallskip
\newtheorem{specialbasislemma}[def CRpshglobal]{Lemma}
\begin{specialbasislemma}
Let $M$ be a compact, smooth, weakly pseudoconvex CR manifold of
dimension $2n-1$ embedded in a complex space $\C^N$ such that $N
\geq n$ and endowed with an induced CR structure.
\label{specialbasis} Around each point $P \in M$, there exists a
small enough neighborhood $U$ and a local orthonormal basis $L_1,
\dots, L_n, \lbar{1}, \dots, \lbar{n}$ of the $n$ dimensional complex
bundle containing $T M$ when restricted to $U$, satisfying:
\begin{enumerate}
\item[(i)] $L_i \big|_P \, = \, \dwi{i}$ for $1 \leq i \leq n $, where $(w_1, \dots, w_N)$ are the
coordinates of $\C^N$, and
\item[(ii)] $[L_i, \lbar{j}] \big|_P \, = \, c_{ij} \,T$, where $T \, = \, L_n - \lbar{n}$ and $c_{ij}$
are the coefficients of the Levi form in $L_1, \dots, L_{n-1},
\lbar{1}, \dots, \lbar{n-1}, T$ which is a local basis for $T M$.
\end{enumerate}
\end{specialbasislemma}
\noindent {\bf{Proof:}} Let $P$ be any point of $M$, $P \, = \,
(P_1, \dots, P_N)$. We can define a new system of coordinates
$z_1, \dots, z_N$ in a neighborhood of $P$ in which $M$ has a much
simpler form. First, we translate $P$ to the origin of $\C^N$.
Next, we will show that we can construct local CR coordinates in
which a neighborhood $U$ of $M$ about the origin is given by
$r(z_1, z_2, \dots, z_n) \, = \, 0, \, z_{n+1} \, = \, 0, \,
\dots, \, z_N \, = \, 0,$ where $r$ is a real-valued smooth
function and the outward pointing normal at $0$, $\nabla r(0) \, =
\, (0, \dots, 0, 1, 0, \dots, 0)$, with $1$ in the
$n^{\underline{th}}$ position. $\C T(U)$, the complexified tangent
bundle, has dimension $2n-1$, so the normal bundle $N (U)$ must be
a subbundle of $\C^N$ of dimension $2(N-n)+1$. Now since $T (\C^N)
\, = \, \C^N$, we can take coordinates $z_{n+1}, \dots, z_N$ so
that the corresponding vectors $\dzi{n+1}, \dots, \dzi{N}$ lie in
$N_0 (U),$ the fiber of $N(U)$ at the origin, and complete them to
a system of coordinates on $\C^N$ by coordinates $z_1, \dots, z_n$
chosen so that $\dzi{1}, \dots, \dzi{N}$ is an orthonormal basis
of $\C^N$. Let $\pi$ be the holomorphic projection of $\C^N$ onto
$\C^n$ given by setting $z_{n+1} \, = \, 0, \dots, z_N \, = \, 0$.
$\pi$ is holomorphic, hence it is also a CR map, and it is smooth.
We want to show that we can shrink $U$ so that $\pi$ is a CR
diffeomorphism of $U$ onto its image. Because of the way the $z$
coordinates are chosen and the fact that $M$ is embedded in
$\C^N$, $d \pi (0)$ has rank $2n-1$, full rank on $T_0 (U)$, so it
is an isomorphism of $T_0 (U)$ onto $T_0 (\pi (U))$, hence we can
apply the Inverse Function Theorem to conclude that if we shrink
$U$, $\pi$ is a local diffeomorphism of $U$ onto its image. $\pi$
is a CR map and a local diffeomorphism on $U$, so it must be a CR
diffeomorphism, i.e. its inverse $\pi^{-1}$ is also a CR map. Now
the image $\pi (U)$ has dimension $2n-1$, which means the
coordinates $z_1, \dots, z_n$ are related by some equation, i.e.
some $r(z_1, \dots, z_n) \, = \, 0$. Moreover, $\nabla r$ is the
normal vector to $\pi(U)$, so we can normalize it at the origin to
be precisely $\nabla r(0) \, = \, (0, \dots, 0, 1, 0, \dots, 0)$,
with $1$ in the $n^{\underline{th}}$ position. $r$ has to be a
smooth real-valued function because it is the defining function of
a smooth hypersurface. We can then parametrize $U$, the
neighborhood of $M$ containing $P$ via $\pi^{-1}$ and then indeed
there are local CR coordinates in which $U$ is given by $r(z_1,
z_2, \dots, z_n) \, = \, 0, \, z_{n+1} \, = \, 0, \, \dots, \, z_N
\, = \, 0,$ where $r$ and $\nabla r(0)$ are as specified above. In
order to avoid confusion, we shall denote from now on by $z_1,
\dots, z_n$ the coordinates on $\C^n$ in which $\pi (U)$ is
sitting and by $w_1, \dots, w_N$ the coordinates on $\C^N$ in
which $U$ is sitting. Next we consider a particularly convenient
basis on $T^{1,0}(\pi(U))$ which we will later pull back to
$T^{1,0}(U)$ via $\pi^{-1}$: $\mathfrak{L}'_i \, = \, \dzi{i} -
\frac {r_i}{r_n} \dzi{n},$ where $1 \leq i \leq n-1,$ $r_i \, = \,
\frac {\partial r}{\partial z_i},$ and $r_n \, = \, \frac
{\partial r}{\partial z_n}.$ Let $\mathfrak{L}'_n \, = \, \frac
{1}{r_n} \dzi{n}$, so $\mathfrak{T}' \, = \,\mathfrak{L}'_n - \bar
{\mathfrak{L}}'_n$ is indeed purely imaginary, as needed.
Therefore, $\mathfrak{L}'_1, \dots, \mathfrak{L}'_{n-1}, \bar
{\mathfrak{L}}'_1, \dots, \bar {\mathfrak{L}}'_{n-1},
\mathfrak{T}'$ is a basis on $T (\pi(U))$. Notice that the
condition on $\nabla r(0)$ means $r_i(0) \, = \, 0 $, $\rb{\imath}
(0)\, = \, \frac {\partial r}{\partial \bar z_i}\big|_{\,0} \, = \, 0$, $r_n
(0) \, = \, 1$, and $\rb{n} (0) \, = \, \frac {\partial
r}{\partial \bar z_n}\big|_{\,0} \, = \, 1$ for $1 \leq i \leq n-1$. This
implies that at the origin $\mathfrak{L}'_i \, = \, \dzi{i}$ and
$\mathfrak{T}' \, = \, \dzi{n} - \dzib{n}$. Note that the basis
$\mathfrak{L}'_1, \dots, \mathfrak{L}'_{n-1}, \bar
{\mathfrak{L}}'_1, \dots, \bar {\mathfrak{L}}'_{n-1},
\mathfrak{T}'$ consists of vectors which are neither all
orthogonal to each other, nor normalized. What is very convenient
about the basis $\mathfrak{L}'_1, \dots, \mathfrak{L}'_{n-1}, \bar
{\mathfrak{L}}'_1, \dots, \bar {\mathfrak{L}}'_{n-1},
\mathfrak{T}'$, however, is that $[\mathfrak{L}'_i, \bar
{\mathfrak{L}}'_j] \, = \, \mathfrak{c}_{ij}\, \mathfrak{T}'$ which
is the simplest possible form for this bracket; by
$\mathfrak{c}_{ij}$ we mean of course the coefficients of the Levi
form in this basis. Let us show this is indeed the case:
\begin{equation*}
\begin{split}
[ \mathfrak{L}'_i, \bar {\mathfrak{L}}'_j] &= \mathfrak{L}'_i \bar
{\mathfrak{L}}'_j - \bar {\mathfrak{L}}'_j \mathfrak{L}'_i
=\left(\dzi{i} - \frac {r_i}{r_n} \dzi{n}\right)\left(\dzib{j} -
\frac {\rb{\jmath}}{\rb{n}} \dzib{n}\right) -\left(\dzib{j} -
\frac {\rb{\jmath}}{\rb{n}} \dzib{n}\right)\left(\dzi{i} - \frac
{r_i}{r_n} \dzi{n}\right) \\ &=  \left(\frac {r_i}{r_n} \dzi{n}
\left(\frac {\rb{\jmath}}{\rb{n}}\right)-\dzi{i} \left(\frac
{\rb{\jmath}}{\rb{n}} \right)\right) \dzib{n} + \left(\dzib{j}
\left(\frac {r_i}{r_n} \right)-\frac {\rb{\jmath}}{\rb{n}}
\dzib{n} \left(\frac {r_i}{r_n}\right)\right) \dzi{n} \\ &= \left(
r_{i\bar \jmath} - \frac{r_i r_{n\bar \jmath}}{r_n} -
\frac{\rb{\jmath} r_{i \bar n}}{\rb{n}} + \frac {\rb{\jmath} r_{n
\bar n} r_i}{\rb{n} r_n} \right) \left( \frac {1}{r_n} \dzi{n}-
\frac {1}{\rb{n}} \dzib{n} \right) =
\mathfrak{c}_{ij}\mathfrak{T}'
\end{split}
\end{equation*}
Notice that the normalization condition on $\nabla r(0)$ implies
that $\mathfrak{c}_{ij}(0) \, = \, r_{i \bar \jmath} (0)$.

We now pull back this basis to $T(U)$ via $\pi^{-1}$. We set
$$L'_i = d \pi^{-1} (\mathfrak{L}'_i)$$ $$\vdots$$ $$L'_n = d
\pi^{-1} (\mathfrak{L}'_n)$$ $$T' =  d \pi^{-1} (\mathfrak{L}'_n)-
d \bar \pi^{-1} (\bar {\mathfrak{L}}'_n)$$ Since $\pi^{-1}$ is a
CR map, $L'_1, \dots, L'_{n-1}$ is a basis for $T^{1,0}(U)$, and
$L'_1, \dots, L'_{n-1}, \lbar{1}', \dots, \lbar{n-1}', T'$ is a
basis for $T(U)$. Moreover, $[\mathfrak{L}'_i, \bar
{\mathfrak{L}}'_j] \, = \, \mathfrak{c}_{ij}\, \mathfrak{T}'$
implies $[L'_i, \lbar{j}'] \, = \, c'_{ij}\, T'$, where $c'_{ij}$
are the coefficients of the Levi form in the basis $L'_1, \dots,
L'_{n-1}, \lbar{1}', \dots, \lbar{n-1}', T'$. Notice that in fact
$L'_1, \dots, L'_n, \lbar{1}', \dots, \lbar{n}'$ spans the
$n$ dimensional complex subbundle containing $T M$. Also by construction,
$d \pi^{-1} (0)$ is given by
$$\begin{pmatrix} I_n \\ - \\ 0 \end{pmatrix}$$ namely an
$n \times n$ identity matrix and an $n \times (N-n)$ zero block,
thus $L'_i \big|_P \, = \, \dwi{i}$ for $1 \leq i \leq n$ and
$T' \big|_P \, = \, \dwi{n} - \dwib{n}$. This means that at $P$, the basis $L'_1, \dots, L'_n, \lbar{1}', \dots,
\lbar{n}'$ is orthonormal. Since the lemma requires this to be
true for all points in $U$, we need to apply the Gram-Schmidt
process to $L'_1, \dots, L'_n, \lbar{1}', \dots, \lbar{n}'$. So we
start with $L'_1$ and we let $$L_1 = \frac {1}{||L'_1||}L'_1,$$
$$L_2 = \frac{L'_2 - \langle L_1, L'_2 \rangle L_1}{||L'_2 -
\langle L_1, L'_2 \rangle L_1||} = \sum_{i=1}^2 a_2^i L'_i,$$
$$\vdots$$ $$L_{n-1} = \sum_{i=1}^{n-1} a_{n-1}^i L'_i,$$
 $$L_n = \sum_{i=1}^n a_n^i L'_i.$$ Now let
us look at the bracket. The most general expression of the bracket
is
$$[L_i,\lbar{j}] = c_{ij} T + \sum_{k=1}^{n-1} d_{ij}^k L_k -
\sum_{k=1}^{n-1} \bar d_{ji}^k \lbar{k},$$ since $\overline{[L_j,
\lbar{i}]} = - [L_i, \lbar{j}]$, where $d_{ij}^k$ are $\smooth$
functions and $c_{ij}$ are the coefficients of the Levi form in
the new orthonormal basis. By construction, however, this basis
has the property that at $P$, $L_i \, = \, L'_i \, = \,
\dwi{i}$, hence (i), and $[L'_i, \lbar{j}'] \, = \, c_{ij}' \, T'$
combined with (i) implies $d_{ij}^k \big|_P = \, \bar d_{ji}^k
\big|_P = \, 0$, hence (ii). \qed

\medskip \noindent In the next lemma, we will use the basis
constructed above to prove that CR plurisubharmonicity has the
more convenient local expression claimed at the beginning of this
section.

\smallskip
\newtheorem{lemmaCRpshloc}[def CRpshglobal]{Lemma}
\begin{lemmaCRpshloc}
Let $M$ be a compact, smooth, weakly pseudoconvex CR manifold of
dimension $2n-1$ embedded in a complex space $\C^N$ such that $N
\geq n$ and endowed with an induced CR structure. If $\lambda$
defined on \label{CRpshdefloc} $M$ is strongly CR
plurisubharmonic, then $\, \exists$ $A_0>0$ and a basis $L_1,
\dots, L_{n-1}, \lbar{1}, \dots, \lbar{n-1}, T$ of $T M$ defined
in a neighborhood $U$ of each point $P \in M$ such that the matrix
with entries $s_{ij} \, = \, \frac {1}{2} \, ( \lbar{j}L_i
(\lambda)+ L_i \lbar{j} (\lambda)) +A_0 \, c_{ij}$ is positive
definite in a neighborhood $U'$ around $P$ which is smaller than
$U$, where $c_{ij}$ are the coefficients of the Levi form in the
basis $L_1, \dots, L_{n-1}, \lbar{1}, \dots, \lbar{n-1}, T$. $A_0$
is of course independent of $P$ or $U$, and the size of $U'$
depends on it. If $\lambda$ is weakly CR plurisubharmonic, then
the matrix with entries $s_{ij} \, = \,\frac {1}{2} \,(
\lbar{j}L_i (\lambda)+ L_i \lbar{j} (\lambda)) +A_0 \, c_{ij}$ is
just positive semi-definite.
\end{lemmaCRpshloc}
\smallskip
\noindent {\bf{Proof:}} This lemma has the same hypotheses about
$M$ as the previous one, so from the conclusion of the previous
lemma we know that in a neighborhood $U$ of each point $P \in M$
there exists a local orthonormal basis $L_1, \dots, L_n, \lbar{1},
\dots, \lbar{n}$ of the $n$ dimensional complex bundle containing $T M$
when restricted to $U$, satisfying $[L_i, \lbar{j}] \big|_P \, =
\, c_{ij} \, T$, where $T \, = \, L_n - \lbar{n}$ and $c_{ij}$ are
the coefficients of the Levi form in the basis $L_1, \dots,
L_{n-1}, \lbar{1}, \dots, \lbar{n-1}, T$. We will need this
shortly, but first, since $\lambda$ is strongly CR
plurisubharmonic on $M$ according to
Definition~\ref{CRpshdefglobal}, the following Hermitian form is
positive definite:$$\langle \frac{1}{2} \, (\partial_b \dbarb
\lambda - \dbarb \partial_b \lambda)+ A_0 \,d\gamma, L \wedge
{\overline{L}} \rangle$$  $\, \forall$ $L \in T^{1,0} (M),$ where
$\langle d\gamma, L \wedge {\overline{L}} \rangle$ is the
invariant expression of the Levi form. Let us compute then
$\langle \frac{1}{2} \, (\partial_b \dbarb \lambda - \dbarb
\partial_b \lambda)+ A_0 \,d\gamma, L \wedge {\overline{L}}
\rangle$ in the basis $L_1, \dots, L_{n-1}, \lbar{1}, \dots,
\lbar{n-1}, T$. $$\partial_b \lambda = \sum_{i=1}^{n-1} L_i
(\lambda) \, \omega_i,$$ where $\omega_1, \dots, \omega_{n-1}$ is
the basis of $(1,0)$ forms dual to $L_1, \dots, L_{n-1}$. We then
have $$\dbarb \partial_b \lambda = \sum_{ij} \lbar{j} L_i
(\lambda) \, \omi{j} \wedge \omega_i + \sum_i L_i (\lambda) \,
\dbarb (\omega_i).$$ Let us see what $\dbarb (\omega_i)$ looks
like. $\dbarb (\omega_i)$ is a $(1,1)$ form which in general has
the expression
$$\dbarb(\omega_i) = \sum_{j,k=1}^{n-1} \rho_{\bar \jmath k}^i \omi{j} \wedge
\omega_k,$$ for some functions $\rho_{\bar \jmath k}^i$. Next we
use the Cartan-Frobenius identity that does indeed hold for
$\dbarb$ (see chapter 8 in \cite{Boggess}):
$$\langle \dbarb \, \omega_i, L_i \wedge \lbar{j} \rangle =
L_i(\langle  \omega_i, \lbar{j} \rangle) - \lbar{j}(\langle
\omega_i, L_i \rangle) - \langle \omega_i, [L_i, \lbar{j}]
\rangle.$$ The first two terms of the right-hand side disappear
and $[ L_i, \lbar{j}] \big|_P \, = \, c_{ij}T$, so $$\langle
\dbarb \, \omega_i, L_i \wedge \lbar{j} \rangle \big|_P= \langle
\omega_i, c_{ij}T \rangle = 0.$$ This means $\dbarb \, \omega_i
\big|_P \, = \, 0. $ In other words, $ \rho_{\bar \jmath k}^i
\big|_P \, = \, 0$ for all $1 \leq i,j,k \leq n-1.$ Clearly,
$\partial_b (\omi{i}) \big|_P \, = \, 0$ follows as well. We then
have
$$\dbarb
\partial_b \lambda = \sum_{ij} \lbar{j} L_i (\lambda) \, \omi{j} \wedge \omega_i+ \sum_{ijk} L_i
(\lambda) \rho_{\bar \jmath k}^i \, \omi{j} \wedge \omega_k.$$
Similarly,
$$\dbarb \lambda = \sum_{i=1}^{n-1} \lbar{i} (\lambda) \,
\omi{i}$$ hence
\begin{equation*}
\begin{split}
 \partial_b \dbarb \lambda &= \sum_{ij} L_j  \lbar{i}(\lambda) \,
\omega_j \wedge  \omi{i} + \sum_i \lbar{i} (\lambda) \, \partial_b
(\omi {i}) \\ &= \sum_{ij} L_j  \lbar{i}(\lambda) \, \omega_j
\wedge  \omi{i}+\sum_{ijk} \lbar{i} (\lambda) \bar \rho_{\bar
\jmath k}^i \, \omega_j \wedge  \omi{k} .
\end{split}
\end{equation*}
Since in Definition~\ref{CRpshdefglobal} $L \in T^{1,0} (M),$
$L \, = \, \sum_i \varepsilon_i L_i$ in the local basis which
implies $\overline{L} \, = \, \sum_i {\overline{\varepsilon}}_i
\lbar{i}$. For a general $K \in T^{1,0} (M) \oplus T^{0,1} (M)$,
$K \, = \, (\varepsilon_1, \dots, \varepsilon_{n-1}, \varrho_1,
\dots, \varrho_{n-1})$ and $\bar K \, = \, (\bar \varrho_1, \dots,
\bar \varrho_{n-1}, \bar \varepsilon_1, \dots, \bar
\varepsilon_{n-1})$. This means,
$$\langle \omi{j} \wedge \omega_i, K \wedge \bar K \rangle =
\begin{vmatrix} \varrho_j & \varepsilon_i \\ \bar \varepsilon_j &
\bar \varrho_i \end{vmatrix} = \varrho_j \bar \varrho_i -
\varepsilon_i \bar \varepsilon_j$$ and $$\langle \omega_j \wedge
\omi{i}, K \wedge \bar K \rangle =
\begin{vmatrix} \varepsilon_j & \varrho_i \\ \bar \varrho_j &
\bar \varepsilon_i \end{vmatrix} = \varepsilon_j \bar
\varepsilon_i - \varrho_i \bar \varrho_j.$$ Now  $L \in T^{1,0}
(M),$ so $\varrho_i \, = \, 0$ for $1 \, \leq \, i \, \leq \, n-1$
which implies $$\langle \omi{j} \wedge \omega_i, L \wedge \bar L
\rangle =  - \varepsilon_i \bar \varepsilon_j$$ and $$\langle
\omega_j \wedge \omi{i}, L \wedge \bar L \rangle = \varepsilon_j
\bar \varepsilon_i.$$ The invariant expression of the Levi form
$\langle d\gamma, L \wedge {\overline{L}} \rangle \, = \,
\sum_{ij} c_{ij} \varepsilon_i \bar \varepsilon_j,$ so we have
that
\begin{equation*}
\begin{split}
\langle  \half (\partial_b \dbarb \lambda - \dbarb
\partial_b \lambda)+ A_0 \, d\gamma, L \wedge {\overline{L}} \rangle
& = \half \sum_{ij}  L_j \lbar{i} (\lambda)  \, \varepsilon_j \bar
\varepsilon_i + \half \sum_{ij}  \lbar{j} L_i (\lambda)  \,
\varepsilon_i \bar \varepsilon_j + \sum_{ij} A_0 \, c_{ij} \,
\varepsilon_i \bar \varepsilon_j \\ & \quad + \half \sum_{ijk}
\lbar{i} (\lambda) \bar \rho_{\bar \jmath k}^i  \, \varepsilon_j
\bar \varepsilon_k + \half \sum_{ijk}   L_i (\lambda) \rho_{\bar
\jmath k}^i  \, \varepsilon_k \bar \varepsilon_j.
\end{split}
\end{equation*}
Now exchange $i$ and $j$ in the first sum and gather the terms to
obtain:
\begin{equation*}
\begin{split}
\langle  \half (\partial_b \dbarb \lambda - \dbarb
\partial_b \lambda)+ A_0 \, d\gamma, L \wedge {\overline{L}} \rangle
& = \half \sum_{ij}  L_i \lbar{j} (\lambda)  \, \varepsilon_i \bar
\varepsilon_j + \half \sum_{ij}  \lbar{j} L_i (\lambda)  \,
\varepsilon_i \bar \varepsilon_j + \sum_{ij} A_0 \, c_{ij} \,
\varepsilon_i \bar \varepsilon_j \\ & \quad + \half \sum_{ijk}
\lbar{i} (\lambda) \bar \rho_{\bar \jmath k}^i  \, \varepsilon_j
\bar \varepsilon_k + \half \sum_{ijk}   L_i (\lambda) \rho_{\bar
\jmath k}^i  \, \varepsilon_k \bar \varepsilon_j\\ &=  \sum_{ij}
\: \left[\half ( \lbar{j} L_i (\lambda)+ L_i \lbar{j} (\lambda) )+
A_0 \, c_{ij}\right] \, \varepsilon_i \bar \varepsilon_j\\ & \quad
+ \half \sum_{ijk} \lbar{i} (\lambda) \bar \rho_{\bar \jmath k}^i
\, \varepsilon_j \bar \varepsilon_k + \half \sum_{ijk}   L_i
(\lambda) \rho_{\bar \jmath k}^i  \, \varepsilon_k \bar
\varepsilon_j
\end{split}
\end{equation*}
As derived above, $ \rho_{\bar \jmath k}^i \big|_P \, = \, 0$ for
all $1 \leq i,j,k \leq n-1,$ which means
\begin{equation}
\begin{split}
\langle  \half (\partial_b \dbarb \lambda - \dbarb
\partial_b \lambda)+ A_0 \, d\gamma, L \wedge {\overline{L}}
\rangle\big|_P & = \sum_{ij} \: \left[\half ( \lbar{j} L_i
(\lambda)+ L_i \lbar{j} (\lambda) )+ A_0 \, c_{ij}\right] \,
\varepsilon_i \bar \varepsilon_j \label{invariantCRinlocnot}
\end{split}
\end{equation}
$\lambda$ being strongly CR plurisubharmonic, the matrix with
entries $s_{ij} \, = \, \frac {1}{2} \, ( \lbar{j}L_i (\lambda)+
L_i \lbar{j} (\lambda)) +A_0 \, c_{ij}$ is positive definite at
$P$, hence in a neighborhood $U'$ around $P$ contained in $U$. If
$\lambda$ is just weakly CR plurisubharmonic, the matrix with
entries $s_{ij}$ is merely positive semi-definite. \qed

\medskip \noindent Now, let us show that any strongly
plurisubharmonic function on $\C^N$ is also strongly CR
plurisubharmonic.
\smallskip
\newtheorem{pshmeansCRpsh}[def CRpshglobal]{Lemma}
\begin{pshmeansCRpsh}
Let $M$ be a compact, smooth, weakly pseudoconvex CR manifold of
dimension $2n-1$ embedded in a complex space $\C^N$ such that $N
\geq n$ and endowed with an induced CR structure. Furthermore, let
$\lambda$ be strongly \label{pshisCRpsh}
 plurisubharmonic on $\C^N$, then $\lambda$ is also strongly CR plurisubharmonic.
\end{pshmeansCRpsh}
\noindent {\bf{Proof:}} We apply Lemma~\ref{specialbasis} to get
the particular local basis in which
Equation~\ref{invariantCRinlocnot} holds at point $P \in M$, i.e.
\begin{equation*}
\begin{split}
\langle  \half (\partial_b \dbarb \lambda - \dbarb
\partial_b \lambda)+ A_0 \, d\gamma, L \wedge {\overline{L}}
\rangle\big|_P & = \sum_{ij} \: \left[\half ( \lbar{j} L_i
(\lambda)+ L_i \lbar{j} (\lambda) )+ A_0 \, c_{ij}\right] \,
\varepsilon_i \bar \varepsilon_j
\end{split}
\end{equation*}

\noindent Now in coordinates
$$L_i = \sum_{k=1}^N a^k_i \dwi{k}$$ which means $$L_i \lbar{j}-
\sum_{k,l=1}^N a^k_i \bar a^l_j \frac{\partial^2}{\partial w_k
\partial \bar w_l} = \sum_{k,l=1}^N a^k_i \dwi{k}(\bar a^l_j)
\dwib{l} \in T^{0,1}( \C^N)$$ and $$ \lbar{j} L_i- \sum_{k,l=1}^N
a^k_i \bar a^l_j \frac{\partial^2}{\partial w_k \partial \bar w_l}
= \sum_{k,l=1}^N \bar a^l_j  \dwib{l}(a^k_i) \dwi{k} \in T^{1,0}(
\C^N).$$ According to Lemma~\ref{specialbasis} (ii), $\, [L_i,
\lbar{j}] \big|_P \, = \,L_i \lbar{j} - \lbar{j}L_i \, = \, c_{ij}
T \, = \, c_{ij} (L_n - \lbar{n})$ so $$ \sum_{k,l=1}^N \bar a^l_j
\dwib{l}(a^k_i) \dwi{k} \Big|_P= - c_{ij} L_n$$ and
$$\sum_{k,l=1}^N a^k_i \dwi{k}(\bar a^l_j) \dwib{l} \Big|_P = - c_{ij}
\lbar{n}.$$ We can then compute
$$\lbar{j}(\lambda (w)) = \sum_{k=1}^N \bar a^k_j \frac{\partial
\lambda}{\partial \bar w_k} ,$$ so
\begin{equation*}
\begin{split}
L_i \lbar{j} (\lambda (w)) \big|_P &=  \sum_{k,l=1}^N a^k_i \bar
a^l_j \frac{\partial^2 \lambda}{\partial w_k \partial \bar
w_l}\Big|_P + \sum_{k,l=1}^N  a^k_i \dwi{k}(\bar a^l_j) \dwib{l}
(\lambda)\Big|_P\\ & =\sum_{k,l=1}^N a^k_i \bar a^l_j
\frac{\partial^2 \lambda}{\partial w_k \partial \bar w_l} \Big|_P-
c_{ij} \lbar{n}(\lambda) \big|_P
\end{split}
\end{equation*}
From Lemma~\ref{specialbasis} (i), $L_i \big|_P \, = \, \dwi{i}$
for $1 \leq i \leq n$ which means
$$L_i \lbar{j} (\lambda (w))\big|_P =  \frac{\partial^2 \lambda}
{\partial w_i \partial \bar w_j}\Big|_P-c_{ij}(P)L_n (\lambda)
\big|_P.$$ Similarly,
$$\lbar{j} L_i (\lambda(w))\big|_P = \frac{\partial^2 \lambda}
{\partial w_i \partial \bar w_j}\Big|_P-c_{ij}(P)\lbar{n}
(\lambda) \big|_P.$$ Therefore, we have that $$\half (\lbar{j} L_i
(\lambda(w))+L_i \lbar{j} (\lambda (w)))\big|_P = \frac{\partial^2
\lambda}{\partial w_i \,
\partial \bar w_j}(w)\Big|_P - c_{ij}(P) \frac {L_n
(\lambda)+ \lbar{n} (\lambda)}{2}\Big|_P$$ which means that if we
take $A_0 \,
> \,  \frac {L_n
(\lambda)+ \lbar{n} (\lambda)}{2}\big|_P$, the Hermitian form
$$\langle \frac{1}{2} \, (\partial_b \dbarb \lambda - \dbarb
\partial_b \lambda)+ A_0 \,d\gamma, L \wedge {\overline{L}}
\rangle$$ will be positive definite at $P$ for any $L \in T^{1,0}
(M)$ hence likewise in an open subset of $U$ containing $P$
because $\lambda$ is strongly plurisubharmonic and $M$ is weakly
pseudoconvex. At each $P \in M$ such a construction works, and we
get a neighborhood $U_P$ and a constant $A_0^P$. $\{ U_P \}_{P \in
M}$ is an open cover of $M$ which is compact, so we get a constant
$A_0$ which renders the Hermitian form $$\langle \frac{1}{2} \,
(\partial_b \dbarb \lambda - \dbarb
\partial_b \lambda)+ A_0 \,d\gamma, L \wedge {\overline{L}}
\rangle$$ positive definite everywhere on $M$. \qed

\smallskip
Clearly, the proof above also shows that any plurisubharmonic
function in the neighborhood of $M$ is also CR plurisubharmonic on
it, but we get for free something even stronger. On $\C^N$ there
exists a canonical strongly plurisubharmonic function, $\lambda
(z) \, = \, \sum_{i=1}^N |z_i|^2$, so the fact that $M$ is
embedded in $\C^N$ guarantees that it has at least one strongly CR
plurisubharmonic function on it by the previous lemma.

For the microlocalization that will be done in the next section,
we need to work with an orthonormal basis which should also have
the property that the bracket of an element of $T^{1,0}$ with an
element of $T^{0,1}$ has a fairly simple expression. The most
convenient such local basis is the one given by the conclusion of
Lemma~\ref{specialbasis}. We thus know that around each point $P
\in M$, there exists a neighborhood $U$ and an orthonormal basis
of vector fields defined on it, $L_1, \dots, L_n, \lbar{1}, \dots,
\lbar{n},$ such that
\begin{equation}
\begin{split}
[L_i,\lbar{j}] &= c_{ij} T + \sum_{k=1}^{n-1} d_{ij}^k L_k -
\sum_{k=1}^{n-1} \bar d_{ji}^k \lbar{k}, \label{primitivebracket}
\end{split}
\end{equation}
 $1 \leq i,j \leq n-1$, $T \, = \, L_n - \lbar{n},$ and
$$[L_i,\lbar{j}] \big|_P = c_{ij} T .$$This means we can show something
stronger. We proceed as follows: In the previous section we
introduced the adjoint of $\lbar{p}$ with respect to the norm with
weight function $e^{-t \lambda}$ which we denoted by $\lbar{p}^{*,
\, t}$ and that with respect to the norm with weight function
$e^{t \lambda}$ denoted by $\lbar{p}^{*, \, -t}$. To take care of
both cases at once, let
$$\lbar{p}^{*, \, \pm t} = - L_p - f_p \pm t L_p (\lambda),$$
where $f_p$ is a $\smooth$ function independent of $t$. We solve
for $L_p$ and plug that into the expression for $[L_i,\lbar{j}]$
from above:
\begin{equation*}
\begin{split}
[L_i,\lbar{j}] &= c_{ij} T - \sum_{k=1}^{n-1} \bar d_{ji}^k
\lbar{k} - \sum_{k=1}^{n-1} d_{ij}^k \lbar{k}^{*, \, \pm t} -
\sum_{k=1}^{n-1} d_{ij}^k f_k \pm t \sum_{k=1}^{n-1} d_{ij}^k
L_k(\lambda)
\end{split}
\end{equation*}
Since $d_{ij}^k \big|_P = \, \bar d_{ji}^k \big|_P  = \, 0$, we
conclude that if we shrink the neighborhood $U$ sufficiently,
there exist $\smooth$ functions $a_{ij}^k,$ $b_{ij}^k,$ $g_{ij},$
and $e_{ij}$ independent of $t$ such that
\begin{equation}
\begin{split}
[L_i,\lbar{j}] &= c_{ij} T + \sum_k a_{ij}^k \lbar{k} + \sum_k
b_{ij}^k \lbar{k}^{*, \, \pm t}+ g_{ij} \pm t \, e_{ij},
\label{bracket}
\end{split}
\end{equation}
where $|e_{ij}|$ is bounded independently of $t$ by a small
positive constant which we denote by $\eps_G$, the subscript
indicating that it is a global quantity: $$|e_{ij}| \leq \eps_G.$$
Let us now fix some small $\eps_G>0$ and the function $\lambda$ to
be the canonical strongly plurisubharmonic function, namely
$\lambda (z) \, = \, \sum_{i=1}^N |z_i|^2$, which is strongly CR
plurisubharmonic by Lemma~\ref{pshisCRpsh}. Around each point $P
\in M$ we have a neighborhood $V$ such that Equation \ref{bracket}
holds. Notice that Equation \ref{bracket} uses the special basis
whose existence is proven in Lemma~\ref{specialbasis}, so each $V$
is small enough that it can be parametrized by a hypersurface in
$\C^n$, as done in the proof of the above-mentioned lemma. We will
employ the covering $\{ V_\mu \}_\mu$ for the microlocalization
that follows.

\bigskip
\section{Microlocalization Result}
\label{microlocsection}

\medskip
\medskip
As mentioned in the introduction, the main proposition will be proven using microlocalization,
namely we will divide the Fourier transform space into three
conveniently chosen regions, two truncated cones $\cp$ and $\cm$
and another region $\co$, with some overlap. Let the coordinates
on the Fourier transform space be $\xi \, = \, (\xi_1, \xi_2,
\dots, \xi_{2n-2}, \xi_{2n-1})$. Write $\xi' \, = \,(\xi_1, \xi_2,
\dots, \xi_{2n-2})$, so then $\xi \, = \, (\xi', \xi_{2n-1})$. The
work is done in coordinate patches on $M,$ each of which has
defined on it local coordinates such that $\xi'$ is dual to the
holomorphic part of the tangent bundle $T^{1,0}(M)\oplus T^{0,1}
(M)$ and $\xi_{2n-1}$ is dual to the totally real part of the
tangent bundle of $M$ spanned by the "bad direction," $T.$ Define
$\cp \, = \, \{ \xi \, | \, \xi_{2n-1} \geq \frac {1}{2} |\xi'| \:
and \: |\xi| \geq \ 1 \}$. Then $\cm \, = \, \{ \xi \, | \, -\xi
\in \cp \}$, and finally $\co \, = \, \{ \xi \, | \, -\frac {3}{4}
|\xi'| \leq \xi_{2n-1} \leq \frac {3}{4} |\xi'| \} \cup \{ \xi \,
| \, |\xi| \leq 1 \}$. Notice that by definition, $\cp$ and $\co$
overlap on two smaller cones and part of the sphere of radius $1$
and similarly  $\cm$ and $\co$, whereas $\cp$ and $\cm$ do not
intersect.

Let us now define three functions on $\{ |\xi'|^2 + |\xi_{2n-1}|^2
\, = \, 1 \}$, which is the unit sphere in $\xi$ space. $\psi^+$,
$\psi^-$, and $\psi^0$ are smooth, take values in $[0,1]$, and
satisfy the condition of symbols of pseudodifferential operators
of order zero. Moreover, $\psi^+$ is supported in $\{ \xi \,  | \,
\xi_{2n-1} \geq \frac {1}{2} |\xi'| \}$ and $\psi^+ \, \equiv \,
1$ on the subset $\{ \xi \, | \, \xi _{2n-1} \geq \frac {3}{4}
|\xi'| \}$. Let then $\psi^- (\xi)\, = \, \psi^+ (-\xi)$ which
means that $\psi^-$ is supported in $\{ \xi \,  | \, \xi_{2n-1}
\leq - \frac {1}{2} |\xi'| \}$ and $\psi^- \, \equiv \, 1$ on the
subset $\{ \xi \, | \, \xi _{2n-1} \leq -\frac {3}{4} |\xi'| \}$.
Finally, let $\psi^0 (\xi)$ satisfy $(\psi^0 (\xi))^2 \, = \, 1 -
(\psi^+ (\xi))^2 - (\psi^- (\xi))^2$ which means that $\psi^0$ is
supported in $\co$ and $\psi^0 \, \equiv \, 1$ on the subset $\{
\xi \, | \, -\frac {1}{2} |\xi'| \leq \xi_{2n-1} \leq \frac {1}{2}
|\xi'| \}$. Next extend $\psi^+$, $\psi^-$, and $\psi^0$
homogeneously by setting $\psi^+ (\xi) \, = \, \psi^+ (\frac
{\xi}{|\xi|})$, for $\xi$ in $\cp$ outside of the unit sphere.
Similarly, let $\psi^- (\xi) \, = \, \psi^- (\frac {\xi}{|\xi|})$,
for $\xi$ in $\cm$ outside of the unit sphere and $\psi^0 (\xi) \,
= \, \psi^0 (\frac {\xi}{|\xi|})$, for $\xi$ in $\co$ outside of
the unit sphere. Extend $\psi^+$, $\psi^-$, and $\psi^0$ inside
the unit sphere in some smooth way so that $(\psi^+)^2 + (\psi^-
)^2+ (\psi^0)^2 \, = \, 1$ still holds. Now we have the functions
$\psi^+$, $\psi^-$, and $\psi^0$ defined everywhere on $\xi$
space. We do, however, want to define symbols of
pseudodifferential operators which can produce factors of $t$ in
the denominator under differentiation in order to cancel factors
of $t$ coming from the differentiation of the weight functions
which we will be using, so we define $\psi^+_t (\xi) \, = \,
\psi^+ (\frac {\xi}{tA})$, $\psi^0_t (\xi) \, = \, \psi^0 (\frac
{\xi}{tA})$, and $\psi^-_t (\xi) \, = \, \psi^- (\frac {\xi}{tA})$
for some positive constant $A$ to be chosen later. Now let $\ppt,$
$\pot,$ and $\pmt$ be the pseudodifferential operators of order
zero with symbols $\psi^+_t$, $\psi^-_t$, and $\psi^0_t$
respectively. By construction, the following is true: $$(\ppt)^*
\ppt + (\pot)^* \pot  +(\pmt)^* \pmt = Id$$ In the absence of
curvature assumptions on the manifold $M,$ such construction
cannot hold globally. Instead, it is only true in a neighborhood
small enough that it can be parametrized by a piece of a
hypersurface in $\C^n$ as in Lemma~\ref{specialbasis}. Let $\{
V_\mu \}_\mu$ be the open covering of $M$ mentioned at the end of
Section~\ref{CRpsh} in which each $V_\mu$ is such a neighborhood.
It follows that the operators $\ppt,$ $\pot,$ and $\pmt$ make
sense only when applied to a function or form restricted to a
neighborhood $V_\mu$ in the covering $\{ V_\mu \}_\mu$.
Furthermore, after one of these operators is applied to a function
or form restricted to $V_\mu,$ the resulting function or form
spreads outside of $V_\mu$ and has to be controlled by multiplying
it by another cutoff function.

Keeping this in mind, we will find an appropriate open cover $\{
U_\nu \}_\nu$ for $M$ which we will use in order to construct a
convenient global norm, ${\langle| \, \cdot \, |\rangle}_t.$ We
will take a smooth partition of unity $\{ \ze_\nu \}_\nu$
subordinate to it, and for each $\ze_\nu,$ a cutoff function
$\dze_\nu$ that dominates it and is also supported in $U_\nu.$ By
$\varphi^\nu$ we will understand a $(0,1)$ form $\varphi$
expressed in the local coordinates on $U_\nu.$ On each $U_\nu$ we
will define pseudodifferential operators $\pptg{\nu},$
$\potg{\nu},$ and $\pmtg{\nu},$ as above. Moreover, for each
$\potg{\nu},$ we will take a pseudodifferential operator
$\potdg{\nu}$ whose symbol dominates that of $\potg{\nu}$ and is
supported in a slightly bigger region of the Fourier transform
space. Let $\adt$ be the adjoint of $\dbarb$ with respect to the
inner product ${\langle| \, \cdot \, , \, \cdot \, |\rangle}_t.$
We define $\qbtg {\varphi} \, = \, \sqtnorm{\dbarb \varphi} +
\sqtnorm{\adt \varphi}.$ Finally, we can state the main estimate
which will be proven close to the end of this section:
\medskip
\newtheorem{lemma}{Proposition}[section]
\begin{lemma}
Let $\varphi$ be a $(0,1)$ form supported on $M$, a compact,
orientable, weakly pseudoconvex \label{global} CR manifold of
dimension at least $5$, and let $M$ be endowed with a strongly CR
plurisubharmonic function $\lambda$ as described in
Section~\ref{CRpsh}. If $\varphi$ satisfies $\varphi \in Dom(\dbarb)
\cap Dom(\ad),$ then there exist constants $K,$ $K_t,$ and
$K'_t,$ and a positive number $T_0$ such that for any $t \geq
T_0$,
\begin{equation}
\label{globaleq} K \qbtg {\varphi}+K_t \sum_\nu \:{||\dze_\nu
\potdg {\nu} \ze_\nu \varphi^\nu||}^2_0 +K'_t \,||\varphi||^2_{-1}
\geq  t {\langle |\varphi|\rangle}^2_t.
\end{equation}
\end{lemma}
\noindent {\bf Note:} Since $\varphi$ is a $(0,1)$ form, it
consists of components $\varphi_i$, but in order to save notation,
we write ${||\dze_\nu \potdg {\nu} \ze_\nu \varphi^\nu||}^2_0$
instead of $\sum_i {||\dze_\nu \potdg {\nu} \ze_\nu
\varphi^\nu_i||}^2_0.$

\bigskip \noindent Going back to the construction of an open
covering for $M,$ let $\{ U_\nu \}_\nu$ be an open cover with the
property that for each $\nu$, $U_\nu$ and all $U_\eta$ such that
$U_\nu \cap U_\eta \, \neq \, \emptyset$ are contained in some
neighborhood $V_\mu$ of $M$ in the open cover $\{ V_\mu \}_\mu,$
the covering mentioned at the end of Section~\ref{CRpsh}. This
property is essential for patching up local estimates via changes
of variables applied to pseudodifferential operators defined in
overlapping sets of the covering of $M$. It is not immediately
obvious such a covering $\{ U_\nu \}_\nu$ exists given the
covering $\{ V_\mu \}_\mu$, so let us prove this fact in the
following lemma:

\medskip
\newtheorem{covlemma}[lemma]{Lemma}
\begin{covlemma}
Let $M$ be a compact Riemannian manifold and let $\{ V_\mu \}_\mu$
be a given open covering of $M$ in the topology induced by $M$'s
metric, then there exists another (finite) open covering of $M$,
$\{ U_\nu \}_\nu$, such that for each $\nu,$ $\exists \: \mu(\nu)$ with the property that $U_\nu$ and all
$U_\eta$ satisfying $U_\nu \cap U_\eta \, \neq \, \emptyset$ are
contained in the neighborhood $V_{\mu(\nu)}$ of $\{ V_\mu \}_\mu$.
\label{coveringlemma}
\end{covlemma}

\smallskip
\noindent {\bf{Proof:}} Since $M$ is a Riemannian manifold, it can
be made into a metric space as follows: Let $\gamma$ be any path
in $M$; we denote by $|| \gamma||$ the length of $\gamma$ in the
metric of $M$. Then for any $x, \: y \in M$ we define $$||x-y|| =
\inf_\gamma \{ || \gamma|| \: \big| \: \gamma:[0,1]\mapsto M, \:
\gamma (0) = x \: and \: \gamma(1) = y \}$$ It is easy to check
this distance function satisfies all the required properties. The
topology induced by the metric of $M$ is then the topology induced
by this distance function. We will denote by $B_{r_P} (P)$ the
open ball in this topology centered at $P \in M$ and of radius
$r_P$.

Since $\{ V_\mu \}_\mu$ is an open cover of $M$, for any point $P$
in $M$, there exists a $\mu$ such that $P \in V_\mu$. Moreover,
let $r_P$ be the radius of the largest ball centered at $P$
contained in $V_\mu$. Consider the open covering $\{B_{r_P}
(P)\}_{P \in M}$. Since $M$ is compact, there exists a finite
subcover $\{B_{r_{P_i}} (P_i)\}_{i=1}^\Omega$. Without loss of
generality, we can assume this covering is minimal, i.e. $\forall
\, j,$ $\, 1 \leq j \leq \Omega,$ $\, \exists \, P \in M$ such
that $P \not\in \, \bigcup_{i \neq j} \, B_{r_{P_i}} (P_i);$
otherwise, we start with $j=1$ and go through $1 \leq j \leq
\Omega,$ eliminating the set $B_{r_{P_j}} (P_j)$ if $M$ is
contained in $\bigcup_{i \neq j} \, B_{r_{P_i}} (P_i)$ until the
covering becomes minimal.

Next, we need to find a positive number $\rho$ by which we can
shrink each $r_{P_i}$ such that $\{B_{r_{P_i}-\rho}
(P_i)\}_{i=1}^\Omega$ is still an open covering of $M$. Start with
$j=1.$ $\, M \setminus \bigcup_{i \neq 1} \, B_{r_{P_i}} (P_i)$ is
a non-empty closed set since the covering $\{B_{r_{P_i}}
(P_i)\}_{i=1}^\Omega$ is minimal. There exists then some $r'_{P_1}
< r_{P_1}$ such that the open ball of radius $r'_{P_1}$ centered
at $P_1$ contains $M \setminus \bigcup_{i \neq 1} \, B_{r_{P_i}}
(P_i)$. Let $\rho_1 = r_{P_1} - r'_{P_1}$. We replace $B_{r_{P_1}}
(P_1)$ by $B_{r'_{P_1}} (P_1)$ in the covering $\{B_{r_{P_i}}
(P_i)\}_{i=1}^\Omega$, and then we do the same for $j=2$ and so on
up to $\Omega$. Finally,
$$\rho = \min_{j=1}^\Omega \rho_j.$$ Now let $r = \frac {\rho}{4}$. The open covering
$\{B_r (P)\}_{P \in M}$ has the property that if $B_r (P) \subset
V_\mu$ for some $\mu$, then for any $P' \in M$ such that $B_r (P)
\cap B_r (P') \neq \emptyset$, $B_r (P') \subset V_\mu$. Indeed,
$\exists j$ such that $P \in B_{r'_{P_j}} (P_j)$ because
$\{B_{r'_{P_i}} (P_i)\}_{i=1}^\Omega$ covers $M$, so $B_r (P)
\subset B_{r_{P_j}} (P_j) \subset V_\mu$ for some $\mu$ by
construction. Also by construction, for all $P' \in M$ such that
$B_r (P) \cap B_r (P') \neq \emptyset$, $B_r (P') \subset
B_{r_{P_j}} (P_j) \subset V_\mu$. Take then a finite subcover of
$M$ from $\{B_r (P)\}_{P \in M}$ and let that be the covering $\{
U_\nu \}_\nu$.\qed

\medskip
We only consider compact CR-manifolds embedded in $\C^N$, hence
Riemannian (endowed with the metric induced by the standard metric
on $\C^N$) so the previous lemma applies giving us a finite open
cover of $M$, $\{ U_\nu \}_\nu$, with the required property. Since
$\ppt$, $\pot$, and $\pmt$ can be defined to act on functions or
forms supported in $V_\mu$, they can certainly be defined to act
on functions or forms supported in $U_\nu$ for each $\nu$, so we
let $\pptg{\nu}$, $\potg{\nu}$, and $\pmtg{\nu}$ be the
pseudodifferential operators of order zero defined on $U_\nu$ and
$\cp_\nu$, $\co_\nu$, and $\cp_\nu$ be the three regions of the
$\xi$ space dual to $U_\nu$ on which the symbol of each of those three
pseudodifferential operators is supported. By construction, the
following is true: $$(\pptg {\nu})^* \pptg{\nu}  + (\potg {\nu})^*
\potg {\nu} + (\pmtg {\nu})^* \pmtg {\nu} = Id$$ There are only
two other conditions that need to be imposed on this covering,
both of which will become clear once we state and prove the
following lemma:
\medskip
\newtheorem{changeofvar}[lemma]{Lemma}
\begin{changeofvar} Let $U_\nu$ and $U_\mu$ be two neighborhoods
of a compact, orientable CR manifold $M$ \label{changevar} which
are defined as above and satisfy $U_\nu \cap \, U_\mu \, \neq \,
\emptyset$. If $U_\nu$ and $U_\mu$ are chosen small enough, then
there exists a diffeomorphism $\vartheta$ between $U_\nu$ and
$U_\mu$ with Jacobian $\mathcal{J}_\vartheta$ such that:
\begin{enumerate}
\item[(i)] ${\,}^t
\hspace{-0.25em} \mathcal{J}_\vartheta (\cp_\mu) \cap \cm_\nu \, =
\, \emptyset$ and $\cp_\nu \cap \, {\,}^t \hspace{-0.25em}
\mathcal{J}_\vartheta (\cm_\mu) \, = \, \emptyset$, where ${{\,}^t
\hspace{-0.25em} \mathcal{J}_\vartheta}$ is the inverse of the
transpose of the Jacobian of $\vartheta$;
\item[(ii)] Let ${\,}^\vartheta \pptg{\mu}$, ${\,}^\vartheta \pmtg{\mu}$,
and ${\,}^\vartheta \potg{\mu}$ be the transfers of $\pptg{\mu}$,
$\pmtg{\mu}$, and $\potg{\mu}$ respectively via $\vartheta$, then
on $\{ \xi \, | \, \xi_{2n-1} \, \geq \, \frac {4}{5} |\xi'| \:
and \: |\xi| \geq (1 + \eps)  \, tA \}$, the principal symbol of
${\,}^\vartheta \pptg{\mu}$ is identically equal to  $1$, on $\{
\xi \, | \, \xi_{2n-1} \, \leq \, - \frac {4}{5} |\xi'| \: and \:
|\xi| \geq (1+ \eps) \, tA \}$, the principal symbol of
${\,}^\vartheta \pmtg{\mu}$ is identically equal to $1$, and on
$\{ \xi \, | \, - \frac {1}{3} |\xi'| \, \leq \, \xi_{2n-1} \,
\leq \, \frac {1}{3} |\xi'| \: and \: |\xi| \geq (1 + \eps) \, tA
\}$, the principal symbol of ${\,}^\vartheta \potg{\mu}$ is
identically equal to $1$, where $\eps > 0$ and can be very small.
\end{enumerate}
\end{changeofvar}

\smallskip
\noindent {\bf{Proof:}} $U_\nu$ and $U_\mu$ are neighborhoods of
the manifold $M$ and if chosen to be very small, a diffeomorphism
$\vartheta$ from $U_\nu$ to $U_\mu$ will be a small perturbation
of a translation in their respective CR local coordinates.
Remember these local CR coordinates exist because they do on each
$V_\eta$ and if $U_\nu \subset V_\eta$, then $U_\mu \subset
V_\eta$ by construction since $U_\nu \cap \, U_\mu \, \neq \,
\emptyset$. A translation has the identity matrix as its Jacobian,
so $\mathcal{J}_\vartheta$ will be very close to the identity
matrix. In fact, the smaller the two neighborhoods $U_\nu$ and
$U_\mu$, the closer $\mathcal{J}_\vartheta$ will be to the
identity matrix. Let $(x, \xi)$ be the coordinates in $U_\nu$ and
its dual space. Also, let $(y, \rho)$ be the coordinates in
$U_\mu$ and its dual space. According to \cite{Treves}, the
diffeomorphism $\vartheta$ induces a linear map on the dual space,
so $$U_\nu \ni (x, \xi) \longrightarrow (y(x), {\,}^t
\hspace{-0.25em} \mathcal{J}_\vartheta^{-1}(x) \xi) = (y, \rho)
\in U_\mu.$$ Since $\cp_\mu$ and $\cm_\mu$ not only do not
intersect, but can be separated by open sets, and similarly
$\cp_\nu$ and $\cm_\nu$, and since ${\,}^t \hspace{-0.25em}
\mathcal{J}_\vartheta^{-1}(x)$ is close to the identity for all $x
\in U_\nu$, $\cp_\mu \cap {{\,}^t \hspace{-0.25em}
\mathcal{J}_\vartheta^{-1}(\cm_\nu)} \, = \, \emptyset$ and
${{\,}^t \hspace{-0.25em} \mathcal{J}_\vartheta^{-1}(\cp_\nu)}
\cap \, \cm_\mu \, = \, \emptyset$ follow for $U_\nu$ and $U_\mu$
small enough. Now, apply ${\,}^t \hspace{-0.25em}
\mathcal{J}_\vartheta$ on both sides of each of these expressions
to get (i). So we see that when doing the global
microlocalization, we do not have to worry about positive
truncated cones intersecting the negative ones after changes of
variables. Now, to show (ii), we need to make use again of
material in \cite{Treves}. By definition, the three operators
$\pptg{\mu}$, $\pmtg{\mu}$, and $\potg{\mu}$ have symbols
identically equal to $1$ on the sets $\{ \rho \, | \, \rho_{2n-1}
\, \geq \, \frac {3}{4} |\rho'| \: and \: |\rho| \geq tA \}$, $\{
\rho \, | \, \rho_{2n-1} \, \leq \, - \frac {3}{4} |\rho'| \: and
\: |\rho| \geq  tA \}$, and $\{ \rho \, | \, - \frac {1}{2}
|\rho'| \, \leq \, \rho_{2n-1} \, \leq \, \frac {1}{2} |\rho'| \:
and \: |\rho| \geq tA \}$ respectively. Since $1$ is homogeneous
of degree $0$ in $\rho$, their principal symbols, which we denote
by $\sigma^+_{t, \, \mu}$, $\sigma^-_{t, \, \mu}$, and
$\sigma^0_{t, \, \mu}$, must also be identically $1$ on these
sets. Now, let ${\,}^\vartheta \sigma^+_{t, \, \mu}$,
${\,}^\vartheta \sigma^-_{t, \, \mu}$, and ${\,}^\vartheta
\sigma^0_{t, \, \mu}$ be the principal symbols of ${\,}^\vartheta
\pptg{\mu}$, ${\,}^\vartheta \pmtg{\mu}$, and ${\,}^\vartheta
\potg{\mu}$ respectively, namely of the transfers via $\vartheta$
of the three operators defined on $U_\mu$. \cite{Treves} allows us
to conclude that ${\,}^\vartheta \sigma^+_{t, \, \mu} (x, \xi) \,
= \,  \sigma^+_{t, \, \mu} (y(x), {\,}^t \hspace{-0.25em}
\mathcal{J}_\vartheta^{-1}(x) \xi) \, = \, \sigma^+_{t, \, \mu}
(y, \rho)$, ${\,}^\vartheta \sigma^-_{t, \, \mu} (x, \xi) \, = \,
\sigma^-_{t, \, \mu} (y(x), {\,}^t \hspace{-0.25em}
\mathcal{J}_\vartheta^{-1}(x) \xi) \, = \, \sigma^-_{t, \, \mu}
(y, \rho)$, and ${\,}^\vartheta \sigma^0_{t, \, \mu} (x, \xi) \, =
\,  \sigma^0_{t, \, \mu} (y(x), {\,}^t \hspace{-0.25em}
\mathcal{J}_\vartheta^{-1}(x) \xi) \, = \, \sigma^0_{t, \, \mu}
(y, \rho)$. ${\,}^t \hspace{-0.25em} \mathcal{J}_\vartheta^{-1}$
is close to the identity for $U_\nu$ and $U_\mu$ small enough, so
it will send the sets specified in (ii) into $\{ \rho \, | \,
\rho_{2n-1} \, \geq \, \frac {3}{4} |\rho'| \: and \: |\rho| \geq
 tA \}$, $\{ \rho \, | \, \rho_{2n-1} \, \leq \, -
\frac {3}{4} |\rho'| \: and \: |\rho| \geq tA \}$, and $\{ \rho \,
| \, - \frac {1}{2} |\rho'| \, \leq \, \rho_{2n-1} \, \leq \,
\frac {1}{2} |\rho'| \: and \: |\rho| \geq  tA \}$ respectively,
sets on which $\sigma^+_{t, \, \mu}$, $\sigma^-_{t, \, \mu}$, and
$\sigma^0_{t, \, \mu}$ are identically equal to $1$. This proves
(ii). Notice that we get the principal symbols identically equal
to $1$ on truncated cones smaller than the ones for $\pptg{\nu}$,
$\pmtg{\nu}$, and $\potg{\nu}$, (namely differing by pieces of
annuli of radius $\eps$) which merely reflects the fact that
whereas ${\,}^t \hspace{-0.25em} \mathcal{J}_\vartheta^{-1}$ is
close to the identity matrix for $U_\nu$ and $U_\mu$ small, it is
not quite equal to it. \qed

\medskip
\medskip
Shrink $r$ a bit in the proof of the Lemma~\ref{coveringlemma}
such that any two sets in $\{ U_\nu \}_\nu$ with non-empty
intersection fulfill the previous lemma -- this is the first
condition promised before the lemma. Also, let $\pptdg{\mu}$ and
$\pmtdg{\mu}$ be pseudodifferential operators which dominate
$\pptg{\mu}$ and $\pmtg{\mu}$ respectively, where $\pptg{\mu}$ and
$\pmtg{\mu}$ are defined on some neighborhood $U_\mu$ in $\{ U_\nu
\}_\nu$. $\pptdg{\mu}$ and $\pmtdg{\mu}$ are defined to be inverse
zero order $t$ dependent like $\pptg{\mu}$ and $\pmtg{\mu}$. We
denote by $\cpd_\mu$ and $\cmd_\mu$ the supports of the symbols of
$\pptdg{\mu}$ and $\pmtdg{\mu}$, respectively. We shrink $r$ once
again (this is just the second time so there is no threat of it
going to zero) in the proof of the Lemma~\ref{coveringlemma} to
get $\{ U_\nu \}_\nu$ and define $\pptdg{\mu}$ and $\pmtdg{\mu}$
in such a way that for all $\nu$ such that $U_\nu \cap \, U_\mu \,
\neq \, \emptyset$, the following hold:
\begin{enumerate}
\item[(i)] Let $\vartheta$ be the diffeomorphism between $U_\nu$
and $U_\mu$ with Jacobian $\mathcal{J}_\vartheta$ given by the
previous lemma, then ${\,}^t \hspace{-0.25em}
\mathcal{J}_\vartheta (\cpd_\mu) \cap \cm_\nu \, = \, \emptyset$
and $\cp_\nu \cap \, {\,}^t \hspace{-0.25em} \mathcal{J}_\vartheta
(\cmd_\mu) \, = \, \emptyset$;
\item[(ii)] Let ${\,}^\vartheta \pptdg{\mu}$, ${\,}^\vartheta
\pmtdg{\mu}$ be the transfers via $\vartheta$ of $\pptdg{\mu}$ and
$\pmtdg{\mu}$, respectively. Then the principal symbol of
${\,}^\vartheta \pptdg{\mu}$ is identically $1$ on $\cp_\nu$ and
the principal symbol of ${\,}^\vartheta \pmtdg{\mu}$ is
identically $1$ on $\cm_\nu$;
\item[(iii)] $\cpd_\mu \cap \cmd_\mu \, = \, \emptyset$.
\end{enumerate}
This is the second condition on the covering promised before the
change of variables lemma, and we shall call \label{propdagger} it
Property $\dag$. Notice that (ii) is not quite (ii) of the
previous lemma, but another look at the proof of the lemma will
show that it is merely a slight variation which is proven the same
way. Notice also that if the CR manifold $M$ were not orientable,
neither Lemma~\ref{changevar}, nor Property $\dag$ would hold, so
orientability is a necessary assumption in
Proposition~\ref{global}.

\smallskip
\noindent {\bf Observation:} In order to save notation, the left
superscript $\vartheta$ indicating the transfer of a
pseudodifferential operator into another local coordinate system
will be henceforth supressed. It will be obvious from the context
which pseudodifferential operators need to be transferred in order
for certain compositions to make sense.

\smallskip
\noindent Let now $\{ \ze_\nu \}_\nu$ be a partition of unity
subordinate to the covering $\{ U_\nu \}_\nu$ satisfying $\sum_\nu
{\ze_\nu}^2 \, = \, 1$, and for each $\nu$ let $\dze_\nu$ be a
cutoff function that dominates $\ze_\nu$ such that $supp \,
(\dze_\nu) \subset U_\nu$. Then we define the global norm as
follows:
$${\langle| \varphi |\rangle}_t^2 =
\sum_\nu \: (|||\dze_\nu \pptg {\nu} \ze_\nu \varphi^\nu|||_t^2 +
||\dze_\nu \potg {\nu} \ze_\nu \varphi^\nu ||_0^2 + |||\dze_\nu
\pmtg {\nu} \ze_\nu \varphi^\nu|||_{-t}^2)$$ where $\varphi$ is a
$(0,1)$ form in $L_2^1 (M)$ and $\varphi^\nu$ is the form
expressed in the local coordinates on $U_\nu$. The reader may
refer back to page \pageref{triplenotation} for an explanation of
the norms appearing in this expression. Since $\pptg {\nu},$
$\potg {\nu}, $ and $\pmtg {\nu} $ actually operate on each
component of the form, the above expression means
$${\langle| \varphi |\rangle}_t^2 = \sum_\nu \sum_i \:
(|||\dze_\nu \pptg {\nu} \ze_\nu \varphi_i^\nu|||_t^2 + ||\dze_\nu
\potg {\nu} \ze_\nu \varphi_i^\nu ||_0^2 + |||\dze_\nu \pmtg {\nu}
\ze_\nu \varphi_i^\nu|||_{-t}^2),$$ where $\varphi \, = \, \sum_i
\varphi_i \, \omi{i},$ but to save notation we will use the first
expression instead of the second. We define the norm in the exact
same way for functions. This norm depends on the covering $\{
U_\nu \}_\nu$ and its subordinate partition of unity $\{ \ze_\nu
\}_\nu$, but that does not matter overall because as we will prove
shortly any such global norm is equivalent to the $L^2$ norm on
$M.$

\smallskip
\noindent Before doing so, however, let us finishing defining the
norms we need on $M:$

\smallskip
\newtheorem{Sobolevdef}[lemma]{Definition}
\begin{Sobolevdef} Let the Sobolev norm of order $s$ for a form $\varphi$ supported on $M$ be
given by: $$||\varphi||^2_s = \sum_\eta \sum_i || \dze_\eta \,
\Lambda^s \, \ze_\eta \, \varphi_i^\eta||^2_0,$$ where as usual
$\Lambda$ is defined to be the pseudodifferential operator with
symbol $(1+|\xi|^2)^\half.$ Again, to save notation we will write
just $$||\varphi||^2_s = \sum_\eta || \dze_\eta \, \Lambda^s \,
\ze_\eta \, \varphi^\eta||^2_0.$$
\end{Sobolevdef}
\noindent The definition for functions is similar. Notice that the
Sobolev norm can be defined in more than one way, but that any two
such definitions are equivalent.

\medskip
\noindent Let us now prove that the global norm is equivalent to
the $L^2$ norm on $M.$


\medskip
\newtheorem{L2normequiv}[lemma]{Lemma}
\begin{L2normequiv} For any $t$, \label{L2equiv} there exist two
positive constants depending on $t$, $C_t$ and $C'_t$ such that
$$ C_t ||\varphi||_0^2 \leq {\langle| \varphi |\rangle}_t^2 \leq
C'_t  ||\varphi||_0^2,$$ where $\varphi$ is a form in $L_2^1 (M)$.
\end{L2normequiv}
\noindent {\bf{Proof:}} $${\langle| \varphi |\rangle}_t^2 =
\sum_\nu \: (|||\dze_\nu \pptg {\nu} \ze_\nu \varphi^\nu|||_t^2 +
||\dze_\nu \potg {\nu} \ze_\nu \varphi^\nu ||_0^2 + |||\dze_\nu
\pmtg {\nu} \ze_\nu \varphi^\nu|||_{-t}^2)$$ Since $M$ is compact,
the two weight functions $e^{-t\lambda}$ and $e^{t\lambda}$ have
maximum and minimum values on it, so for each $t$, the two norms
$||| \cdot |||_t$ and $||| \cdot |||_{-t}$ are equivalent to the
$L^2$ norm. Then there exist two constants $C_{1, \, t}$ and
$C_{2, \, t}$ such that
\begin{equation*}
\begin{split}
&\sum_\nu \: (|||\dze_\nu \pptg {\nu} \ze_\nu \varphi^\nu|||_t^2 +
||\dze_\nu \potg {\nu} \ze_\nu \varphi^\nu ||_0^2 + |||\dze_\nu
\pmtg {\nu} \ze_\nu \varphi^\nu|||_{-t}^2)\\ &\leq \sum_\nu \:
(C_{1, \, t}||\dze_\nu \pptg {\nu} \ze_\nu \varphi^\nu||_0^2 +
||\dze_\nu \potg {\nu} \ze_\nu \varphi^\nu ||_0^2+ C_{2, \,
t}||\dze_\nu \pmtg {\nu} \ze_\nu \varphi^\nu||_0^2)
\end{split}
\end{equation*}
Take the term $||\dze_\nu \pptg {\nu} \ze_\nu \varphi^\nu||_0^2$
and express it as
$$||\dze_\nu \pptg {\nu} \ze_\nu \varphi^\nu||_0^2 = ||(1 - (1 -\dze_\nu)) \pptg {\nu} \ze_\nu
\varphi^\nu||_0^2 \leq 2 ||\pptg {\nu} \ze_\nu \varphi^\nu||_0^2 +
2||(1-\dze_\nu) \pptg {\nu} \ze_\nu \varphi^\nu||_0^2$$
Technically $\pptg {\nu} \ze_\nu \varphi^\nu$ does not make sense
because although $\ze_\nu \varphi^\nu$ is supported in $U_\nu$,
$\pptg {\nu} \ze_\nu \varphi^\nu$ spreads out in the whole $2n-1$
dimensional real space used to parametrize $U_\nu$. One gets
around this difficulty by noticing that outside the support of
$\ze_\nu$, $\pptg{\nu}$ is a smoothing operator, hence the outside
which is given by $||(1-\dze_\nu) \pptg {\nu} \ze_\nu
\varphi^\nu||_0^2$ can be controlled by any negative Sobolev norm,
in particular by $||\ze_\nu \varphi^\nu||_{-1}^2$. Therefore,
$$||\dze_\nu \pptg {\nu} \ze_\nu \varphi^\nu||_0^2 \leq 2 ||\pptg {\nu} \ze_\nu
\varphi^\nu||_0^2+ 2 C_{-1}||\ze_\nu \varphi^\nu||_{-1}^2,$$ for
some constant $C_{-1}>0$. Next, apply Plancherel's Theorem to see
that $$|| \pptg {\nu} \ze_\nu \varphi^\nu||_0^2 \, = \, ||
\psi^+_{\nu, \, t} \widehat{\ze_\nu \varphi^\nu}||_0^2 \, \leq \,
||\widehat{\ze_\nu \varphi^\nu}||_0^2 \, = \, ||\ze_\nu
\varphi^\nu||_0^2,$$ because $\psi^+_{\nu, \, t}$ is the symbol of
$\pptg {\nu}$ and takes values in $[0,1]$ by construction. We
conclude that there exists a constant $C_{\nu, \, +}$ such that
$||\dze_\nu \pptg {\nu} \ze_\nu \varphi^\nu||_0^2  \, \leq \,
C_{\nu, \, +}||\ze_\nu \varphi^\nu||_0^2$. Similarly, there exist
constants $C_{\nu, \, 0}$ and $C_{\nu, \, -}$ such that
$||\dze_\nu \potg {\nu} \ze_\nu \varphi^\nu||_0^2  \, \leq \,
C_{\nu, \, 0}||\ze_\nu \varphi^\nu||_0^2$ and $||\dze_\nu \pmtg
{\nu} \ze_\nu \varphi^\nu||_0^2 \, \leq \, C_{\nu, \, -}|| \ze_\nu
\varphi^\nu||_0^2$ which implies that
\begin{equation*}
\begin{split}
{\langle| \varphi |\rangle}_t^2 &\leq \sum_\nu \: (C_{1, \,
t}||\dze_\nu \pptg {\nu} \ze_\nu \varphi^\nu||_0^2 + ||\dze_\nu
\potg {\nu} \ze_\nu \varphi^\nu ||_0^2+ C_{2, \, t}||\dze_\nu
\pmtg {\nu} \ze_\nu \varphi^\nu||_0^2)\\ &\leq \sum_\nu \: (C_{1,
\, t} C_{\nu, \, +}|| \ze_\nu \varphi^\nu||_0^2 + C_{\nu, \, 0} ||
\ze_\nu
\varphi^\nu ||_0^2+C_{2, \, t} C_{\nu, \, -}||\ze_\nu \varphi^\nu||_0^2)\\
&\leq C_{3, \, t} \sum_\nu || \ze_\nu \varphi^\nu ||_0^2,
\end{split}
\end{equation*}
where $C_{3, \, t} \, = \, \max \, \{C_{\nu, \, 0}, C_{1, \, t}
C_{\nu, \, +}, C_{2, \, t} C_{\nu, \, -} \}$. $\{ \ze_\nu \}_\nu$
is a partition of unity subordinate to the finite covering $\{
U_\nu \}_\nu$ of $M$, so the number of functions $\ze_\nu$ is
finite, and each $\ze_\nu$ takes values in $[0,1]$, which means
$|| \ze_\nu \varphi^\nu ||_0^2 \, \leq \, ||\varphi||_0^2$. Thus,
there is exists a constant $C'_t$ such that
$${\langle| \varphi |\rangle}_t^2 \leq C_{3, \, t} \sum_\nu || \ze_\nu \varphi^\nu
||_0^2 \leq C'_t ||\varphi||_0^2.$$

\noindent To prove $||\varphi||^2_0$ bounds ${\langle| \varphi
|\rangle}_t^2$ from below, we use the partition of unity $\{
\ze_\nu \}_\nu$ and the corresponding set of functions dominating
each function in the partition, namely $\{ \dze_\nu \}_\nu$. Given
$\sum_\nu \ze_\nu^2 \, = \, 1 \, = \, \sum_\nu \dze_\nu \ze_\nu^2
,$ we have
\begin{equation*}
\begin{split}
||\varphi||^2_0  &=  (\sum_\nu  \ze_\nu^2 \varphi, \varphi)_0 =
\sum_\nu ||\ze_\nu \varphi^\nu||^2_0  \\ &= \sum_\nu ( ((\pptg
{\nu})^* \pptg{\nu}  + (\potg {\nu})^* \potg {\nu} + (\pmtg
{\nu})^* \pmtg {\nu}) \ze_\nu \varphi^\nu, \ze_\nu \varphi^\nu)_0
\\ &= \sum_\nu \: (||(\dze_\nu + (1-\dze_\nu)) \pptg{\nu} \ze_\nu
\varphi^\nu||^2_0 + ||(\dze_\nu + (1-\dze_\nu)) \potg{\nu} \ze_\nu
\varphi^\nu||^2_0
 \\ & \quad + ||(\dze_\nu+ (1-\dze_\nu)) \pmtg{\nu} \ze_\nu
\varphi^\nu||^2_0),
\end{split}
\end{equation*}
since by construction $(\pptg {\nu})^* \pptg{\nu}  + (\potg
{\nu})^* \potg {\nu} + (\pmtg {\nu})^* \pmtg {\nu} \, = \, Id,$
where $Id$ is the identity operator. As above, $(1-\dze_\nu)
\pptg{\nu} \ze_\nu,$ $(1-\dze_\nu) \potg{\nu} \ze_\nu,$ and
$(1-\dze_\nu) \pmtg{\nu} \ze_\nu$ are smoothing pseudodifferential
operators which implies $|| (1-\dze_\nu) \pptg{\nu} \ze_\nu
\varphi^\nu||^2_0,$ $|| (1-\dze_\nu) \potg{\nu} \ze_\nu
\varphi^\nu||^2_0,$ and $|| (1-\dze_\nu) \pmtg{\nu} \ze_\nu
\varphi^\nu||^2_0$ can be controlled by $||\dze_\nu \pptg{\nu}
\ze_\nu \varphi^\nu||^2_0,$ $|| \dze_\nu \potg{\nu} \ze_\nu
\varphi^\nu||^2_0,$ and $|| \dze_\nu \pmtg{\nu} \ze_\nu
\varphi^\nu||^2_0$ respectively. This means the previous
expression can be rewritten as follows:
\begin{equation*}
\begin{split}
||\varphi||^2_0  &\leq 4 \sum_\nu \: (||\dze_\nu \pptg{\nu}
\ze_\nu \varphi^\nu||^2_0 + ||\dze_\nu  \potg{\nu}  \ze_\nu
\varphi^\nu||^2_0 + ||\dze_\nu \pmtg{\nu} \ze_\nu
\varphi^\nu||^2_0).
\end{split}
\end{equation*}
The two norms $||| \cdot |||_{\, t}$ and $||| \cdot |||_{\, -t}$
are equivalent to the $L^2$ norm, so there exist two constants
$C_{4, \, t}$ and $C_{5, \, t}$ for which
\begin{equation*}
\begin{split}
||\varphi||^2_{\, 0}  &\leq 4 \sum_\nu \: (C_{4, \, t}|||\dze_\nu
\pptg{\nu} \ze_\nu \varphi^\nu|||^2_t + ||\dze_\nu  \potg{\nu}
\ze_\nu \varphi^\nu||^2_0 + C_{5, \, t} |||\dze_\nu \pmtg{\nu}
\ze_\nu \varphi^\nu|||^2_{-t})
\\ &\leq C_t \sum_\nu \: (|||\dze_\nu \pptg {\nu} \ze_\nu \varphi^\nu|||_t^2 +
||\dze_\nu \potg {\nu} \ze_\nu \varphi^\nu ||_0^2 + |||\dze_\nu
\pmtg {\nu} \ze_\nu \varphi^\nu|||_{-t}^2) \\ &= C_t \, {\langle|
\varphi |\rangle}_t^2,
\end{split}
\end{equation*}
where $C_t \, = \, \max \, \{4, 4 \, C_{4, \, t}, 4 \, C_{5, \, t}
\}$. \qed

\medskip
\noindent We thus know the global norm is indeed equivalent to the
$L^2$ norm. This result has the following corollary:

\smallskip
\newtheorem{equivoperator}[lemma]{Corollary}
\begin{equivoperator}
There exists a self-adjoint operator $G_t$ \label{equivop} such
that $$ (\varphi, \phi)_0 = {\langle| \varphi, G_t \, \phi
|\rangle}_t,$$ for any two forms $\varphi$ and $\phi$ in $L_2^1
(M)$. $G_t$ is the inverse of $$\sum_\nu \: \left(\ze_\nu (\pptg
{\nu})^* \dze_\nu e^{-t \lambda} \dze_\nu \pptg {\nu} \ze_\nu +
\ze_\nu (\potg {\nu})^* \dze_\nu^2 \potg {\nu} \ze_\nu + \ze_\nu
(\pmtg {\nu})^* \dze_\nu e^{t \lambda} \dze_\nu \pmtg {\nu}
\ze_\nu\right).$$
\end{equivoperator}
\noindent {\bf{Proof:}} Use Lemma~\ref{L2equiv} and the Riesz
representation theorem to see that such a self-adjoint operator
$G_t$ must exist. Then notice that $${\langle| \varphi, \phi
|\rangle}_t = (\varphi, \sum_\nu \: (\ze_\nu (\pptg {\nu})^*
\dze_\nu e^{-t \lambda} \dze_\nu \pptg {\nu} \ze_\nu + \ze_\nu
(\potg {\nu})^* \dze_\nu^2 \potg {\nu} \ze_\nu + \ze_\nu (\pmtg
{\nu})^* \dze_\nu e^{t \lambda} \dze_\nu \pmtg {\nu} \ze_\nu) \,
\phi )_0, $$ which implies that $G_t$ and $$F_t = \sum_\nu \:
\left(\ze_\nu (\pptg {\nu})^* \dze_\nu e^{-t \lambda} \dze_\nu
\pptg {\nu} \ze_\nu + \ze_\nu (\potg {\nu})^* \dze_\nu^2 \potg
{\nu} \ze_\nu + \ze_\nu (\pmtg {\nu})^* \dze_\nu e^{t \lambda}
\dze_\nu \pmtg {\nu} \ze_\nu\right)$$ are inverses. Notice also
that the symbol of $F_t$ is never zero, so an inverse does exist.
\qed

\medskip \noindent For the sake of completeness, we include as a
lemma a very elementary fact which will be used repeatedly:
\smallskip
\newtheorem{sclclemma}[lemma]{Lemma}
\begin{sclclemma}
Let $( \ \cdot \ , \ \cdot \ )$ be any Hermitian inner product.
For any two functions $\alpha$ and $\beta$ \label{sclc} and any
small positive number $0 < \epsilon \ll 1$, the following holds
$$- \epsilon {||\alpha||}^2 - \frac {1}{\epsilon} \, {||\beta||}^2 \leq 2 \Re \{ (\alpha,\beta) \}
\leq \epsilon {||\alpha||}^2 +\frac {1}{\epsilon}
\,{||\beta||}^2,$$ where $|| \cdot ||$ is the norm corresponding
to the given inner product.
\end{sclclemma}
\noindent {\bf{Proof:}} Expand ${||\sqrt{\epsilon} \, \alpha -
\sqrt{\frac {1}{\epsilon}} \, \beta||}^2 \geq 0$ to obtain that $2
\Re \{ (\alpha,\beta) \} \leq \epsilon {||\alpha||}^2 + \frac
{1}{\epsilon} \, {||\beta||}^2$. For the rest of the inequality,
expand $-{||\sqrt{\epsilon} \, \alpha + \sqrt{\frac {1}{\epsilon}}
\, \beta||}^2 \leq 0$. \qed

\bigskip

\noindent Now we need to look at the adjoint of $\dbarb$
corresponding to the global norm defined above. Notice that since
such an adjoint exists for the $L^2$ norm, it follows from
Lemma~\ref{L2equiv} and the Riesz representation theorem that an
adjoint also exists for this global norm.

\medskip
\newtheorem{adjointclaim}[lemma]{Lemma}
\begin{adjointclaim}
Let $f$ be a smooth function and $\varphi$ a smooth $(0,1)$ form
on $M$.\label{adclaim}
\begin{equation*}
\begin{split}
\adt &= \ad -t\sum_\rho \ze^2_\rho \pptdg{\rho} [\ad , \lambda] +
t \sum_\rho \ze^2_\rho \pmtdg{\rho} [\ad, \lambda]\\ &\quad
+\sum_\rho \: (\dze_\rho [\dze_\rho \pptg {\rho} \ze_\rho,
\dbarb]^* \dze_\rho \pptg {\rho} \ze_\rho + \ze_\rho (\pptg
{\rho})^* \dze_\rho [\ad -t [\ad, \lambda], \dze_\rho \pptg {\rho}
\ze_\rho] \dze_\rho+ \\
&\quad+\ze_\rho (\pmtg {\rho})^* \dze_\rho [\ad +t [\ad, \lambda],
\dze_\rho \pmtg {\rho} \ze_\rho] \dze_\rho +\dze_\rho [\dze_\rho
\pmtg {\rho} \ze_\rho, \dbarb]^* \dze_\rho \pmtg {\rho} \ze_\rho+
E_t
\end{split}
\end{equation*}
where the error given by the psedodifferential operator $E_t$ is
of order zero and up to terms of one order lower, its symbol is
supported in $\co_\rho$ for each $\rho$.
\end{adjointclaim}

\smallskip
\noindent {\bf Proof of Lemma ~\ref{adclaim}:} This is a
straightforward and rather tedious computation, so the details
have been relegated to Section~\ref{Compdetails}. \qed

\medskip
\noindent Finally, we can compute the global energy form $\qbtg
{\varphi}\, = \, \sqtnorm{\dbarb \varphi} + \sqtnorm{\adt
\varphi}.$

\medskip
\newtheorem{globalenergycalc}[lemma]{Lemma}
\begin{globalenergycalc}
There exist constants $K > 1$ and $K_t >0$ such that for each
smooth $(0,1)$ form $\varphi$ on $M$ \label{globalenergy}
\begin{equation*}
\begin{split}
&K \qbtg {\varphi}+K_t \sum_\nu \:{||\dze_\nu \potdg {\nu} \ze_\nu
\varphi^\nu||}^2_0+O({\langle |\varphi|\rangle}^2_t) +O_t
(||\varphi||^2_{-1})\\&\geq \sum_\nu \: {|||\dze_\nu \pptg {\nu}
\adp\ze_\nu \varphi^\nu|||}^2_{\,t}+ \sum_\nu \:{||\dze_\nu \potg
{\nu} \ze_\nu \ado \varphi^\nu||}^2_0+ \sum_\nu \:{|||\dze_\nu
\pmtg {\nu} \adm\ze_\nu \varphi^\nu|||}^2_{-t}\\&\quad + \sum_\nu
\: {|||\dze_\nu \pptg {\nu} \ze_\nu \dbarb
\varphi^\nu|||}^2_{\,t}+ \sum_\nu \:{||\dze_\nu \potg {\nu}
\ze_\nu \dbarb \varphi^\nu||}^2_0+ \sum_\nu \:{|||\dze_\nu \pmtg
{\nu} \ze_\nu \dbarb \varphi^\nu|||}^2_{-t}.
\end{split}
\end{equation*}
\end{globalenergycalc}
\smallskip \noindent {\bf Proof:} Once again, the details can be
found in Section~\ref{Compdetails} for the same reason as above.
\qed

\bigskip\noindent We will be working in each neighborhood $U_\nu$
separately, so it is necessary to pull $\dbarb$ and the three
adjoints $\adp,$ $\ado,$ and $\adm$ outside the cutoff functions
and the pseudodifferential operators $\pptg{\nu},$ $\potg{\nu},$
and $\pmtg{\nu}$ in the previous expression. Notice that the
commutators $[\dze_\nu \pptg{\nu} \ze_\nu, \dbarb],$ $[\dze_\nu
\potg{\nu} \ze_\nu, \dbarb],$ $[\dze_\nu \pmtg{\nu} \ze_\nu,
\dbarb],$ $[\dze_\nu \pptg{\nu} , \adp],$ $[\dze_\nu \potg{\nu}
\ze_\nu, \ado],$ and $[\dze_\nu \pmtg{\nu} , \adm],$ are of order
zero and bounded independently of $t$ since $\pptg{\nu},$
$\potg{\nu},$ and $\pmtg{\nu}$ are inverse zero order $t$
dependent and the cutoff functions are independent of $t,$ so the
errors can all be absorbed in $O({\langle |\varphi|\rangle}^2_t)$.
Thus,
\begin{equation*}
\begin{split}
\sum_\nu \: {|||\dze_\nu \pptg {\nu}\ze_\nu \dbarb
\varphi|||}^2_{\,t} &= \sum_\nu \: {|||\dbarb \dze_\nu \pptg
{\nu} \ze_\nu \varphi+ [\dze_\nu \pptg{\nu} \ze_\nu, \dbarb]\dze_\nu \varphi|||}^2_{\,t}\\
&\geq (1-\eps) \, \sum_\nu \: {|||\dbarb \dze_\nu \pptg {\nu}
\ze_\nu \varphi|||}^2_{\,t}+O({\langle |\varphi|\rangle}^2_t),
\end{split}
\end{equation*}
by the $\epsilon$/$\frac {1}{\epsilon}$ trick from
Lemma~\ref{sclc}. The other five sums on the right-hand side of
the expression for $\qbtg{\varphi}$ are handled in the same way.
Therefore, we enlarge $K$ and $K_t$ slightly again and obtain the
following expression for $\qbtg{\varphi},$
\begin{equation}
\begin{split}
&K \qbtg {\varphi}+K_t \sum_\nu \:{||\dze_\nu \potdg {\nu} \ze_\nu
\varphi^\nu||}^2_0+O({\langle |\varphi|\rangle}^2_t) +O_t
(||\varphi||^2_{-1})\\&\geq \sum_\nu \: {|||\adp\dze_\nu \pptg
{\nu} \ze_\nu \varphi^\nu|||}^2_{\,t}+ \sum_\nu \:{||\ado\dze_\nu
\potg {\nu} \ze_\nu  \varphi^\nu||}^2_0+ \sum_\nu \:{|||\adm
\dze_\nu \pmtg {\nu} \ze_\nu \varphi^\nu|||}^2_{-t}\\&\quad +
\sum_\nu \: {|||\dbarb\dze_\nu \pptg {\nu} \ze_\nu
\varphi^\nu|||}^2_{\,t}+ \sum_\nu \:{||\dbarb\dze_\nu \potg {\nu}
\ze_\nu \varphi^\nu||}^2_0+ \sum_\nu \:{|||\dbarb\dze_\nu \pmtg
{\nu} \ze_\nu \varphi^\nu|||}^2_{-t}. \label{qbtexpr}
\end{split}
\end{equation}

\medskip
\noindent Before tackling the main estimate, we first prove some
local ones, namely for a $(0,1)$ form $\varphi$ supported in a
neighborhood $U$ of $M$. We choose an orthonormal basis of $(1,0)$
vector fields $L_1, \dots, L_{n-1}$ on $T^{1,0} (U)$, which means
the dual basis of $(1,0)$ forms $\omega_1, \dots, \omega_{n-1}$
will also be orthonormal. Since the local results will be applied
in each of the neighborhoods $U_\nu$, it is necessary to note that
by construction, each $U_\nu$ has an orthonormal basis of $(1,0)$
vector fields $L_1, \dots, L_{n-1}$ which satisfies
Equation~\ref{bracket}, namely there exist $\smooth$ functions
$a_{ij}^k,$ $b_{ij}^k,$ and $g_{ij}$ independent of $t$ such that
\begin{equation}
\begin{split}
[L_i,\lbar{j}] &= c_{ij} T + \sum_k a_{ij}^k \lbar{k} + \sum_k
b_{ij}^k \lbar{k}^{*, \, \pm t}+ g_{ij} + t \, e_{ij},
\end{split}
\end{equation}
where $|e_{ij}|$ is bounded independently of $t$ by a small
positive global constant, $\eps_G.$ Henceforth, we will assume the
neighborhood $U$ has this property, too.

Let $\ze$ be a cutoff function with support in $U$ and $\dze$ a
cutoff function which dominates it, whose support is contained in
a slightly larger open set $U'$. $U'$ is taken to be small enough
so that we can define on it the three pseudodifferential operators
of order zero $\ppt$, $\pot$, and $\pmt$ supported in $\cp$,
$\co$, and $\cm$ respectively, as detailed at the beginning of
this section. For each of the three norms, the $t$ norm, the $-t$
norm, and the zero norm, we need to define equivalents of
$\qbtg{\cdot}$ as follows: $$\qbtp {\phi} = {|||\dbarb
\phi|||}_{\, t}^2 + {|||\adp \phi|||}_{\, t}^2,$$ $$\qbto {\phi} =
{||\dbarb \phi||}_{\, 0}^2 + {||\ado \phi||}_{\, 0}^2,$$ and
$$\qbtm {\phi} = {|||\dbarb \phi|||}_{-t}^2 + {|||\adm
\phi|||}_{-t}^2.$$ Given this notation, we can rewrite
Equation~\ref{qbtexpr} as
\begin{equation}
\begin{split}
&K \qbtg {\varphi}+K_t \sum_\nu \:{||\dze_\nu \potdg {\nu} \ze_\nu
\varphi^\nu||}^2_0+O({\langle |\varphi|\rangle}^2_t) +O_t
(||\varphi||^2_{-1})\\&\geq \sum_\nu \: \qbtp{\dze_\nu \pptg {\nu}
\ze_\nu \varphi^\nu}+ \sum_\nu \:\qbto{\dze_\nu \potg {\nu}
\ze_\nu  \varphi^\nu}\\&\quad + \sum_\nu \: \qbtm{\dze_\nu \pmtg
{\nu} \ze_\nu \varphi^\nu}. \label{qbtexpression}
\end{split}
\end{equation}

\medskip
\noindent Next, we find expressions for $\qbtp{\cdot}$ and
$\qbtm{\cdot}$ in a series of two lemmas.
\smallskip
\newtheorem{qbtpcalc}[lemma]{Lemma}
\begin{qbtpcalc}
Let $\varphi$ be a $(0,1)$ form \label{qbtp} supported in $U',$
$\varphi \in Dom(\dbarb) \cap Dom(\ad).$ If $[L_i, \lbar{j}]=
c_{ij}T + \sum_k a^k_{ij}\lbar{k} + \sum_k b^k_{ij} \lbar{k}^{*,
\, t} + g_{ij} +t \, e_{ij},$ where $a^k_{ij}$, $b^k_{ij}$,
$g_{ij},$ and $e_{ij}$ are \smooth functions independent of $t$
and $c_{ij}$ are the coefficients of the Levi form, then there
exists $1 \gg \eps' > 0$ such that
\begin{equation*}
\begin{split}
\qbtp{\varphi} &\geq (1-\eps')\, \sum_{ij} |||\lbar{j}
\varphi_i|||^2_{\, t} + \sum_{ij} \Re \{(c_{ij} T \varphi_i,
\varphi_j)_t \}+ t \sum_{ij}\Re \{ ( \lbar{j} L_i (\lambda)
\varphi_i, \varphi_j)_t \}\\&\quad+ t \sum_{ij}\Re \{ ( e_{ij}
\varphi_i, \varphi_j)_t \} + O(|||\varphi|||^2_{\, t} ).
\end{split}
\end{equation*}
\end{qbtpcalc}
\noindent {\bf{Proof:}} Since $\varphi \in Dom(\ad),$ clearly
$\varphi \in Dom(\adp).$ Given $\varphi \, = \, \sum_i \varphi_i \
\omi{i}$, we have that $\adp \varphi = \sum_i \lbar{i}^{*, \,t}
(\varphi_i)$, therefore $(\adp \varphi, \adp \varphi)_t \, = \,
|||\adp \varphi|||^2_{\, t} \, = \, \sum_{i,j} (\lbar{i}^{*, \, t}
\varphi_i, \lbar{j}^{*, \, t} \varphi_j)_t$. Now, since
$$ \dbarb \  \varphi = \sum_i \dbarb (\varphi_i \wedge \omi{i})
=  \sum_{i<j} (\lbar {i} \varphi_j - \lbar {j} \varphi_i + \sum_k
m_{ij}^k \varphi_k)  \ \omi {i} \wedge \omi {j} \ ,$$ where
$m_{ij}^k$ are in $\smooth (U)$, the square of the norm of $\dbarb
\varphi$ equals the following: $$|||\dbarb \varphi|||^2_{\, t}  =
\sum_{i<j} |||\lbar{i} \varphi_j - \lbar{j} \varphi_i|||^2_{\, t}
+ O(|||\varphi|||^2_{\, t} + (\sum_{i,j} |||\lbar{i}
\varphi_j|||^2_{\, t} )^{\frac{1}{2}}|||\varphi|||_{\, t} )$$ Next
we calculate,
\begin{equation*}
\begin{split}
\sum_{i<j} |||\lbar{i} \varphi_j - \lbar{j} \varphi_i|||^2_{\, t}
&= \sum_{i \neq j} |||\lbar{j} \varphi_i|||^2_{\, t}  - \sum_{i
\neq j}
(\lbar{j} \varphi_i, \lbar{i} \varphi_j)_t \\
&= \sum_{i, j} |||\lbar{j} \varphi_i|||^2_{\, t}  - \sum_{i, j}
(\lbar{j} \varphi_i, \lbar{i} \varphi_j)_t \\ &= \sum_{i, j}
|||\lbar{j} \varphi_i|||^2_{\, t}  - \sum_{i, j} ( \lbar{i}^{*, \,
t} \lbar{j} \varphi_i, \varphi_j)_t \\ &= \sum_{i, j} |||\lbar{j}
\varphi_i|||^2_{\, t}  - \sum_{i, j} ( \lbar{i}^{*, \, t}
\varphi_i, \lbar{j}^{*, \, t} \varphi_j)_t - \sum_{i, j} (
[\lbar{i}^{*, \, t}, \lbar{j}] \varphi_i, \varphi_j)_t
\end{split}
\end{equation*}
We now combine the results gotten this far,
\begin{equation*}
\begin{split}
\qbtp{\varphi} &=|||\dbarb \varphi|||^2_{\, t}+|||\adp
\varphi|||^2_{\, t} \\
&= \sum_{i, j} |||\lbar{j} \varphi_i|||^2_{\, t} - \sum_{i, j} (
[\lbar{i}^{*, \, t}, \lbar{j}] \varphi_i, \varphi_j)_t - \sum_{i,
j} ( \lbar{i}^{*, \, t} \varphi_i, \lbar{j}^{*, \, t}
\varphi_j)_t+\sum_{i, j} ( \lbar{i}^{*, \, t} \varphi_i,
\lbar{j}^{*, \, t} \varphi_j)_t\\ & \quad+
O\bigg(|||\varphi|||^2_{\, t} + (\sum_{i,j} |||\lbar{i}
\varphi_j|||^2_{\, t} )^{\frac{1}{2}}|||\varphi|||_{\, t} \bigg)\\
&=\sum_{i, j} |||\lbar{j} \varphi_i|||^2_{\, t}  - \sum_{i, j} (
[\lbar{i}^{*, \, t}, \lbar{j}] \varphi_i, \varphi_j)_t +
O\bigg(|||\varphi|||^2_{\, t} + (\sum_{i,j} |||\lbar{i}
\varphi_j|||^2_{\, t} )^{\frac{1}{2}}|||\varphi|||_{\, t} \bigg)
\end{split}
\end{equation*}
Let $\eps$ be a very small positive number, then $(\sum_{i,j}
|||\lbar{i} \varphi_j|||^2_{\, t} )^{\frac{1}{2}}|||\varphi|||_{\,
t} \geq - \eps \sum_{i,j}|||\lbar{i} \varphi_j|||^2_{\, t} - \frac
{1}{\eps}|||\varphi|||^2_{\, t} \, ,$ which yields the following:
\begin{equation*}
\begin{split}
\qbtp{\varphi} &\geq (1 - \eps) \, \sum_{i, j} |||\lbar{j}
\varphi_i|||^2_{\, t} - \sum_{i, j} ( [\lbar{i}^{*, \, t},
\lbar{j}] \varphi_i, \varphi_j)_t + O(|||\varphi|||^2_{\, t})
\end{split}
\end{equation*}
Since $\lbar {i}^{*, \, t} \,  = \,  -  L_i  - f_i + t L_i
(\lambda)$, $[\lbar{i}^{*, \, t}, \lbar{j}] \, = \, [-  L_i  - f_i
+ t L_i (\lambda), \lbar{j}] \, = \, - [L_i, \lbar{j}] + \lbar{j}
(f_i) - t \lbar{j} L_i (\lambda)$. Moreover, once we split up
$[\lbar{i}^{*, \, t}, \lbar{j}]$, the terms which result are no
longer real-valued, so we have to take their real parts. Using
this expression for $[\lbar{i}^{*, \, t}, \lbar{j}]$ and the
hypothesis of the lemma about $[L_i, \lbar{j}]$ we thus obtain
\begin{equation*}
\begin{split}
\qbtp{\varphi} &\geq (1 - \eps) \,\sum_{i, j} |||\lbar{j}
\varphi_i|||^2_{\, t} - \sum_{i, j} ( (- [L_i, \lbar{j}] +
\lbar{j} (f_i) - t \lbar{j} L_i (\lambda)) \varphi_i, \varphi_j)_t
+ O(|||\varphi|||^2_{\, t}) \\ &= (1 - \eps) \,\sum_{i, j}
|||\lbar{j} \varphi_i|||^2_{\, t} + \sum_{i, j} \Re \{ ([L_i,
\lbar{j}] \varphi_i, \varphi_j)_t \} - \sum_{i, j}\Re \{ (
\lbar{j} (f_i)\varphi_i, \varphi_j)_t \}  \\ &  \quad+ t \sum_{i,
j}\Re \{( \lbar{j} L_i (\lambda) \varphi_i, \varphi_j)_t\}+
O(|||\varphi|||^2_{\, t})\\ &= (1 - \eps) \,\sum_{i, j}
|||\lbar{j} \varphi_i|||^2_{\, t} + \sum_{i, j} \Re \{(c_{ij}T
\varphi_i, \varphi_j)_t \} + \sum_{i, j} \Re \{ ( \sum_k
a^k_{ij}\lbar{k} \varphi_i, \varphi_j)_t \} \\& \quad+ \sum_{i, j}
\Re \{(\sum_k b^k_{ij} \lbar{k}^{*, \, t} \varphi_i,
\varphi_j)_t\}+ \sum_{i, j} \Re \{(g_{ij} \varphi_i, \varphi_j)_t
\} + t \sum_{ij}\Re \{ ( e_{ij} \varphi_i, \varphi_j)_t \} \\&
\quad- \sum_{i, j} \Re \{( \lbar{j} (f_i) \varphi_i,
\varphi_j)_t\}+ t \sum_{i, j}\Re \{(\lbar{j} L_i (\lambda)
\varphi_i, \varphi_j)_t \} + O(|||\varphi|||^2_{\, t})
\\&= (1 - \eps) \,\sum_{i, j} |||\lbar{j} \varphi_i|||^2_{\, t} + \sum_{i, j} \Re
\{ (c_{ij}T \varphi_i, \varphi_j)_t \}+ t \sum_{i, j}\Re \{
(\lbar{j} L_i (\lambda)  \varphi_i, \varphi_j)_t\}\\& \quad+ t
\sum_{ij}\Re \{ ( e_{ij} \varphi_i, \varphi_j)_t \}+ \sum_{i, j,
k} \Re \{(\lbar{k} \varphi_i, \bar a^k_{ij} \varphi_j)_t \}+
\sum_{i, j, k} \Re \{(\varphi_i, \lbar{k}(\bar b^k_{ij}
\varphi_j))_t \}\\ & \quad+ \sum_{i, j} \Re \{(g_{ij} \varphi_i,
\varphi_j)_t \}- \sum_{i, j} \Re \{( \lbar{j} (f_i) \varphi_i,
\varphi_j)_t \} + O(|||\varphi|||^2_{\, t})
\end{split}
\end{equation*}
Notice that $\sum_{i, j} Re \{(g_{ij} \varphi_i, \varphi_j)_t \}
\, = \, O(|||\varphi|||^2_{\, t}),$ $\sum_{i, j} Re \{( \lbar{j}
(f_i) \varphi_i, \varphi_j)_t \} \, = \, O(|||\varphi|||^2_{\,
t}),$ and for the same very small positive number $\eps,$ the
following two inequalities hold:
$$\sum_{i, j, k} Re \{(\lbar{k} \varphi_i, \bar a^k_{ij}
\varphi_j)_t \} \geq - \eps \sum_{i, k} |||\lbar{k}
\varphi_i|||^2_{\, t} +  O(|||\varphi|||^2_{\, t})$$ and
$$\sum_{i, j, k} Re \{(\varphi_i, \lbar{k}(\bar b^k_{ij}
\varphi_j))_t \} \geq - \eps  \sum_{j, k} |||\lbar{k}
\varphi_j|||^2_{\, t} +  O(|||\varphi|||^2_{\, t}).$$ Let $j$ be
$i$ and $k$ be $j$ on the right-hand sides of these two
inequalities. Then, we plug these into the expression for $\qbtp
{\varphi},$
\begin{equation*}
\begin{split}
\qbtp{\varphi} &\geq (1 - 3 \eps)\, \sum_{i, j} |||\lbar{j}
\varphi_i|||^2_{\, t} + \sum_{i, j} \Re \{(c_{ij}T \varphi_i,
\varphi_j)_t \}+ t \sum_{i, j}\Re \{( \lbar{j} L_i (\lambda)
\varphi_i, \varphi_j)_t \}\\&\quad+ t \sum_{ij}\Re \{ ( e_{ij}
\varphi_i, \varphi_j)_t \}+ O(|||\varphi|||^2_{\, t}),
\end{split}
\end{equation*}
and we let $\eps' \, = \, 3 \eps.$ \qed
\smallskip
\newtheorem{qbtmcalc}[lemma]{Lemma}
\begin{qbtmcalc}
Let $\varphi$ be a $(0,1)$ form \label{qbtm} supported in $U'$
such that $\varphi \in Dom(\dbarb) \cap Dom(\ad).$ If $[L_i,
\lbar{j}]= c_{ij}T + \sum_k a^k_{ij}\lbar{k} + \sum_k b^k_{ij}
\lbar{k}^{*,-t} + g_{ij} - t \, e_{ij},$ where $a^k_{ij}$,
$b^k_{ij}$, $g_{ij},$ and $e_{ij}$ are \smooth functions
independent of $t$ and $c_{ij}$ are the coefficients of the Levi
form, then there exists $1 \gg \eps" >0$ for which
\begin{equation*}
\begin{split}
\qbtm{\varphi} &\geq (1- \eps")\, \sum_{ij} |||\lbar{j}^{*, -t}
\varphi_i|||^2_{-t} + \sum_{ij} \Re \{(c_{ij} T \varphi_i,
\varphi_j)_{-t} \} - \sum_{j} \Re \{((\sum_i c_{ii}) T \varphi_j,
\varphi_j)_{-t}\}\\ & \quad- t \sum_{ij} \Re \{(\lbar{j} L_i
(\lambda) \varphi_i, \varphi_j)_{-t}\}+ t \sum_{j} \Re \{((\sum_i
\lbar{i} L_i (\lambda)) \varphi_j, \varphi_j)_{-t}\}\\ & \quad- t
\sum_{ij} \Re \{(e_{ij} \varphi_i, \varphi_j)_{-t}\}+ t \sum_{j}
\Re \{((\sum_i e_{ii}) \varphi_j, \varphi_j)_{-t}\} +
O(|||\varphi|||^2_{-t} ).
\end{split}
\end{equation*}
\end{qbtmcalc}
\noindent {\bf{Proof:}} Again, since $\varphi \in Dom(\ad),$
$\varphi \in Dom(\adm).$ Notice that we can use the same proof as
in the previous lemma with $t$ replaced by $-t$ and integrate the cross terms by parts to show that
\begin{equation*}
\begin{split}
\qbtm{\varphi} &=\sum_{i, j} |||\lbar{j} \varphi_i|||^2_{- t}  -
\sum_{i, j} ( [\lbar{i}^{*, \,- t}, \lbar{j}] \varphi_i,
\varphi_j)_{-t}\\&\quad + O\bigg(|||\varphi|||^2_{- t} +
\big(\sum_{i,j} |||\lbar{i}^{*,-t} \varphi_j|||^2_{- t}
\big)^{\frac{1}{2}}\, |||\varphi|||_{-t} \bigg).
\end{split}
\end{equation*}
Then we rewrite $|||\lbar{j} \varphi_i|||^2_{-t}$ in a convenient
form,
\begin{equation*}
\begin{split}
|||\lbar{j} \varphi_i|||^2_{-t} &= (\lbar{j} \varphi_i, \lbar{j}
\varphi_i)_{-t} = (\lbar{j}^{*, -t} \lbar{j} \varphi_i,
\varphi_i)_{-t} = (\lbar{j} \lbar{j}^{*, -t} \varphi_i,
\varphi_i)_{-t} + ([\lbar{j}^{*,-t} ,\lbar{j}] \varphi_i,
\varphi_i)_{-t}\\ &= (\lbar{j}^{*, -t} \varphi_i, \lbar{j}^{*, -t}
\varphi_i)_{-t} + ([\lbar{j}^{*,-t}, \lbar{j}] \varphi_i,
\varphi_i)_{-t}.
\end{split}
\end{equation*}
Since $\lbar {j}^{*,-t} \,  = \,  -  L_j  - f_j - t L_j
(\lambda)$, $[\lbar{j}^{*,-t}, \lbar{j}] \, = \, [-  L_j  - f_j -
t L_j (\lambda), \lbar{j}] \, = \, - [L_j, \lbar{j}] + \lbar{j}
(f_j) + t \lbar{j} L_j (\lambda)$, so using the hypothesis of the
lemma about $[L_j, \lbar{j}]$ and the same method as in the
previous lemma we obtain that
\begin{equation*}
\begin{split}
|||\lbar{j} \varphi_i|||^2_{-t} &= |||\lbar{j}^{*,-t}
\varphi_i|||^2_{-t} + ([\lbar{j}^{*,-t}, \lbar{j}] \varphi_i,
\varphi_i)_{-t} \\ &\geq (1-2\eps')\, |||\lbar{j}^{*,-t}
\varphi_i|||^2_{-t} - \Re \{ (c_{jj} T \varphi_i,
\varphi_i)_{-t}\}+ t \, \Re \{(\lbar{j} L_j (\lambda) \varphi_i,
\varphi_i)_{-t}\}\\&\quad+ t \, \Re \{(e_{jj} \varphi_i,
\varphi_i)_{-t}\} + O(|||\varphi|||^2_{-t}).
\end{split}
\end{equation*}
Plug into the expression for $\qbtm{\varphi},$ unravel the commutator $[\lbar{i}^{*,-t}, \lbar{j}]$ as above,
and again use the techniques from the previous lemma to get:
\begin{equation*}
\begin{split}
\qbtm{\varphi} &\geq (1-\eps")\, \sum_{ij} |||\lbar{j}^{*, -t}
\varphi_i|||^2_{-t} - \sum_{ij} \Re \{(c_{jj} T \varphi_i,
\varphi_i)_{-t}\} + \sum_{ij}\Re \{ (c_{ij} T \varphi_i,
\varphi_j)_{-t}\}\\& \quad+ t \sum_{ij}\Re \{ (\lbar{j} L_j
(\lambda) \varphi_i, \varphi_i)_{-t}\}- t \sum_{ij}\Re \{
(\lbar{j} L_i (\lambda)\varphi_i, \varphi_j)_{-t} \}\\
& \quad+ t \sum_{ij}\Re \{ (e_{jj} \varphi_i, \varphi_i)_{-t}\}- t
\sum_{ij}\Re \{ (e_{ij} \varphi_i, \varphi_j)_{-t} \}+
O(|||\varphi|||^2_{-t} ).
\end{split}
\end{equation*}
Exchange $i$ and $j$ in the second, the fourth, and the sixth
terms of the right-hand side and pull the summation under $i$
inside the parentheses to finish off the proof of this lemma. \qed

\medskip
\medskip
\noindent The upcoming two lemmas concern applications of the G{\aa}rding
inequality for positive semi-definite matrix operators on $\cp$
and $\cm.$

\medskip
\medskip
\smallskip
\newtheorem{cpGarding}[lemma]{Lemma}
\begin{cpGarding}
Let $M$ be a weakly pseudoconvex CR-manifold and $\varphi$ be a
$(0,1)$ form supported on $U'$ such that up to a smooth term,
$\hat \varphi$ is supported in $\cp,$ then \label{Gardcp} the
following is true:
$$\sum_{i,j} \Re \{(c_{ij} T \varphi_i, \varphi_j)_t\} \geq
 tA \sum_{i,j} (c_{ij} \varphi_i,
\varphi_j)_t+O(|||\varphi|||^2_{\, t})+O_t(||\dze \potd
\varphi||^2_0),$$ where $c_{ij}$ are the coefficients of the Levi
form.
\end{cpGarding}
\noindent {\bf{Proof:}} Let $\pptd$ be a
pseudodifferential operator of order zero whose symbol dominates
$\hat \varphi$ (up to a smooth error) and is supported in $\cpd$, a slightly
larger truncated cone in the $\xi$ (Fourier transform) space than
$\cp$. We proceed as follows:
\begin{equation*}
\begin{split}
\sum_{i,j} (c_{ij} T \varphi_i, \varphi_j)_t &= \sum_{i,j} (\dze e^{-t\lambda}
c_{ij} T \pptd \varphi_i, \dze \pptd \varphi_j)_{\,0}+ smooth \: terms
\\&=\sum_{i,j} (\dze (\pptd)^* \dze^2 e^{-t\lambda}
c_{ij} T \pptd \varphi_i, \varphi_j)_{\,0}+ smooth \: terms
\end{split}
\end{equation*}
We need to analyze the symbol of $ (\pptd)^* \dze^2 e^{-t\lambda}
c_{ij} T \pptd.$ We look first at $\sigma(T \pptd)$ and then at $\sigma((\pptd)^* \dze^2 e^{-t\lambda}
c_{ij}).$ Let ${\widetilde{\psi}}_t^+ (x, \xi)$ be the symbol of $\pptd.$ We note that $\sigma(T) = \xi_{2n-1}.$
The symbol of the composition $T \pptd$ is then given by
\begin{equation*}
\sum_{\beta} \frac {1}{\beta !} \:
\partial_{\xi}^{\beta} (\xi_{2n-1}) D_x^{\beta} {\widetilde{\psi}}_t^+ = \xi_{2n-1}
{\widetilde{\psi}}_t^+ (x, \xi) + D_x^{(0, \dots, 0,1)}
 {\widetilde{\psi}}_t^+ (x, \xi)
\end{equation*}
because for any other $\beta$, $\partial_{\xi}^{\beta}
(\xi_{2n-1}) \, = \, 0$. Now, since $\hat \varphi$ is supported in
$\cp$ (up to a smooth term) and ${\widetilde{\psi}}_t^+ (x, \xi)
\, \equiv \, 1$ on $\cp$, any of its derivatives will be zero, so
$ \sigma(T \pptd) \, =  \, \xi_{2n-1} {\widetilde{\psi}}_t^+
(x, \xi)$ up to smooth terms when applied to $\varphi.$
Let us show that $\sigma((\pptd)^* \dze^2 e^{-t\lambda} c_{ij}) \, = \,
{\widetilde{\psi}}_t^+ (x, \xi)\dze^2 e^{-t\lambda} c_{ij} $ on $\cp.$ First, since the symbol
of $\pptd$ is identically $1$ on $\cp,$ it follows the symbol of
$(\pptd)^*$ must also be ${\widetilde{\psi}}_t^+ (x, \xi) $ up to
terms supported in $\co \setminus \cp.$ This implies
\begin{equation*}
\sigma((\pptd)^* \dze^2 e^{-t\lambda} c_{ij}) =  \sum_{\beta} \frac {1}{\beta !} \:
\partial_{\xi}^{\beta} {\widetilde{\psi}}_t^+ (x, \xi) D_x^{\beta} \,(\dze^2 e^{-t\lambda} c_{ij}).
\end{equation*}
$\partial_{\xi}^{\beta} {\widetilde{\psi}}_t^+ (x, \xi) \, = \, 0$
if $\beta \neq 0$ on $\cp,$ so indeed $\sigma((\pptd)^* \dze^2 e^{-t\lambda} c_{ij}) \,
= \, {\widetilde{\psi}}_t^+ (x, \xi) \dze^2 e^{-t\lambda} c_{ij}$ up to errors on $\co
\setminus \cp.$ Next, let us compute $\sigma((\pptd)^* \dze^2 e^{-t\lambda} c_{ij} T
\pptd)$:
\begin{equation*}
\begin{split}
\sigma((\pptd)^* \dze^2 e^{-t\lambda} c_{ij} T \pptd) &=  \sum_{\alpha} \frac
{1}{\alpha !} \: \partial_{\xi}^{\alpha} \sigma((\pptd)^* \dze^2 e^{-t\lambda} c_{ij})
D_x^{\alpha}\sigma(T \pptd) \\&=  \sum_{\alpha} \frac {1}{\alpha !}
\: \partial_{\xi}^{\alpha} ({\widetilde{\psi}}_t^+ (x, \xi)\dze^2 e^{-t\lambda}
c_{ij}) D_x^{\alpha}(\xi_{2n-1}{\widetilde{\psi}}_t^+ (x, \xi))\\
&={\widetilde{\psi}}_t^+ (x, \xi)\dze^2 e^{-t\lambda}
c_{ij}\xi_{2n-1}{\widetilde{\psi}}_t^+ (x, \xi) =\dze^2 e^{-t\lambda} c_{ij}\xi_{2n-1}
\end{split}
\end{equation*}
on $\cp$ because $\dze^2 e^{-t\lambda} c_{ij}$ is independent of $\xi$ and
${\widetilde{\psi}}_t^+ (x, \xi) \,  \equiv \, 1.$ So now
since $\sigma(T) \, = \, \xi_{2n-1} \geq tA$ on $\cp$ and $(\dze^2 e^{-t\lambda} c_{ij})$ is positive semi-definite,
we can apply Lemma~\ref{Gardsecond} with $T$ as $R$ and $(\dze^2 e^{-t\lambda} c_{ij})$ as
$(h_{ij})$ to conclude that there exists a constant $C$
independent of $t$ such that:
\begin{equation*}
\begin{split}
\sum_{i,j}\Re \{ (c_{ij} T \varphi_i, \varphi_j)_t\} &\geq tA
\sum_{i,j} (\dze^2 e^{-t\lambda} c_{ij} \varphi_i,
\varphi_j)_0+O(|||\varphi|||^2_{\, t})-C|||\varphi|||^2_{\, t}
+O_t(||\dze \potd \varphi||^2_0)\\
&= tA \sum_{i,j}
(c_{ij}\varphi_i,\varphi_j)_t+O(|||\varphi|||^2_{\, t})+O_t(||\dze
\potd \varphi||^2_0),
\end{split}
\end{equation*}
where all the errors whose Fourier transforms are supported in
$\co \setminus \cp$ have been included in $O_t(||\dze \potd
\varphi||^2_0),$ an error term which we can afford given the setup of this argument, and all the smooth
errors have been included in $O(||\varphi||^2_{\, t}).$ \qed

\medskip
\medskip
\smallskip
\newtheorem{cmGarding}[lemma]{Lemma}
\begin{cmGarding}
Let $\varphi$ be a $(0,1)$ form supported on $U'$ such that up to
a smooth term, $\hat \varphi$ is supported in $\cm$ and $(g_{ij})$
a Hermitian, positive semi-definite matrix, then \label{Gardcm} the
following is true:
$$\sum_{i,j} \Re \{(g_{ij} (-T) \varphi_i, \varphi_j)_{-t}\} \geq  t A
\sum_{i,j} (g_{ij} \varphi_i, \varphi_j)_{-t}+O(|||\varphi|||^2_{-
t})+O_t(||\dze \potd \varphi||^2_0).$$
\end{cmGarding}
\noindent {\bf{Proof:}} Notice that we can use the proof of the
previous lemma with the obvious modifications which do not affect
any steps there: Use $(\, \cdot \, , \, \cdot \,)_{-t}$ instead of
$(\, \cdot \, , \, \cdot \,)_t$, $(g_{ij})$ instead of the
coefficients of the Levi form, and $-T$ instead of $T$ because on
$\cm$, $\sigma(-T) \, = \, - \xi_{2n-1} \geq tA$. \qed

\medskip
\medskip
\noindent Lastly, here is a very elementary linear algebra lemma which will
be needed to prove the estimate on $\cm$.
\smallskip
\newtheorem{linalglemma}[lemma]{Lemma}
\begin{linalglemma}
If $n>2$ and $ H \, = \, (h_{ij})$ is a Hermitian, positive
definite matrix, \label{linalg} $1 \leq i,j \leq n-1,$ then
$(\delta_{ij} \sum_i h_{ii} - h_{ij})$ is also Hermitian and
positive definite. For any $n$, if $H$ is positive semi-definite,
then $(\delta_{ij} \sum_i h_{ii} - h_{ij})$ is also positive
semi-definite.
\end{linalglemma}
\smallskip
\noindent {\bf{Proof:}} Diagonalize $(h_{ij})$ at a point to get a
matrix:
$$ G =
\begin{pmatrix}
g_{11} & 0 & \dots & 0 \\
0 & g_{22} & \dots & 0 \\
\vdots & \vdots & \ddots & \vdots \\
0 & 0 & \dots & g_{n-1 n-1}\\
\end{pmatrix} $$
So $\sum_i g_{ii} \, = \, \sum_i h_{ii}$. Let $C$ be the
transition matrix and $I$ the identity matrix, then we have that:
\begin{equation*}
\begin{split}
C^{-1}((\sum_i h_{ii})I - H)C &= C^{-1}(\sum_i h_{ii})IC -
C^{-1}HC = (\sum_i h_{ii})C^{-1}C - G = (\sum_i g_{ii})I - G \\ &=
\begin{pmatrix}
(\sum_i g_{ii})-g_{11} & 0 & \dots & 0 \\
0 & (\sum_i g_{ii})-g_{22} & \dots & 0 \\
\vdots & \vdots & \ddots & \vdots \\
0 & 0 & \dots & (\sum_i g_{ii})-g_{n-1 n-1}\\
\end{pmatrix}
\end{split}
\end{equation*}
Therefore, if $H$ has more than one eigenvalue, i.e. $n-1>1 \,
\Leftrightarrow \, n>2$, then $(h_{ij})$ positive definite implies
$\sum_i h_{ii} > 0 \, \Leftrightarrow \, \sum_i g_{ii}>0$, so
$\sum_i g_{ii} - g_{jj} > 0$ for at least one $j, \, 1 \leq j \leq
n-1,$ which means $(\delta_{ij} \sum_i h_{ii} - h_{ij})$ is also
positive definite. Now, for any $n$, if we just have that $H$ is
positive semi-definite, i.e. $\sum_i h_{ii} \geq 0 \,
\Leftrightarrow \, \sum_i g_{ii} \geq 0$, then $\sum_i g_{ii} -
g_{jj} \geq 0$ for all $j$, which implies $(\delta_{ij} \sum_i
h_{ii} - h_{ij})$ is positive semi-definite. \qed

\medskip
\medskip
Finally, all the necessary lemmas have been proven and we can
state and prove the two local results that will be used in the
proof of the main estimate.
\smallskip
\newtheorem{locplusmicro}[lemma]{Proposition}
\begin{locplusmicro}
Let $\varphi  \in Dom(\dbarb) \cap Dom(\ad)$ be a $(0,1)$ form
supported in $U$, a neighborhood of a compact, weakly pseudoconvex
\label{localplus} CR-manifold $M.$ $U$ is taken small enough that
the pseudodifferential operators $\ppt$, $\pot$, and $\pmt$
described at the beginning of this section can be defined. Let
there exist an orthonormal basis $L_1, \dots, L_{n-1}$ of vector
fields on $T^{1,0}(U)$ such that $[L_i, \lbar{j}]= c_{ij}T +
\sum_k a^k_{ij}\lbar{k} + \sum_k b^k_{ij} \lbar{k}^{*,t} + g_{ij}+
t e_{ij},$ where $a^k_{ij},$ $b^k_{ij},$ $g_{ij},$ and $e_{ij}$
are \smooth functions independent of $t,$ $|e_{ij}| < \eps_G$ for
a very small $\eps_G,$ and $c_{ij}$ are the coefficients of the
Levi form. Also, let $M$ be endowed with a strongly CR
plurisubharmonic function $\lambda$, i.e. $\, \exists$ $A_0>0$
such that the matrix with entries $s_{ij} \, = \, \frac {1}{2} \,
( \lbar{j}L_i (\lambda)+ L_i \lbar{j} (\lambda)) +A_0\, c_{ij}$ is
positive definite at each point $x \in M$, as described in
Section~\ref{CRpsh}. Then, there exists a constant $C_1$
independent of $t$ such that
\begin{equation}
\qbtp{\dze \ppt \varphi} + O(|||\dze \ppt \varphi|||^2_{\,
t})+O_t(||\dze \potd \varphi||^2_0) \geq C_1 \, t|||\dze \ppt
\varphi|||^2_{\, t}. \label{localpluseq}
\end{equation}
\end{locplusmicro}
\noindent {\bf{Proof:}} The three pseudodifferential operators
$\ppt$, $\pot$, and $\pmt$ are defined using $A \, = \,A_0$,
where $A_0$ is the positive constant given by the definition of
the CR-plurisubharmonicity of $\lambda$. Next, since $\varphi \in
Dom(\dbarb) \cap Dom(\ad),$ it easily follows that $\dze \ppt
\varphi \in Dom(\dbarb) \cap Dom(\ad).$ Moreover, $\dze \ppt
\varphi$ is supported in $U'$ and the bracket $[L_i, \lbar{j}]$
has the required property, so we apply Lemma~\ref{qbtp} to $\qbtp
{\dze \ppt \varphi}$ to conclude that for some $1 \gg \eps'
>0$:
\begin{equation*}
\begin{split}
\qbtp{\dze \ppt \varphi} &\geq (1-\eps')\, \sum_{ij} |||\lbar{j}
\dze \ppt \varphi_i|||^2_{\, t} + \sum_{ij} \Re \{(c_{ij} T \dze
\ppt \varphi_i, \dze \ppt \varphi_j)_t \}\\&\quad+ t \sum_{ij}\Re
\{ ( \lbar{j} L_i (\lambda) \dze \ppt \varphi_i, \dze \ppt
\varphi_j)_t \}+ t \sum_{ij}\Re \{ ( e_{ij} \dze \ppt \varphi_i,
\dze \ppt \varphi_j)_t \}\\&\quad + O(|||\dze \ppt
\varphi|||^2_{\, t} ).
\end{split}
\end{equation*}
We rewrite the term $t \sum_{ij} \Re \{(\lbar{j} L_i (\lambda)
\dze \ppt \varphi_i, \dze \ppt \varphi_j)_t \}$ using the properties of the Hermitian inner product as follows:
\begin{equation*}
\begin{split}
t \sum_{ij} \Re \{(\lbar{j} L_i (\lambda) \dze \ppt \varphi_i,
\dze \ppt \varphi_j)_t \} &= t \sum_{ij}(\half (\lbar{j} L_i
(\lambda) +L_i \lbar{j} (\lambda)) \dze \ppt \varphi_i, \dze \ppt
\varphi_j)_t
\end{split}
\end{equation*}
Then apply Lemma~\ref{Gardcp} to the term $\sum_{ij} \Re \{(c_{ij}
T \dze \ppt \varphi_i, \dze \ppt \varphi_j)_t\}$ since $\dze$ is
supported in $U',$ $M$ is weakly pseudoconvex and up to a smooth
term, the Fourier transform of $\dze \ppt \varphi$ is supported in
$\cp$:
$$\sum_{i,j}\Re \{ (c_{ij} T \dze \ppt \varphi_i, \dze \ppt \varphi_j)_t \} \geq
 t\, A_0 \sum_{i,j} (c_{ij} \dze \ppt \varphi_i, \dze
\ppt \varphi_j)_t+O(|||\dze \ppt \varphi|||^2_{\, t})+O_t(||\dze
\potd \varphi||^2_0).$$ We put these two together:
\begin{equation*}
\begin{split}
\qbtp{\dze \ppt \varphi} &\geq (1 -\eps')\,\sum_{ij} |||\lbar{j}
\dze \ppt \varphi_i|||^2_{\, t}+t A_0 \sum_{i,j} (c_{ij} \dze \ppt
\varphi_i, \dze \ppt \varphi_j)_t\\ & \quad+  t \sum_{ij}(\half
(\lbar{j} L_i (\lambda) +L_i \lbar{j} (\lambda)) \dze \ppt
\varphi_i, \dze \ppt \varphi_j)_t+ O(|||\dze \ppt \varphi|||^2_{\,
t} )\\&\quad + t \sum_{ij}\Re \{ ( e_{ij} \dze \ppt \varphi_i,
\dze \ppt \varphi_j)_t \}+O_t(||\dze \potd \varphi||^2_0)\\ &\geq
t \sum_{ij}((\half (\lbar{j} L_i (\lambda) +L_i \lbar{j}
(\lambda))+
A_0 \, c_{ij}) \dze \ppt \varphi_i, \dze \ppt \varphi_j)_t \\
& \quad+ t \sum_{ij}\Re \{ ( e_{ij} \dze \ppt \varphi_i, \dze \ppt
\varphi_j)_t \}+O(|||\dze \ppt \varphi|||^2_{\, t}
)\\&\quad+O_t(||\dze \potd \varphi||^2_0) \\&\geq t (C-\eps_G)\,
|||\dze \ppt \varphi|||^2_{\, t} + O(|||\dze \ppt \varphi|||^2_{\,
t})+O_t(||\dze \potd \varphi||^2_0),
\end{split}
\end{equation*}
because $\lambda$ is strongly CR plurisubharmonic and $|e_{ij}| <
\eps_G,$ where $1 \gg \eps_G >0.$ Thus, there exists a constant
$C_1$ independent of $t$ such that Equation~\ref{localpluseq} holds:
\begin{equation*}
 \qbtp{\dze \ppt \varphi}  + O(|||\dze \ppt \varphi|||^2_{\,
t})+O_t(||\dze \potd \varphi||^2_0) \geq C_1 \, t|||\dze \ppt
\varphi|||^2_{\, t}
\end{equation*} \qed

\medskip
\newtheorem{locminusmicro}[lemma]{Proposition}
\begin{locminusmicro}
Let $\varphi  \in Dom(\dbarb) \cap Dom(\ad)$ be a $(0,1)$ form
supported in $U$, a neighborhood of a compact, weakly pseudoconvex
\label{localminus} CR-manifold $M$ of dimension at least $5.$ $U$
is taken small enough that the pseudodifferential operators
$\ppt$, $\pot$, and $\pmt$ described at the beginning of this
section can be defined. Let there exist an orthonormal basis $L_1,
\dots, L_{n-1}$ of vector fields on $T^{1,0}(U)$ such that $[L_i,
\lbar{j}]= c_{ij}T + \sum_k a^k_{ij}\lbar{k} + \sum_k b^k_{ij}
\lbar{k}^{*,-t} + g_{ij}- t e_{ij},$ where $a^k_{ij},$ $b^k_{ij},$
$g_{ij},$ and $e_{ij}$ are \smooth functions independent of $t,$
$|e_{ij}|<\eps_G$ for a very small positive number $\eps_G,$ and
$c_{ij}$ are the coefficients of the Levi form. Also, let $M$ be
endowed with a strongly CR plurisubharmonic function $\lambda$,
i.e. $\, \exists$ $A_0>0$ such that the matrix with entries
$s_{ij} \, = \, \frac {1}{2} \, ( \lbar{j}L_i (\lambda)+ L_i
\lbar{j} (\lambda)) +A_0\, c_{ij}$ is positive definite at each
point $x \in M$, as described in Section~\ref{CRpsh}. Then, there
exists a constant $C_2$ independent of $t$ such that
\begin{equation}
\label{localminuseq} \qbtm{\dze \pmt \varphi} + O(|||\dze \pmt
\varphi|||^2_{- t})+O_t(||\dze \potd \varphi||^2_0) \geq C_2 \, t
|||\dze \pmt \varphi |||^2_{-t}.
\end{equation}
\end{locminusmicro}
\noindent {\bf{Proof:}} Just as in the previous proof, the three
pseudodifferential operators $\ppt$, $\pot$, and $\pmt$ are
defined using $A \, = \, A_0$, where $A_0$ is the positive
constant given by the definition of the CR-plurisubharmonicity of
$\lambda.$ Next, $\varphi \in Dom(\dbarb) \cap Dom(\ad)$ implies
$\dze \pmt \varphi \in Dom(\dbarb) \cap Dom(\ad).$ Also $\dze \pmt
\varphi$ is supported in $U',$ so we can apply Lemma~\ref{qbtm} to
it which gives some $1 \gg \eps" >0$ for which:
\begin{equation*}
\begin{split}
\qbtm{\dze \pmt \varphi} &\geq (1- \eps")\, \sum_{ij}
|||\lbar{j}^{*, -t} \dze \pmt \varphi_i|||^2_{-t} + \sum_{ij} \Re
\{(c_{ij} T \dze \pmt \varphi_i, \dze \pmt \varphi_j)_{-t} \}\\ &
\quad - \sum_{j} \Re \{((\sum_i c_{ii}) T \dze \pmt \varphi_j,
\dze \pmt \varphi_j)_{-t}\}\\ & \quad- t \sum_{ij} \Re \{(\lbar{j}
L_i (\lambda) \dze \pmt \varphi_i, \dze \pmt \varphi_j)_{-t}\}\\ &
\quad+ t \sum_{j} \Re
\{((\sum_i \lbar{i} L_i (\lambda)) \dze \pmt \varphi_j, \dze \pmt \varphi_j)_{-t}\}\\
& \quad- t \sum_{ij} \Re \{(e_{ij} \dze \pmt \varphi_i, \dze \pmt
\varphi_j)_{-t}\}+ t \sum_{j} \Re \{((\sum_i e_{ii}) \dze \pmt
\varphi_j, \dze \pmt \varphi_j)_{-t}\}\\ & \quad + O(|||\dze \pmt
\varphi|||^2_{-t} ).
\end{split}
\end{equation*}
Rewriting two of the terms on the right-hand side gives:
\begin{equation*}
\begin{split}
&- \sum_{j} \Re \{ ((\sum_i c_{ii}) T \dze \pmt \varphi_j, \dze
\pmt \varphi_j)_{-t}\}+ \sum_{ij} \Re \{(c_{ij} T \dze \pmt
\varphi_i, \dze \pmt \varphi_j)_{-t} \} \\ &= \sum_{j} \Re\{
((\sum_i c_{ii})( -T) \dze \pmt \varphi_j, \dze \pmt
\varphi_j)_{-t} \}- \sum_{ij} \Re \{(c_{ij} (-T) \dze \pmt
\varphi_i, \dze \pmt \varphi_j)_{-t} \} \\ &= \sum_{ij} \Re \{
((\delta_{ij} \sum_i c_{ii} - c_{ij})( -T) \dze \pmt \varphi_i,
\dze \pmt \varphi_j)_{-t} \}
\end{split}
\end{equation*}
Since $(\delta_{ij} \sum_i c_{ii} - c_{ij})$ is positive
semi-definite by Lemma~\ref{linalg}, we can apply
Lemma~\ref{Gardcm} to conclude that
\begin{equation*}
\begin{split}
\sum_{ij} \Re \{((\delta_{ij} \sum_i c_{ii} - c_{ij})( -T) \dze
\pmt \varphi_i, \dze \pmt \varphi_j)_{-t} \} &\geq t \,A_0
\sum_{ij} ((\delta_{ij} \sum_i c_{ii} - c_{ij}) \dze \pmt
\varphi_i, \dze \pmt \varphi_j)_{-t} \\ & \quad+ O(|||\dze \pmt
\varphi|||^2_{-t} )+O_t(||\dze \potd \varphi||^2_0).
\end{split}
\end{equation*}
Also, by a similar manipulation to the one in
Proposition~\ref{localplus},
\begin{equation*}
\begin{split}
&t\sum_{j} \Re \{((\sum_i \lbar{i} L_i (\lambda)) \dze \pmt
\varphi_j, \dze \pmt \varphi_j)_{-t} \}- t \sum_{ij}\Re \{
(\lbar{j} L_i (\lambda) \dze \pmt \varphi_i, \dze \pmt
\varphi_j)_{-t} \} \\ &= t \sum_{ij} \Re \{ ((\delta_{ij}\sum_i
\lbar{i} L_i (\lambda)-\lbar{j} L_i (\lambda)) \dze \pmt
\varphi_i, \dze \pmt \varphi_j)_{-t} \} \\& = t \sum_{ij}
([\delta_{ij}\sum_i \half ( \lbar{i} L_i (\lambda)+ L_i \lbar{i}
(\lambda))- \half (\lbar{j} L_i (\lambda)+ L_i \lbar{j} (\lambda))
] \dze \pmt \varphi_i, \dze \pmt \varphi_j)_{-t}
\end{split}
\end{equation*}
We then plug these two into the expression for $\qbtm{\dze \pmt
\varphi}$:
\begin{equation*}
\begin{split}
&\qbtm{\dze \pmt \varphi} \geq (1-\eps") \, \sum_{ij}
|||\lbar{j}^{*, -t} \dze \pmt \varphi_i|||^2_{-t}+ O(|||\dze \pmt
\varphi|||^2_{-t} )+O_t(||\dze \potd \varphi||^2_0) \\& \quad+t
\sum_{ij} ([\delta_{ij}\sum_i \half ( \lbar{i} L_i (\lambda)+ L_i
\lbar{i} (\lambda))- \half (\lbar{j} L_i (\lambda)+ L_i \lbar{j}
(\lambda)) ] \dze \pmt \varphi_i, \dze \pmt \varphi_j)_{-t}\\&
\quad- t \sum_{ij} \Re \{(e_{ij} \dze \pmt \varphi_i, \dze \pmt
\varphi_j)_{-t}\}+ t \sum_{j} \Re \{((\sum_i e_{ii}) \dze \pmt
\varphi_j, \dze \pmt \varphi_j)_{-t}\} \\ & \quad+t \, A_0
\sum_{ij} ((\delta_{ij} \sum_i
c_{ii} - c_{ij}) \dze \pmt \varphi_i, \dze \pmt \varphi_j)_{-t}\\
&\geq t \sum_{ij} ([\delta_{ij}\sum_i \half ( \lbar{i} L_i
(\lambda)+ L_i \lbar{i} (\lambda))- \half (\lbar{j} L_i (\lambda)+
L_i \lbar{j} (\lambda)) ] \dze \pmt \varphi_i, \dze \pmt
\varphi_j)_{-t}\\&\quad+ t \, A_0 \sum_{ij} ((\delta_{ij} \sum_i
c_{ii} - c_{ij}) \dze \pmt \varphi_i, \dze \pmt \varphi_j)_{-t}- t
\sum_{ij} \Re \{(e_{ij} \dze \pmt
\varphi_i, \dze \pmt \varphi_j)_{-t}\}\\
& \quad+ t \sum_{j} \Re \{((\sum_i e_{ii}) \dze \pmt \varphi_j,
\dze \pmt \varphi_j)_{-t}\} + O(|||\dze \pmt \varphi|||^2_{-t}
)+O_t(||\dze \potd \varphi||^2_0)
\end{split}
\end{equation*}
Since $\lambda$ is strongly CR plurisubharmonic, the matrix with
entries $w_{ij} \, = \,A_0 (\delta_{ij} \sum_i c_{ii} - c_{ij})
+\delta_{ij}\sum_i \half ( \lbar{i} L_i (\lambda)+ L_i \lbar{i}
(\lambda))- \half (\lbar{j} L_i (\lambda)+ L_i \lbar{j}
(\lambda))$ is positive definite for $n>2$ by Lemma~\ref{linalg},
i.e. if the CR-manifold $M$ has at least dimension $5$ which holds
by the hypothesis, and $|e_{ij}|<\eps_G \ll 1,$ therefore there
exists a constant $C_2$ independent of $t$ such that:
\begin{equation*}
\qbtm{\dze \pmt \varphi}+ O(|||\dze \pmt \varphi|||^2_{-t}
)+O_t(||\dze \potd \varphi||^2_0) \geq C_2 \, t|||\dze \pmt
\varphi|||^2_{-t}.
\end{equation*}
This is exactly Equation~\ref{localminuseq}.
\qed


\bigskip
\noindent We can now prove the global estimate stated at the
beginning of the section:

\medskip
\noindent {\bf{Proof of the main estimate:}} The norm $\tnorm {\,
\cdot \,}{\, \cdot \,}$ is defined using a covering of $M$,
$\{U_\nu \}_\nu$ consisting of neighborhoods on which all the
extra hypotheses of the local results, Proposition~\ref{localplus}
and Proposition~\ref{localminus}, held. Moreover, we take $A \, =
\,  A_0,$ where $A_0$ is the constant of strong CR
plurisubharmonicity of $\lambda,$ in the construction of the
pseudodifferential operators in each $U_\nu$ as explained at the
beginning of the proofs of these two propositions. We can thus
apply these two results to conclude that for each $\nu$, there
exist constants $C_{\nu, 1}$ and $C_{\nu, 2}$ independent of $t$
such that
\begin{equation*}
\qbtp{\dze \ppt \varphi}  + O(|||\dze \ppt \varphi|||^2_{\,
t})+O_t(||\dze \potd \varphi||^2_0) \geq C_{\nu, 1} \, t|||\dze
\ppt \varphi|||^2_{\, t}
\end{equation*}
and
\begin{equation*}
\qbtm{\dze \pmt \varphi} + O(|||\dze \pmt \varphi|||^2_{-
t})+O_t(||\dze \potd \varphi||^2_0) \geq C_{\nu, 2} \, t |||\dze
\pmt \varphi |||^2_{-t}.
\end{equation*}
Equation~\ref{qbtexpression} which is the expression for
$\qbtg{\varphi}$ can be rewritten as
\begin{equation*}
\begin{split}
&K \qbtg {\varphi}+K_t \sum_\nu \:{||\dze_\nu \potdg {\nu} \ze_\nu
\varphi^\nu||}^2_0+O({\langle |\varphi|\rangle}^2_t) +O_t
(||\varphi||^2_{-1})\\&\geq \sum_\nu \: \qbtp{\dze_\nu \pptg {\nu}
\ze_\nu \varphi^\nu} + \sum_\nu \: \qbtm{\dze_\nu \pmtg {\nu}
\ze_\nu \varphi^\nu}.
\end{split}
\end{equation*}
Notice that since $\potdg{\nu}$ dominates $\potg{\nu}$ and $K_t$
grows exponentially with $t,$ $K_t \sum_\nu \:{||\dze_\nu \potdg
{\nu} \ze_\nu \varphi^\nu||}^2_0$ controls $\sum_\nu \: t \,
{||\dze_\nu \potg {\nu} \ze_\nu \varphi^\nu||}^2_0$ up to errors
of the type $ O_t (||\varphi||^2_{-1}).$ Let $$C \, = \, \min_\nu
\{ \min \{ C_{\nu,1}, C_{\nu,2},1 \}\},$$ which exists and is
positive because the covering $\{U_\nu\}_\nu$ is finite, then
Proposition~\ref{localplus} and Proposition~\ref{localminus} along
with the observation made above imply that if we increase $K_t$
slightly:
\begin{equation*}
\begin{split}
&K \qbtg {\varphi}+K_t \sum_\nu \:{||\dze_\nu \potdg {\nu} \ze_\nu
\varphi^\nu||}^2_0+O({\langle |\varphi|\rangle}^2_t) +O_t
(||\varphi||^2_{-1})\\& \geq C \, \sum_\nu \: t ({|||\dze_\nu
\pptg {\nu} \ze_\nu {\varphi}|||}^2_t + {||\dze_\nu \potg {\nu}
\ze_\nu {\varphi}||}^2_0 + {|||\dze_\nu \pmtg {\nu} \ze_\nu
{\varphi}|||}^2_{-t}) = C \, t {\langle |\varphi|\rangle}^2_t
\end{split}
\end{equation*}
For $t$ large enough, the term $O({\langle |\varphi|\rangle}^2_t)$
can be absorbed into the right-hand side. Take $T_0$ to be the smallest $t$ for which this happens.
Then for all $t \geq T_0$,
\begin{equation*}
K \qbtg {\varphi}+K_t \sum_\nu \:{||\dze_\nu \potdg {\nu} \ze_\nu
\varphi^\nu||}^2_0 +O_t (||\varphi||^2_{-1}) \geq \half \, C \, t
{\langle |\varphi|\rangle}^2_t
\end{equation*}
Divide through by $\frac {C}{2}$ and call the new constants also
$K$ and $K_t$ respectively. Moreover, let $K'_t$ be a large enough
$t$ dependent constant such that for all $t \geq T_0,$
\begin{equation*}
\begin{split}
K \qbtg {\varphi}+K_t \sum_\nu \:{||\dze_\nu \potdg {\nu} \ze_\nu
\varphi^\nu||}^2_0 +K'_t \,||\varphi||^2_{-1}& \geq  t {\langle
|\varphi|\rangle}^2_t.
\end{split}
\end{equation*} \qed

\medskip
\medskip
\noindent Before we conclude this section, we will prove a lemma
that will be useful in the next section for establishing the
existence and regularity results for $\dbarb,$ namely that
$\qbtg{\, \cdot \,}$ controls the Sobolev $1$-norm of a form whose
Fourier transform is basically supported in $\co$ up to a smooth
error. First, however, let us show $\qbto{\, \cdot \,}$ controls
the $1$ norm of such a form, and then we will relate $\qbto{\,
\cdot \,}$ and $\qbtg{\, \cdot \,}.$
\smallskip
\newtheorem{coSob1}[lemma]{Lemma}
\begin{coSob1}
Let $\varphi$ be a $(0,1)$ form \label{Sob1loc} supported in $U'$
such that up to a smooth term, $\hat \varphi$ is supported in
$\co.$ There exist positive constants $C$ and $\Upsilon$
independent of $t$ for which
$$C \qbto{\varphi} + \Upsilon ||\varphi||_{\, 0}^2 \geq ||\varphi||_1^2.$$
\end{coSob1}
\noindent {\bf{Proof:}} Since $T$, $L_j$, and $\lbar{j}$ for $1
\leq j \leq n-1$ span the tangent space of $M$, we have to show
that $\qbto{\varphi}$ controls the $T$, $L$, and $\lbar{\,}$
derivatives modulo the square of the $L^2$ norm of $\varphi$ on
$\co$. We use the method in Lemma~\ref{qbtp} with $t = 0$ to
deduce that for some small positive $\eps'$
\begin{equation*}
\begin{split}
\qbto{\varphi} &\geq (1-\eps')\, \sum_{ij} ||\lbar{j}
\varphi_i||^2_0 + \sum_{ij} \Re \{(c_{ij} T \varphi_i,
\varphi_j)_0 \} + O(||\varphi||^2_0 ).
\end{split}
\end{equation*}
For some small positive $\eps$ $$\sum_{ij} \Re \{(c_{ij} T
\varphi_i, \varphi_j)_0 \} \geq -\eps \sum_i ||T \varphi_i||^2_0 +
O(||\varphi||^2_0 ).$$ Consider the sum
\begin{equation*}
\begin{split}
\sum_{ij} (T \varphi_i, T \varphi_j)_0 &= \sum_{ij} (T^* T
\varphi_i, \varphi_j)_0.
\end{split}
\end{equation*}
Since $\sigma (T) \, = \,  \xi_{2n-1}$, up to smooth terms
controlled by $O(||\varphi||^2_0 ) ,$ $\sigma (T^*) \, = \,
 \sum_{\alpha} \frac {1}{\alpha !} \: \partial_{\xi}^{\alpha}
D_x^{\alpha} ( {\overline{ \xi_{2n-1}}}) \, = \, \xi_{2n-1}$
because $ \xi_{2n-1}$ is real and does not depend on $x$, so
$D_x^{\alpha} (\overline { \xi_{2n-1}}) \, = \, 0$ for $\alpha
\neq 0$. By Plancherel's Theorem and the definition of the $\co$
region, the previous sum equals:
\begin{equation*}
\begin{split}
\sum_{ij} (T \varphi_i, T \varphi_j)_0 &= \sum_{ij} (T^* T
\varphi_i, \varphi_j)_0 = \sum_{ij} (|\xi_{2n-1}|^2
\widehat{\varphi_i}, \widehat{\varphi_j})_0\\&\leq C \sum_{ij}
(\sum_{l = 1}^{2n-2} |\xi_l|^2 \widehat{\varphi_i},
\widehat{\varphi_j})_0= C \sum_{ij} (|\xi'|^2 \widehat{\varphi_i},
\widehat{\varphi_j})_0,
\end{split}
\end{equation*}
where $C$ is a positive constant independent of $t$. Since $\xi'$
is dual to the holomorphic part of the tangent bundle
$T^{1,0}(M)\oplus T^{0,1} (M)$ spanned by $L_1, \dots,
L_{n-1},\lbar{1}, \dots, \lbar{n-1}$ and $\lbar{i}^* \, = \, -
L_i-f_i,$ there exists some $t$ independent constant $C'$ so
$$\sum_i ||T \varphi_i||^2_0 \leq C'\sum_{ij} (||\lbar{j}
\varphi_i||_0^2+||\lbar{j}^*
\varphi_i||_0^2)+O(||\varphi||_0^2).$$ In other words, $L$ and
$\bar L$ derivatives control $T$ derivatives up to
$O(||\varphi||^2_0 ) $ terms. This implies that
$$ \sum_{ij} \Re \{(c_{ij} T \varphi_i, \varphi_j)_0 \} \geq
-\eps\, C'\sum_{ij} ||\lbar{j} \varphi_i||_0^2 -\eps\, C'\sum_{ij}
||\lbar{j}^* \varphi_i||_0^2+ O(||\varphi||^2_0 ).$$ Just as in
the proof of the Lemma~\ref{qbtm},
\begin{equation*}
\begin{split}
||\lbar{j}^* \varphi_i||_0^2 &= {(\lbar{j}^* \lbar{j} \varphi_i,
\varphi_i)}_0 +{( [\lbar{j},\lbar{j}^*] \varphi_i,
\varphi_i)}_0\\&=||\lbar{j} \varphi_i||_0^2+{( [L_j,\lbar{j}]
\varphi_i, \varphi_i)}_0+O(||\varphi||^2_0 )\\&\leq
(1+2\eps)\,||\lbar{j} \varphi_i||_0^2+\Re \{{( c_{jj}T \varphi_i,
\varphi_i)}_0\}+O(||\varphi||^2_0 ),
\end{split}
\end{equation*}
where we have also used Equation~\ref{primitivebracket} from
Section~\ref{CRpsh} along with $\lbar{j}^* \, = \, - L_j-f_j.$
Thus, $$ \sum_{ij} \Re \{(c_{ij} T \varphi_i, \varphi_j)_0 \} \geq
-2\, C' \,\eps (1+\eps) \, \sum_{ij} ||\lbar{j} \varphi_i||_0^2
-\eps\, C'\sum_{ij} \Re \{{( c_{jj} T \varphi_i, \varphi_i)}_0\}+
O(||\varphi||^2_0 ).$$ By the same method we derive
\begin{equation*}
\begin{split}
\sum_{ij}\Re \{(c_{jj} T \varphi_i, \varphi_i)_0 \} &\leq \eps\,
C'\sum_{ij} ||\lbar{j} \varphi_i||_0^2 +\eps\, C'\sum_{ij}
||\lbar{j}^* \varphi_i||_0^2+ O(||\varphi||^2_0 )\\&\leq 2 \,
C'\,\eps(1+\eps)\,\sum_{ij} ||\lbar{j} \varphi_i||_0^2 +\eps\,
C'\sum_{ij}\Re \{(c_{jj} T \varphi_i, \varphi_i)_0 \}+
O(||\varphi||^2_0 ),
\end{split}
\end{equation*}
which means
$$\sum_{ij} \Re \{(c_{ij} T \varphi_i, \varphi_j)_0 \} \leq
\frac{2\, C' \,\eps (1+\eps)}{1-\eps\, C'} \, \sum_{ij} ||\lbar{j}
\varphi_i||_0^2 +O(||\varphi||^2_0 ).$$ Therefore, if $\eps$ is
chosen to be small enough, there exists a very small constant
$\eps"$ such that
$$\qbto{\varphi} \geq (1-\eps")\, \sum_{ij} ||\lbar{j}
\varphi_i||^2_0  + O(||\varphi||^2_0 ),$$ which means $\qbto {\,
\cdot \,}$ controls all the $\bar L$ derivatives up to
$O(||\varphi||^2_0 ) $ terms. Similarly, one proves that
$$\qbto{\varphi} \geq (1-\eps")\, \sum_{ij} ||\lbar{j}^*
\varphi_i||^2_0  + O(||\varphi||^2_0 ),$$ hence that $\qbto {\,
\cdot \,}$ also controls all the $L$ derivatives up to
$O(||\varphi||^2_0 ) $ terms. \qed

\medskip \noindent Before showing how $\qbto {\, \cdot \,}$ relates to $\qbtg
{\, \cdot \,},$ let us relate $\ad$ with $\adt$ using the
operators $F_t$ and $G_t,$ defined in Corollary~\ref{equivop}:
\begin{equation*}
\begin{split}
\tnorm{\varphi}{\dbarb \phi} &= (F_t \varphi, \dbarb \phi)_0 =
(F_t \ad \varphi, \phi)_0 + ([\ad, F_t] \varphi, \phi)_0\\&=
\tnorm{\ad \varphi}{\phi}+\tnorm{G_t [\ad, F_t] \varphi}{\phi}
\end{split}
\end{equation*}
We conclude that
\begin{equation}
\adt = \ad + G_t [\ad, F_t] \label{adtusingad}.
\end{equation}

\smallskip
\newtheorem{coSob1global}[lemma]{Lemma}
\begin{coSob1global}
Let $\varphi$ be a $(0,1)$ form \label{Sob1global} supported in
$U_\nu$ for some $\nu$ such that up to a smooth term, $\hat
\varphi$ is supported in $\cod_\nu.$ There exist positive
constants $C'>1$ and $\Upsilon'$ independent of $t$ for which
$$C' \qbtgs{\varphi}{G_t \varphi} + \Upsilon' ||\varphi||_{\, 0}^2 \geq ||\varphi||_1^2.$$
\end{coSob1global}
\noindent {\bf{Proof:}} Keeping in mind that $F_t$ and $G_t$ are both self-adjoint and
inverses of each other and that $\sum_\nu \ze^2_\nu \, = \,1,$ we compute $\qbtg{\, \cdot \,}$ in terms
of $\qbto{\, \cdot \,}$ using Equation~\ref{adtusingad}:
\begin{equation*}
\begin{split}
\qbtgs{\varphi}{\phi} &= \tnorm{\dbarb \varphi}{\dbarb \phi} +
\tnorm {\adt \varphi}{\adt \phi} \\&= (\dbarb \varphi, F_t \dbarb
\phi)_0 + \sum_\nu (\ze_\nu(\ad + G_t [\ad, F_t])\varphi, \ze_\nu F_t
(\ad + G_t [\ad, F_t]) \phi)_0\\&= (\dbarb \varphi, \dbarb F_t
\phi)_0 + (\dbarb \varphi, [F_t, \dbarb] \phi)_0 + (\ad \varphi,
\ad F_t \phi)_0 + (\ad \varphi, [F_t, \ad] \phi)_0\\&\quad+ \sum_\nu
(\ze_\nu G_t [\ad, F_t] \varphi, \ze_\nu F_t \ad
\phi)_0+ \sum_\nu (\ze_\nu G_t [\ad, F_t] \varphi, \ze_\nu [\ad, F_t]
\phi)_0+(\ad \varphi, [\ad,F_t] \phi)_0\\&= \qbtos{\varphi}{F_t
\phi}+ (\dbarb \varphi, [F_t, \dbarb] \phi)_0 + (\ad \varphi,
[F_t, \ad] \phi)_0+ \sum_\nu (\ze_\nu G_t [\ad, F_t] \varphi, \ze_\nu
F_t \ad \phi)_0\\&\quad+ \sum_\nu (\ze_\nu G_t [\ad, F_t] \varphi,
\ze_\nu [\ad, F_t] \phi)_0+(\ad \varphi, [\ad,F_t] \phi)_0
\end{split}
\end{equation*}
Take $\phi \, = \, G_t \varphi$ in the previous expression to
conclude that
\begin{equation*}
\begin{split}
\qbtgs{\varphi}{G_t \varphi} &= \qbto{\varphi}+ (\dbarb \varphi,
[F_t, \dbarb]G_t \varphi)_0 + (\ad \varphi, [F_t, \ad] G_t
\varphi)_0\\&\quad+ \sum_\nu (\ze_\nu G_t [\ad, F_t] \varphi, \ze_\nu
F_t \ad G_t \varphi)_0+ \sum_\nu (\ze_\nu G_t [\ad, F_t] \varphi,
\ze_\nu [\ad, F_t] G_t \varphi)_0\\&\quad+(\ad \varphi,
[\ad,F_t] G_t \varphi)_0.
\end{split}
\end{equation*}
Since both $F_t$ and $G_t$ are of order zero, $$|(\dbarb \varphi,
[F_t, \dbarb]G_t \varphi)_0| \leq \eps ||\dbarb \varphi||_0^2 +
\Gamma ||\varphi||^2_0,$$ $$|(\ad \varphi, [F_t, \ad] G_t
\varphi)_0| \leq  \eps ||\ad \varphi||_0^2 + \Gamma
||\varphi||^2_0,$$ $$\Big|\sum_\nu \, (\ze_\nu G_t [\ad, F_t] \varphi, \ze_\nu
[\ad, F_t] G_t \varphi)_0\Big| \leq \Gamma ||\varphi||^2_0,$$ $$|(\ad
\varphi, [\ad,F_t] G_t \varphi)_0| \leq \eps ||\ad \varphi||_0^2 +
\Gamma ||\varphi||^2_0,$$ and
\begin{equation*}
\begin{split}
\Big|\sum_\nu \, (\ze_\nu G_t [\ad, F_t] \varphi, \ze_\nu F_t \ad G_t
\varphi)_0\Big| &\leq \Big|\sum_\nu \, (\ze_\nu G_t [\ad, F_t] \varphi, \ze_\nu F_t [\ad
,G_t ]\varphi)_0\Big|\\&\quad+ \Big|\sum_\nu \, (\ze^2_\nu G_t [\ad, F_t] \varphi,   \ad
\varphi)_0\Big| \\&\leq \eps ||\ad \varphi||_0^2 + \Gamma
||\varphi||^2_0,
\end{split}
\end{equation*}
for some small $\eps>0$ and some large positive $\Gamma.$ This
means that for some small $\eps'>0$ and some large $\Gamma'>0,$
$$\qbtgs{\varphi}{G_t \varphi} + \eps' \qbto{\varphi} + \Gamma'
||\varphi||^2_0\geq \qbto{\varphi}, $$ so there exist positive
constants $C'>1$ and $\Gamma'$ independent of $t$ such that
$$C'\qbtgs{\varphi}{G_t \varphi}+ \Gamma'
||\varphi||^2_0\geq \qbto{\varphi}.$$ Now we can apply
Lemma~\ref{Sob1loc} with $U_\nu$ as $U'$ and $\cod_\nu$ as $\co$
to see that $$C \qbto{\varphi} + \Upsilon ||\varphi||_{\, 0}^2
\geq ||\varphi||_1^2.$$ Putting these two together and enlarging
$C',$ it follows there exists some large positive $\Upsilon'$
independent of $t$ such that $$C' \qbtgs{\varphi}{G_t \varphi} +
\Upsilon' ||\varphi||_{\, 0}^2 \geq ||\varphi||_1^2.$$ \qed

\bigskip
\section{Existence and Regularity Results for the \dbarb\, Operator}
\label{regdbar}

\medskip
\medskip
\noindent We start by using Proposition~\ref{global} to prove an
estimate similar to the one in \cite{Kohnrange}.
\medskip
\newtheorem{proposition}{Proposition}[section]
\begin{proposition}
Let $M$ be a compact, orientable, weakly pseudoconvex
\label{estimate} CR manifold of dimension at least $5$, endowed
with a strongly CR plurisubharmonic function $\lambda$ as
described in Section~\ref{CRpsh}. Given some $s>0,$ a $(0,1)$ form
$\alpha \in H^s (M),$ and a $(0,1)$ form $\varphi$ supported on
$M$ satisfying $\varphi \in Dom(\dbarb) \cap Dom(\ad)$ and
$$\qbtgs{\varphi}{\phi} = \tnorm{\alpha}{\phi},$$ for all $\phi
\in Dom(\dbarb) \cap Dom(\ad),$ there exists a positive number
$T_s$ such that for any $t \geq T_s,$ $\varphi \in H^s$ and
\begin{equation}
\label{estimateexpr} ||\varphi||_s \leq C_t \, (||\alpha||_s +
||\varphi||_{\, 0}),
\end{equation}
for some $t$ dependent constant $C_t.$
\end{proposition}

\medskip
\noindent The main ingredient of the proof is the following lemma:
\smallskip
\newtheorem{qbtcomm}[proposition]{Lemma}
\begin{qbtcomm}
\label{qbtcommlemma} Let $P$ be a pseudodifferential operator of
order $p$ and $\varphi$ and $\phi$ be two $(0,1)$ forms such that
$\varphi, \phi \in Dom(\dbarb) \cap Dom(\ad),$ then
\begin{equation}
\begin{split}
\qbtgs{P \varphi}{\phi} &= \qbtgs{\varphi}{P^{*, \, t} \phi}+ 2\,
\tnorm{[\dbarb,P] \varphi}{\dbarb \phi}+2\,\tnorm{[\adt,P]
\varphi}{\adt \phi}\\&\quad+\tnorm{[[\dbarb,P],\adt]
\varphi}{\phi}+\tnorm{[[\adt,P],\dbarb] \varphi}{\phi},
\label{qbtcommexpr}
\end{split}
\end{equation}
where $P^{*, \, t}$ is the adjoint of $P$ with respect to
${\langle| \, \cdot \,|\rangle}_t.$ Moreover, if $P$ is inverse
zero order $t$ dependent,
\begin{enumerate}
\item[(i)] for each small $\eps >0,$ there exist constants $C$ and
$C_t$ such that $$\qbtgs{P \varphi}{\phi} \leq
\qbtgs{\varphi}{P^{*, \, t} \phi}+ \eps \, \qbtg{\phi} + C \,
\sqtnorm{\Lambda^p \varphi}+C_t\, \sqtnorm{\Lambda^{p-1}
\varphi}+\eps \, \sqtnorm{\phi};$$
\item[(ii)]$$\qbtg{P \varphi} \leq
C \, \qbtgs{\varphi}{P^{*, \, t} P \varphi} + C \, \sqtnorm{\Lambda^p
\varphi}+  C_t\, \sqtnorm{\Lambda^{p-1} \varphi};$$
\item[(iii)] and
\begin{equation*}
\begin{split}
\qbtgs{P \varphi}{R \varphi} &\leq \qbtgs{\varphi}{P^{*, \, t}R
\varphi}+ C \, \qbtgs{\varphi}{R^{*, \, t}R \varphi} + C \,
\sqtnorm{\Lambda^p \varphi}+  C_t\, \sqtnorm{\Lambda^{p-1}
\varphi}\\&\quad + C \, \sqtnorm{\Lambda^r \varphi}+  C_t\,
\sqtnorm{\Lambda^{r-1} \varphi},
\end{split}
\end{equation*}
where $R$ is a pseudodifferential operator of order $r.$
\end{enumerate}
\end{qbtcomm}
\smallskip
\noindent {\bf Proof:} $$\qbtgs{P \varphi}{\phi} = \tnorm{\dbarb P
\varphi}{\dbarb \phi}+\tnorm{\adt P \varphi}{\adt \phi}$$ Now,
\begin{equation*}
\begin{split}
\tnorm{\dbarb P \varphi}{\dbarb \phi}&= \tnorm{ P\dbarb
\varphi}{\dbarb \phi}+\tnorm{[\dbarb, P] \varphi}{\dbarb
\phi}\\&=\tnorm{ \dbarb \varphi}{P^{*, \, t}\dbarb
\phi}+\tnorm{[\dbarb, P] \varphi}{\dbarb \phi}\\&=\tnorm{ \dbarb
\varphi}{\dbarb P^{*, \, t} \phi}+\tnorm{ \dbarb \varphi}{[P^{*,
\, t},\dbarb] \phi}+\tnorm{[\dbarb, P] \varphi}{\dbarb
\phi}\\&=\tnorm{ \dbarb \varphi}{\dbarb P^{*, \, t}
\phi}+\tnorm{[\dbarb, P] \varphi}{\dbarb \phi}\\&\quad+\tnorm{
[P^{*, \, t},\dbarb]^{*, \, t} \varphi}{\adt \phi}+\tnorm{[[P^{*,
\, t},\dbarb]^{*, \, t}, \dbarb] \varphi}{ \phi}.
\end{split}
\end{equation*}
Since $$[P^{*, \, t},\dbarb]^{*, \, t}=(P^{*, \, t}\dbarb)^{*, \,
t}-(\dbarb P^{*, \, t})^{*, \, t} = \adt P - P \adt = [\adt ,P],$$
\smallskip
\begin{equation*}
\begin{split}
\tnorm{\dbarb P \varphi}{\dbarb \phi}&=\tnorm{ \dbarb
\varphi}{\dbarb P^{*, \, t} \phi}+\tnorm{[\dbarb, P]
\varphi}{\dbarb \phi}+\tnorm{ [\adt ,P] \varphi}{\adt
\phi}\\&\quad+\tnorm{[[\adt ,P], \dbarb] \varphi}{ \phi}.
\end{split}
\end{equation*}
\smallskip
Equation~\ref{qbtcommexpr} thus follows. $P$ is inverse zero order
$t$ dependent, and the top order term of $\adt$ is independent of
$t,$ so (i) follows as well. To show (ii) apply (i) with $\phi \,
= \, P \varphi,$ and for (iii) use both (i) and (ii). \qed

\medskip
\noindent {\bf Proof of Proposition~\ref{estimate}:} First, we
prove Expression~\ref{estimateexpr} as an a priori estimate,
assuming $\varphi \in H^s.$ By Proposition~\ref{global},
\begin{equation*}
t {\langle |\Lambda^s \varphi|\rangle}^2_t \leq C \qbtg {\Lambda^s
\varphi}+C_t (\: \sum_\nu \:{||\dze_\nu \potdg {\nu} \ze_\nu
\Lambda^s \varphi^\nu||}^2_0 +||\varphi||^2_{s-1}).
\end{equation*}
Since
\begin{equation*}
\begin{split}
\sum_\nu \:{||\dze_\nu \potdg {\nu} \ze_\nu \Lambda^s
\varphi^\nu||}^2_0 &\leq 2 \sum_\nu \:{||\dze_\nu \Lambda^{-1}
\potdg {\nu} \ze_\nu \Lambda^s \varphi^\nu||}^2_1 + C \,
||\varphi||^2_{s-1}\\&\leq C \, (\: \sum_\nu \: \qbtgs{\dze_\nu
\Lambda^{-1} \potdg {\nu} \ze_\nu \Lambda^s \varphi^\nu}{G_t
\dze_\nu \Lambda^{-1} \potdg {\nu} \ze_\nu \Lambda^s
\varphi^\nu}+||\varphi||^2_{s-1})
\end{split}
\end{equation*}
by Lemma~\ref{Sob1global} proven in the previous section and since
$G_t$ is not inverse zero order $t$ dependent, we can only apply
Equation~\ref{qbtcommexpr} twice, not part (iii) of
Lemma~\ref{qbtcommlemma}, to obtain that $$\sum_\nu \:{||\dze_\nu
\potdg {\nu} \ze_\nu \Lambda^s \varphi^\nu||}^2_0 \leq C_t \,
(||\alpha||^2_{s-1}+||\varphi||^2_{s-1}).$$ Part (ii) of
Lemma~\ref{qbtcommlemma} applied to $\qbtg {\Lambda^s \varphi}$
together with the property of $\varphi$ that
$$\qbtgs{\varphi}{\phi} = \tnorm{\varphi}{\phi},$$ for all $\phi \in
Dom(\dbarb) \cap Dom(\ad),$ imply that
\begin{equation*}
\begin{split}
\qbtg {\Lambda^s \varphi} &\leq \tnorm{\Lambda^s \alpha}{\Lambda^s
\varphi} + C \, \sqtnorm{\Lambda^s \varphi} + C_t \,
\sqtnorm{\Lambda^{s-1} \varphi} \\&\leq C_t \, (||\alpha||^2_s +
||\varphi||^2_{s-1}) + (C+ \eps) \, \sqtnorm{\Lambda^s \varphi},
\end{split}
\end{equation*}
for some small $\eps >0.$ Putting everything together, we see that
$$t {\langle |\Lambda^s \varphi|\rangle}^2_t \leq  C_t \, (||\alpha||^2_s +
||\varphi||^2_{s-1}) + (C+ \eps) \, \sqtnorm{\Lambda^s \varphi}.$$
In other words, there exists some $T_s \geq 1$ such that for all
$t \geq T_s,$ $$||\varphi||_s^2\leq  C_t \, (||\alpha||^2_s +
||\varphi||^2_{s-1}).$$ The rest of the argument is done by
induction. To establish the base case, notice that for any $0 < s
\leq 1,$ $$||\varphi||_s^2\leq  C_t \, (||\alpha||^2_s +
||\varphi||^2_{s-1}) \leq  C_t \, (||\alpha||^2_s +
||\varphi||^2_0).$$ For any $s >1,$ the induction step is that
$$||\varphi||_{s-1}^2\leq  C_t \, (||\alpha||^2_{s-1} +
||\varphi||^2_0)$$ which implies immediately that
$$||\varphi||_s^2\leq  C_t \, (||\alpha||^2_s +
||\varphi||^2_0),$$ for all $t \geq T_s.$

\smallskip
\noindent Second, we show that indeed $\varphi \in H^s.$ We define
a family of mollifiers in the $\xi$ space as follows: Let $\rho$
be a smooth function on $[0, \infty)$ satisfying
\begin{enumerate}
\item[(i)] $\rho(0) \, = \, 1;$
\item[(ii)] $\exists \, \tau>0$ such that $\rho(y) \, = \, 1,$ for
all $0<y < \tau;$
\item[(iii)] $\int_0^\infty \rho(y)\, dy \, = \, 1;$
\item[(iv)] $\rho$ has compact support.
\end{enumerate}
For each $0< \delta <1,$ define $$\rho_\delta(\xi) = \rho(\delta
|\xi|).$$ Now let $S_\delta$ be the pseudodifferential operator
with symbol $\rho_\delta.$ The operators $S_\delta$ are smoothing,
but of order zero with respect to $\delta^{-1}.$ Notice also that
as $\delta \rightarrow 0,$ $S_\delta \varphi \rightarrow \varphi$
in $L^2.$ If we then proved that there existed some $\Delta <1$
such that $||S_\delta \varphi||_s$ was bounded independently of
$\delta$ for all $\delta<\Delta,$ this would imply by
the Monotone Convergence Theorem that $\varphi \in H^s.$ This proof
follows the same line as the one of the a priori estimate: By
Proposition~\ref{global},
\begin{equation*}
t {\langle |\Lambda^s S_\delta \varphi|\rangle}^2_t \leq C \qbtg
{\Lambda^s S_\delta\varphi}+C_t (\: \sum_\nu \:{||\dze_\nu \potdg
{\nu} \ze_\nu \Lambda^s S_\delta\varphi^\nu||}^2_0
+||S_\delta\varphi||^2_{s-1}),
\end{equation*}
for all $t \geq T_0.$ Then we apply Equation~\ref{qbtcommexpr} to
obtain that
\begin{equation*}
\begin{split}
\qbtg{\Lambda^s S_\delta \varphi} &= \tnorm {\Lambda^s S_\delta
\alpha}{\Lambda^s S_\delta \varphi}+ 2\, \tnorm{[\dbarb,\Lambda^s
S_\delta] \varphi}{\dbarb \Lambda^s S_\delta
\varphi}\\&\quad+2\,\tnorm{[\adt,\Lambda^s S_\delta] \varphi}{\adt
\Lambda^s S_\delta \varphi}+\tnorm{[[\dbarb,\Lambda^s
S_\delta],\adt] \varphi}{\Lambda^s S_\delta \varphi
}\\&\quad+\tnorm{[[\adt,\Lambda^s S_\delta],\dbarb]
\varphi}{\Lambda^s S_\delta \varphi}.
\end{split}
\end{equation*}
Notice that $$[\dbarb,\Lambda^s S_\delta] = [\dbarb, \Lambda^s]
S_\delta + \Lambda^s[\dbarb, S_\delta],$$ $$[[\dbarb,\Lambda^s
S_\delta],\adt] = [[\dbarb, \Lambda^s], \adt] S_\delta + [\dbarb,
\Lambda^s][S_\delta, \adt] + [\Lambda^s, \adt][\dbarb,
S_\delta]+\Lambda^s [[\dbarb, S_\delta], \adt],$$ and similarly
for $[\adt,\Lambda^s S_\delta]$ and $[[\adt,\Lambda^s
S_\delta],\dbarb].$ Since the symbol of $S_\delta$ has compact
support and is identically equal to $1$ for all $|\xi| <
\frac{\tau}{\delta},$ the commutators $[\dbarb, S_\delta]$ and
$[S_\delta, \adt]$ have symbols that go to zero as $\delta
\rightarrow 0.$ The top order terms of $\adt$ are independent of
$t$ so it follows that there exists some $\Delta <1$ such that for
all $\delta < \Delta$ $$\qbtg{\Lambda^s S_\delta \varphi} \leq C
\, \sqtnorm{\Lambda^s S_\delta \alpha} + C' \,\sqtnorm{\Lambda^s
S_\delta \varphi}+ C_t \, ||S_\delta \varphi||^2_{s-1}.$$ Also, by
the same procedure as above, namely Lemma~\ref{Sob1global} and
Equation~\ref{qbtcommexpr}, $$\sum_\nu \:{||\dze_\nu \potdg {\nu}
\ze_\nu \Lambda^s S_\delta\varphi^\nu||}^2_0 \leq C_t \, ( \,
||S_\delta \alpha||^2_{s-1} + ||S_\delta \varphi||^2_{s-1}).$$
Altogether, there exists some $T" \geq 1$ such that for all $t
\geq T"$ and all $\delta < \Delta,$
\begin{equation*}
||\Lambda^s S_\delta \varphi||^2_s \leq C_t \, ( \, ||S_\delta \alpha||^2_s + ||S_\delta
\varphi||^2_{s-1}) \leq C_t \, ( \, || \alpha||^2_s + ||S_\delta \varphi||^2_{s-1}).
\end{equation*}
The rest follows by induction. For $s \, = \, 0,$ $||S_\delta \varphi||_0 \leq ||\varphi||_0$ by construction.
For $0 < s \leq 1,$
\begin{equation*}
||\Lambda^s S_\delta \varphi||^2_s \leq C_t \, ( \, || \alpha||^2_s + ||S_\delta \varphi||^2_{s-1})
\leq C_t \, ( \, || \alpha||^2_s + || \varphi||^2_0).
\end{equation*}
By induction, $||S_\delta \varphi||_{s-1} \leq C_t,$ for some $C_t$ independent of $\delta,$ then
$$||\Lambda^s S_\delta \varphi||^2_s \leq C_t \, ( \, || \alpha||^2_s +C'_t)\leq C"_t,$$
for some $C"_t$ independent of $\delta$ and all $t \geq T"$ and
$\delta < \Delta.$ \qed

\bigskip
\bigskip
\noindent Next, for each $t$ we define the corresponding harmonic space to be $$\harmon = \{ \varphi \in Dom(\dbarb)
\cap Dom(\ad) \; \big| \; \dbarb \varphi = 0 \: and \: \adt \varphi=0 \}$$ $$\Updownarrow$$
$$\harmon = \{ \varphi \in Dom(\dbarb) \cap Dom(\ad) \; \big| \; \qbtg{\varphi}=0 \},$$ and we prove the following
lemma similar to Lemma $5.7$ of \cite{Kohnglobalreg}:

\smallskip
\newtheorem{harmonic}[proposition]{Lemma}
\begin{harmonic}
Let $M$ be a compact, orientable, weakly pseudoconvex CR manifold of dimension
at least $5$, endowed with a strongly CR plurisubharmonic function
$\lambda$ as described in Section~\ref{CRpsh}.\label{harmsp} There exists some $T_0 >1$
such that for each $t \geq T_0$ the space $\harmon$ is finite dimensional and some $t$ independent positive
constant $C$ such that all $\varphi \in Dom(\dbarb) \cap Dom(\ad)$ satisfying $\varphi \perp \harmon$
with respect to the ${\langle | \, \cdot \, | \rangle}_t$ norm also satisfy
\begin{equation}
\label{harmperp} \sqtnorm{\varphi} \leq C \, \qbtg{\varphi}.
\end{equation}
\end{harmonic}
\smallskip
\noindent {\bf Proof:} By Proposition~\ref{global}, for all $t \geq T_0$ and $\varphi \in \harmon,$
\begin{equation*}
\begin{split}
t {\langle | \varphi|\rangle}^2_t &\leq C \qbtg {
\varphi}+C_t (\: \sum_\nu \:{||\dze_\nu \potdg {\nu} \ze_\nu
 \varphi^\nu||}^2_0 +||\varphi||^2_{-1})\\&\leq C_t (\: \sum_\nu \:{||\dze_\nu \potdg {\nu} \ze_\nu
 \varphi^\nu||}^2_0 +||\varphi||^2_{-1}).
\end{split}
\end{equation*}
Just as in the proof of Proposition~\ref{estimate}, by Lemma~\ref{Sob1global},
\begin{equation*}
\begin{split}
\sum_\nu \:{||\dze_\nu \potdg {\nu} \ze_\nu
\varphi^\nu||}^2_0 &\leq C \, (\: \sum_\nu \: \qbtgs{\dze_\nu
\Lambda^{-1} \potdg {\nu} \ze_\nu  \varphi^\nu}{G_t
\dze_\nu \Lambda^{-1} \potdg {\nu} \ze_\nu
\varphi^\nu}+||\varphi||^2_{-1})\\&=  C \, (\: \sum_\nu \: \tnorm{\dbarb \dze_\nu
\Lambda^{-1} \potdg {\nu} \ze_\nu  \varphi^\nu}{\dbarb G_t
\dze_\nu \Lambda^{-1} \potdg {\nu} \ze_\nu
\varphi^\nu}\\&+ \sum_\nu \: \tnorm{\adt \dze_\nu
\Lambda^{-1} \potdg {\nu} \ze_\nu  \varphi^\nu}{\adt G_t
\dze_\nu \Lambda^{-1} \potdg {\nu} \ze_\nu
\varphi^\nu} +||\varphi||^2_{-1})\\&\leq C_t \,||\varphi||^2_{-1},
\end{split}
\end{equation*}
since the commutators of $\dze_\nu \Lambda^{-1} \potdg {\nu} \ze_\nu$ and $G_t
\dze_\nu \Lambda^{-1} \potdg {\nu} \ze_\nu $ with $\dbarb$ and $\adt$ are of order
$-1$ and $\varphi \in \harmon.$ This implies that for all $t \geq T_0$ and $\varphi \in \harmon,$
$$t {\langle | \varphi|\rangle}^2_t\leq C_t \,||\varphi||^2_{-1}.$$ In other words,
the identity map from $H^{-1} \big|_{\harmon}$ to $L^2$ is bounded. By Rellich's theorem,
$L^2$ is compactly embedded in $H^{-1}$ via the identity map, so the same will be true of $L^2 \big|_{\harmon}.$
The composition of a bounded operator with a compact one is compact, which implies the identity
map from $L^2 \big|_{\harmon}$ to $L^2$ is compact. It follows that the unit sphere of $\harmon$ is compact,
hence that $\harmon$ is finite dimensional for each $t \geq T_0.$

In order to prove Equation~\ref{harmperp}, we assume it does not hold, namely that
there exists a sequence $\{\varphi_k\}_k$ such that ${\langle | \varphi_k |\rangle}_t \, = \, 1,$
$\varphi_k \perp \harmon$ with respect to the ${\langle | \, \cdot \, | \rangle}_t$ norm, and
\begin{equation}
\label{harmperpseqprop1} {\langle | \varphi_k |\rangle}^2_t \geq k \, \qbtg{\varphi_k}.
\end{equation}
From this last expression, one obtains by the same procedure as in the first part of the proof of this lemma that
\begin{equation}
\label{harmperpseqprop2} {\langle | \varphi_k |\rangle}^2_t \leq C_t \, ||\varphi_k||^2_{-1}
\end{equation}
for $k$ large enough. Just as above, Expression~\ref{harmperpseqprop2} then implies a subsequence of $\{ \varphi_k\}_k$
converges in $L^2$ in the ${\langle | \, \cdot \, |\rangle}_t$ norm. Moreover, notice that
$[\qbtg{\, \cdot \,} +{\langle | \, \cdot \, |\rangle}^2_t\,]^\half$ is a norm, so
Expressions~\ref{harmperpseqprop1} and ~\ref{harmperpseqprop2} together imply that a subsequence of $\{ \varphi_k\}_k$
converges in $L^2$ in this norm as well. Thus, we can find a subsequence $\{ \varphi_{k_j}\}$ that
converges in both norms to some $\varphi.$ Since $\varphi$ is the limit in the
${\langle | \, \cdot \, |\rangle}_t$ norm, it must satisfy ${\langle | \varphi |\rangle}_t \, = \, 1$
and $\varphi \perp \harmon.$ The fact that it is also the limit in the
$[\qbtg{\, \cdot \,} +{\langle | \, \cdot \, |\rangle}^2_t\,]^\half$ norm along with Expression~\ref{harmperpseqprop1}
implies that $\varphi \in \harmon.$ We thus have the contradiction we sought. \qed

\medskip
\noindent {\bf Note:} In order to avoid confusion, we shall specify that from now on we assume all orthogonality
and all projections to be with respect to the ${\langle | \, \cdot \, |\rangle}_t$ norm.

\bigskip
\noindent Let $$\harperp = \{ \varphi \in L^2 \; \big| \;
{\langle | \varphi, \phi | \rangle}_t=0 \; \forall \, \phi \in \harmon \},$$ and on $\harperp$ we define
$$\boxbt = \dbarb \, \adt + \adt \, \dbarb.$$ Notice that
$$Dom(\boxbt) = \{ \varphi \in L^2 \; \big| \; \varphi \in Dom(\dbarb) \cap Dom(\ad), \: \dbarb \varphi \in Dom(\ad),
\: \ad \varphi \in Dom(\dbarb)\},$$ since the equivalence of the $L^2$ norm and ${\langle | \, \cdot | \rangle}_t$
implies that $Dom(\ad) \, = \, Dom(\adt).$ We claim that $\boxbt$ is self adjoint.
Given that $\dbarb$ is a closed, densely defined operator, it is a standard fact
in functional analysis that both $\dbarb \, \adt$ and $\adt \, \dbarb$ are self-adjoint.
Still, the sum of two self-adjoint operators may not be self-adjoint. In this case, however,
$\dbarb^2 \, = \, 0$ implies that $\ran(\dbarb) \subset \nul(\dbarb)$ which along with another elementary
fact that $\ranperp(\dbarb) \, = \, \nul(\adt)$ implies that $\nulperp(\dbarb) \subset \nul(\adt).$ On \harperp,
$\nul(\boxbt) \, = \, 0,$ hence it follows that if we decompose each $\varphi \in Dom(\boxbt)$ into the orthogonal
projection onto $\nul(\dbarb)$ which we shall call $\varphi_1$ and the rest which we shall call $\varphi_2,$ then
$\varphi_1 \in Dom( \dbarb \, \adt) \cap \nul(\dbarb)$ so $\boxbt \varphi_1 \, = \, \dbarb \, \adt \varphi_1$ and
$\varphi_2 \in Dom( \adt \, \dbarb) \cap \nul(\adt)$ so $\boxbt \varphi_2 \,= \, \adt \, \dbarb \varphi_2.$
The claim follows.
Next, we shall use following version of the Friedrichs' Extension Theorem
from \cite{FollandKohn} or \cite{Krantz} to define $\boxbt$ on a larger subset of
$\harperp$ than $Dom(\boxbt):$

\smallskip
\newtheorem{Friedrichslemma}[proposition]{Lemma}
\begin{Friedrichslemma}
\label{Friedrichs} Let $H$ be a Hilbert space equipped with the inner product $(\, \cdot \, , \, \cdot \,)$ and
corresponding norm $|| \, \cdot \,||$ and $Q$ a positive definite Hermitian form defined
on a dense $D \subset H$ satisfying $$|| \varphi||^2 \leq C \, Q(\varphi, \varphi)$$ for all $\varphi \in D.$ $D$ and $Q$
are such that $D$ is a Hilbert space under the inner product $Q(\, \cdot \, , \, \cdot \,).$ Then there exists a unique
self-adjoint operator $F$ with $Dom(F) \subseteq D$ satisfying $$Q(\varphi, \phi) = (F \varphi, \phi),$$ for all
$\varphi \in Dom(F)$ and $\phi \in D.$ $F$ is called the Friedrichs representative.
\end{Friedrichslemma}

\smallskip
\noindent {\bf Proof:} For each $\alpha \in H,$ we define a linear functional $\Omega_\alpha$ on $D$ given by
$$\Omega_\alpha \, : \, D \ni \phi \rightarrow (\alpha, \phi).$$ Since
\begin{equation*}
\begin{split}
|\Omega_\alpha (\phi)| &\leq ||\alpha|| \: ||\phi|| \\&\leq C \, ||\alpha|| \: Q(\phi, \phi)^\half,
\end{split}
\end{equation*}
by the Riesz Representation Theorem there exists some $\varphi_\alpha \in D$ such that
$$Q(\varphi_\alpha, \phi) = (\alpha, \phi).$$ In other words, that there exists an operator
$$T \, : \, H \ni \alpha \rightarrow \varphi_\alpha \in D.$$ The inverse of $T$ is the operator $F$ we seek. Let us
show $T$ has an inverse. First, notice that $T$ is bounded:
\begin{equation*}
\begin{split}
||T \alpha||^2 &\leq C \, Q(T\alpha, T \alpha) \\&= C \, (\alpha, T \alpha)\\&\leq C \, ||\alpha|| \: ||T \alpha||
\end{split}
\end{equation*}
which implies $||T \alpha|| \leq C \, ||\alpha||.$ If $T \alpha \, = \, 0,$ $(\alpha, \phi) \, = \,
Q (T \alpha, \phi) \, = \, 0$ for all $\phi \in D,$ but $D$ is dense in $H,$ so $\alpha \, = \, 0,$ namely $T$ is injective.
$T$ is also self-adjoint as follows:
\begin{equation*}
\begin{split}
(T \alpha, \beta) &= \overline{(\beta, T \alpha)} \\&= \overline{Q(T \beta, T \alpha)}\\&= Q (T \alpha, T \beta)
\\&=(\alpha, T \beta).
\end{split}
\end{equation*}
So we set $U \, = \, \ran( T) \subset D,$ and define
$$F = T^{-1} \, : \, U \rightarrow H$$ which implies
that $F$ is self-adjoint, $Dom (F) \subset D,$ and $$Q(\varphi, \phi) = (F \varphi, \phi)$$
for all $\varphi \in Dom(F)$ and all $\phi \in D$ as claimed.\qed

\medskip
\noindent We obtain as a corollary the main existence and regularity results for operators $\boxbt$ and $\dbarb:$

\smallskip
\newtheorem{existencereg}[proposition]{Corollary}
\begin{existencereg}
\label{existreg} Let $M$ be a compact, orientable, weakly pseudoconvex CR manifold of dimension
at least $5$, endowed with a strongly CR plurisubharmonic function
$\lambda$ as described in Section~\ref{CRpsh}. Let $\alpha$ be a $(0,1)$ form such that
$\alpha \perp \harmon.$ If $\alpha \in H^s$ for some $s\geq 0,$ then there exists
a positive constant $T_s$ and a unique $(0,1)$ form $\varphi_t \perp \harmon$ such that
$$\qbtgs{ \varphi_t}{\phi} =\tnorm{ \alpha}{\phi},$$ $\forall \, \phi \in Dom(\dbarb) \cap Dom(\ad)$ and
$\varphi_t \in H^s$  for all $t \geq T_s.$ Let $N_t$ be the $\dbarb$-Neumann operator that maps $\alpha$ into
$\varphi_t$. We define $N_t$ to be identically $0$ on $\harmon.$ $N_t$ is a bounded operator, and if $\alpha$
is closed, then $u_t \, = \, \adt \, N_t \alpha$
satisfies $$\dbarb u_t = \alpha.$$ Moreover, $u_t \in H^s$ whenever $\alpha \in H^s.$
\end{existencereg}
\smallskip
\noindent {\bf Proof:} We apply Lemma~\ref{Friedrichs} with $\harperp$ as
$H,$ ${\langle| \, \cdot \, | \rangle}_t$ as $|| \, \cdot \, ||,$ $\qbtg{\, \cdot \,}$
as $Q,$ and $Dom(\dbarb) \cap Dom(\ad) \cap \harperp$
as $D.$ Since $\qbtg{\, \cdot \,}^\half$ restricted to $\harperp$
is a norm, and given the way $\dbarb$ and $\ad$ are defined in Section~\ref{defnot},
$Dom(\dbarb) \cap Dom(\ad) \cap \harperp$ is a
Hilbert space under the corresponding inner product $\qbtg{\, \cdot \,}.$
Furthermore, smooth forms are dense in $\harperp$
and in $Dom(\dbarb) \cap Dom(\ad),$ and for $t$ large enough, Lemma~\ref{harmsp}
guarantees that all $\varphi \perp \harmon$ satisfy $$ \sqtnorm{\varphi} \leq C \, \qbtg{\varphi}.$$
By Lemma~\ref{Friedrichs} we obtain the Friedrichs representative $F$
and clearly $$F = \boxbt$$ on $Dom(\boxbt) \cap \harperp$
including all the smooth forms in $\harperp$ and
$Dom (\boxbt) \cap \harperp \subset Dom(F).$ Moreover, it is easy to see that
$Dom(\boxbt^*) \cap \harperp \supset Dom(F^*).$ Since $\boxbt$ and $F$ have been shown to be self-adjoint,
$Dom(\boxbt) \cap \harperp \, = \, Dom(F).$
Finally, it follows from Lemma~\ref{Friedrichs} that there exists a unique $\varphi_t \perp \harmon$
such that $$\qbtgs{ \varphi_t}{\phi} =\tnorm{ \alpha}{\phi},$$
$\forall \, \phi \in Dom(\dbarb) \cap Dom(\ad) \cap \harperp.$ Since
$\alpha \perp \harmon,$ however, $$\qbtgs{ \varphi_t}{\phi} =\tnorm{ \alpha}{\phi},$$
for all $\phi \in Dom(\dbarb) \cap Dom(\ad),$ as needed. This proves the existence of a weak solution
for $\boxbt.$ The regularity follows from Lemma~\ref{estimate}, namely that there exists some $T_s \geq 1$
such that $\varphi_t \in H^s$ for all $t \geq T_s.$

\smallskip
\noindent Now, from the proof of Lemma~\ref{Friedrichs} it is clear that the $\dbarb$-Neumann operator $N_t$ is
bounded in $L^2.$ By definition, $\boxbt \, N_t \alpha \, = \, \alpha$ when $\alpha \perp \harmon.$ $N_t \alpha \in Dom(\boxbt)$ because as shown
above, the domain of the Friedrichs representative $F$ equals $Dom(\boxbt) \cap \harperp,$ so
 $$\dbarb \, \adt \, N_t \alpha + \adt \, \dbarb \, N_t  \alpha = \alpha$$ $$\Downarrow$$
$$\dbarb (\dbarb \, \adt \, N_t \alpha + \adt \, \dbarb \, N_t \alpha) = \dbarb \, \alpha=0$$ $$\Downarrow$$
$$\dbarb \, \adt \, \dbarb \, N_t \alpha = 0.$$ Since $$\tnorm{\dbarb \, \adt \, \dbarb \, N_t \alpha}{\dbarb \, N_t \alpha}
= \sqtnorm{\adt \, \dbarb \, N_t \alpha} = 0.$$ It follows that $$\dbarb \, \adt \, N_t \alpha = \alpha.$$ Let then
$u_t \, = \, \adt \, N_t \alpha.$ $u_t$ satisfies $$\dbarb u_t = \alpha.$$ Furthermore,
 $\alpha \perp \harmon$ implies $u_t \perp \nul(\dbarb).$ To show $\alpha \in H^s$ implies $u_t \in H^s$
first we claim that $N_t$ is a bounded map from $H^s$ to $H^s$ for all $s \geq
0.$ The case $s \, = \, 0$ was shown above,
whereas for $s >0$ we realize that given the definition of $N_t$, Proposition~\ref{estimate} applies
with $\varphi \, = \, N_t \, \alpha,$ so
Equation~\ref{estimateexpr} can be rewritten as
$$||N_t \, \alpha||_s \leq C_t \, (||\alpha||_s + ||N_t \, \alpha||_0) \leq C'_t \, ||\alpha||_s,$$ for some
$t$ dependent constants $C_t$ and $C'_t,$ which proves the claim. Next, we show that $\adt \, N_t$ is
bounded as a map from $H^s$ to $H^s$ when $s \geq 0$ by showing simultaneously that both it and $\dbarb \, N_t$
are bounded. For the
$(0,1)$ form $\alpha \perp \harmon$
\begin{equation*}
\begin{split}
||\dbarb \, N_t \, \alpha||^2_s + ||\adt \, N_t \, \alpha||^2_s &=
\sum_\eta \: \big(\, ||\dze_\eta \Lambda^s \ze_\eta \dbarb \, N_t
\, \alpha||^2_0+ ||\dze_\eta \Lambda^s \ze_\eta \adt \, N_t \,
\alpha||^2_0 \, \big)\\ &= \sum_\eta \: \big(\, (\dze_\eta
\Lambda^s \ze_\eta \adt \, \dbarb \, N_t \,\alpha,\dze_\eta
\Lambda^s \ze_\eta  \, N_t \,\alpha)_0\\&\quad+ (\dze_\eta
\Lambda^s \ze_\eta \dbarb \, \adt \, N_t \, \alpha,\dze_\eta
\Lambda^s \ze_\eta  N_t \, \alpha )_0 \, \big)+errors\\&\leq
||\alpha||_s \, ||N_t \, \alpha||_s + |errors|,
\end{split}
\end{equation*}
where for some small $\eps>0$ and $C_t$ and $C'_t$ dependent of
$t,$
$$|errors| \leq \eps \, \big(||\dbarb \, N_t \, \alpha||^2_s + ||\adt \, N_t \, \alpha||^2_s \big)+ C_t \,
||N_t \, \alpha||_s^2 + C'_t \, \big(||\dbarb \, N_t \,
\alpha||^2_{s-1} + ||\adt \, N_t \, \alpha||^2_{s-1} \big).$$
Given the result for $N_t$ proved above and noticing that for $s
\, = \, 0$
$$||\dbarb \, N_t \, \alpha||_{-1} \leq C \, ||N_t \, \alpha||_0 \leq C_t \, ||\alpha||_0$$ and $$||\adt \, N_t \,
\alpha||_{-1} \leq C_t \, ||N_t \, \alpha||_0 \leq C'_t \,
||\alpha||_0,$$ we obtain inductively that
$$||\dbarb \, N_t \, \alpha||^2_s + ||\adt \, N_t \, \alpha||^2_s \leq C_t \, ||\alpha||^2_s.$$ \qed

\bigskip

\bigskip
\newtheorem{closerangedef}[proposition]{Definition}
\begin{closerangedef}
An operator $P$ has closed range if $\forall \: \alpha \in \overline{\ran(P)},$ where $\overline{\ran(P)}$ is the closure
of the range of $P,$ $\alpha \in \ran(P).$
\end{closerangedef}

\medskip
\noindent We shall prove that the results obtained above impy $\dbarb,$ $\adt,$ and
$\boxbt$ all have closed range in this case.

\smallskip
\newtheorem{closedrangelemma}[proposition]{Lemma}
\begin{closedrangelemma}
\label{closedrange} Let $M$ be a compact, orientable, weakly pseudoconvex CR manifold of dimension
at least $5$, endowed with a strongly CR plurisubharmonic function
$\lambda$ as described in Section~\ref{CRpsh}. For $t$ large enough, the operators $\dbarb$ on functions, $\adt$ on
$(0,2)$ forms, and $\boxbt$ on $(0,1)$ forms have closed range in $L^2$ and all $H^s$ spaces for $s>0.$
\end{closedrangelemma}
\smallskip
\noindent {\bf Proof:} We set out by claiming that the following two statements are equivalent:
\begin{equation*}
\begin{aligned}
&\dbarb \, \alpha = 0 \\ &\alpha \perp \harmon
\end{aligned}
\qquad \Longleftrightarrow \qquad
\begin{aligned}
&\dbarb \, \alpha = 0 \\ &\alpha \perp \nul(\adt)
\end{aligned}
\end{equation*}
"$\Longrightarrow$" $\dbarb \, \alpha \, = \, 0$ means $\alpha \in \nul(\dbarb),$ whereas $\alpha \perp \harmon$
is equivalent to $\alpha \perp \nul(\dbarb) \cap \nul(\adt).$ It follows $\alpha \perp \nul(\adt).$

\noindent "$\Longleftarrow$" $\alpha \perp \nul(\adt)$ implies $\alpha \perp \nul(\dbarb) \cap \nul(\adt),$ i.e.
$\alpha \perp \harmon.$ The claim is proven.

\smallskip \noindent By Corollary~\ref{existreg}, the $\dbarb$ problem can be solved, so
\begin{equation*}
\begin{aligned}
&\dbarb \, \alpha = 0 \\ &\alpha \perp \nul(\adt)
\end{aligned}
\quad \Longleftrightarrow \quad \alpha \in \ran(\dbarb).
\end{equation*}
As in Section~\ref{defnot}, let $L^1_2(M)$ be the set of $(0,1)$ forms in $L^2$ of $M.$ Thus, the set of closed, $(0,1)$
forms in $L^2$ decomposes as follows:
$$L^1_2(M) \cap \nul (\dbarb) = \ran(\dbarb) \oplus \harmon.$$ This implies that
the range of $\dbarb$ is closed. By the same method with $\dbarb$ and $\adt$ exchanged, we see that the following
decomposition holds:
$$L^2_2(M) \cap \nul(\adt) = \ran(\adt) \oplus \harmon$$ for $t$ large enough,
hence the range of $\adt$ is also closed. Finally, by Corollary~\ref{existreg}, $\boxbt$ has a solution for $t$ large
enough, more exactly:
\begin{equation*}
\alpha \perp \harmon \quad \Longleftrightarrow \quad \alpha \in \ran(\boxbt).
\end{equation*}
Therefore, $$L^1_2(M) = \ran(\boxbt) \oplus \harmon,$$ hence $\boxbt$ has closed range in $L^2.$ Notice that the
regularity results contained in Corollary~\ref{existreg} imply that the ranges of $\dbarb,$ $\adt,$ and $\boxbt$
are also closed in $H^s$ for $s>0.$ \qed

\bigskip
\bigskip
\noindent Let $$\hol(M) = \{f \in Dom(\dbarb) \: \big| \: \dbarb f = 0\}.$$ Then the Szeg\"{o} projection $S_{b,t}$ is
the operator that projects $L^2(M)$ to $\hol(M)$ in the ${\langle | \, \cdot \, |\rangle}_t$ norm.

\medskip
\noindent {\bf Facts about the Szeg\"{o} projection:}
\begin{enumerate}
\item[(i)] For all $t \geq T_0,$ $$S_{b,t} = I - \adt \, N_t \, \dbarb,$$
where $I$ is the identity operator;
\item[(ii)] $S_{b,t}$ is exact regular for $t$ large enough, i.e. it is a bounded operator from $H^s$ to $H^s$ for all
$s \geq 0.$
\end{enumerate}

\bigskip
\noindent We shall now provide proofs for both of these facts. The first one was proven by J. J. Kohn for the Bergman
projection, which is an object similar to the Szeg\"{o} projection, except that the former pertains to domains
whereas the latter pertains to the boundaries of domains or to CR manifolds. See for example \cite{Kohnsurvey}.

\smallskip
\newtheorem{Szego1lemma}[proposition]{Lemma}
\begin{Szego1lemma}
\label{Szego1} Let $M$ be a compact, orientable, weakly pseudoconvex CR manifold of dimension at least $5$,
endowed with a strongly CR plurisubharmonic function $\lambda$ as described in Section~\ref{CRpsh}, then
$$S_{b,t} = I - \adt \, N_t \, \dbarb,$$ for all $t \geq T_0.$
\end{Szego1lemma}

\smallskip
\noindent {\bf Proof:} If $f \in \hol(M),$ then $(I - \adt \, N_t \, \dbarb)\, f \, = \, f,$ so the expression for
$S_{b,t}$ holds. Now, if $f \perp \hol(M),$ we want to solve the problem $$\dbarb \, u = \dbarb \, f.$$ According
to Corollary~\ref{existreg}, for all $t \geq T_0,$ $u \, = \, \adt \, N_t \, \dbarb \, f$ is a solution, which implies
that $$\dbarb \, (u-f) = 0,$$ i.e. that $u-f \in \hol(M).$ I claim that $u \perp \hol(M).$ Indeed,
$\forall \: g \in \hol(M)$
$${\langle |\adt \, N_t \, \dbarb \, f , g |\rangle}_t = {\langle | N_t \, \dbarb \, f, \dbarb \, g |\rangle}_t = 0.$$
Since $f \perp \hol(M),$ it turns out that $u-f \perp \hol(M);$ therefore, $u-f \, = \, 0.$
Thus, $(I - \adt \, N_t \, \dbarb)\, f \, = \, 0,$ as expected. The formula for $S_{b,t}$ thus holds on the two
orthogonal subspaces composing $Dom(\dbarb),$ hence on all of it. \qed

\medskip
\noindent For the second fact about the Szeg\"{o} projection, we point the reader to \cite{ChenShaw} for a
proof of the similar result in the case of the Bergman projection, noting that we shall give a different proof here:

\smallskip
\newtheorem{Szego2lemma}[proposition]{Lemma}
\begin{Szego2lemma}
\label{Szego2} Let $M$ be a compact, orientable, weakly pseudoconvex CR manifold of dimension at least $5$,
endowed with a strongly CR plurisubharmonic function $\lambda$ as described in Section~\ref{CRpsh}, then for
$t$ large enough, $S_{b,t}$ is bounded as a map from $H^s$ to $H^s$ for all $s \geq 0.$
\end{Szego2lemma}

\smallskip
\noindent {\bf Proof:} If $f \in \hol(M),$ then $S_{b,t}$ is the identity, and there is nothing to prove. Otherwise, let
$\dbarb \, f \, = \, \alpha$ and by Corollary~\ref{existreg}, for each $t \geq T_0$ there exist a function
$u_t \perp \hol(M)$ and a $(0,1)$ form $\varphi_t \perp \harmon$ such that
$\dbarb \, f \, = \, \dbarb \, u_t \, = \, \alpha,$ $N_t \, \alpha \, = \, \varphi_t,$ and $\adt \, \varphi_t \, = \, u_t.$
We need to show that there exists a $t$ dependent constant $C_t$ for each $s \geq 0$ such that
$$||u_t ||_s \leq C_t ||f||_s.$$ The first step is to show that there exist constants $C$ and $C_t$ such that
$$\sqtnorm{\Lambda^s  u_t} \leq C \, \sqtnorm{\Lambda^s  f}+C \, \sqtnorm{\Lambda^s  \varphi_t}
+C_t \, \sqtnorm{\Lambda^{s-1} \varphi_t}.$$ We proceed as follows:
\begin{equation*}
\begin{split}
\sqtnorm{\Lambda^s  u_t} &= \tnorm{\Lambda^s \adt \, \varphi_t}{\Lambda^s u_t}
\\&=\tnorm{[\Lambda^s, \adt]\, \varphi_t}{\Lambda^s u_t}+\tnorm{\Lambda^s  \varphi_t}{[\dbarb,\Lambda^s] \, u_t}
+\tnorm{\Lambda^s  \varphi_t}{\Lambda^s \dbarb f}\\&=\tnorm{[\Lambda^s, \adt]\, \varphi_t}{\Lambda^s u_t}+\tnorm{\Lambda^s  \varphi_t}{[\dbarb,\Lambda^s] \, u_t}
+\tnorm{\Lambda^s  \varphi_t}{[\Lambda^s, \dbarb] f}\\&\quad+\tnorm{[\adt, \Lambda^s]  \varphi_t}{\Lambda^s f}
+\tnorm{\Lambda^s u_t}{\Lambda^s f}\\&\leq C \, \sqtnorm{\Lambda^s  f}+C \, \sqtnorm{\Lambda^s  \varphi_t}
+C_t \, \sqtnorm{\Lambda^{s-1} \varphi_t}+ \eps \,\sqtnorm{\Lambda^s  u_t},
\end{split}
\end{equation*}
for some small positive $\eps >0,$ which concludes the proof of the first step. As a second step, we would like to show
that for $t$ large enough and some small positive $\eps >0$ independent of $t,$
$$\sqtnorm{\Lambda^s  \varphi_t} \leq C_t \, \sqtnorm{\Lambda^s  f}+C_t \, \sqtnorm{\Lambda^{s-1} \varphi_t}
+ \eps \,\sqtnorm{\Lambda^s  u_t}.$$ We start by using Proposition~\ref{global} to conclude that for $t$ large enough,
\begin{equation*}
t \sqtnorm{\Lambda^s \varphi_t} \leq C \qbtg {\Lambda^s
\varphi_t}+C_t (\: \sum_\nu \:{||\dze_\nu \potdg {\nu} \ze_\nu
\Lambda^s \varphi_t^\nu||}^2_0 +||\varphi_t||^2_{s-1}).
\end{equation*}
Using that $\qbtgs{ \varphi_t}{\phi} =\tnorm{ \alpha}{\phi}$ $\forall \, \phi \in Dom(\dbarb) \cap Dom(\ad),$
part (ii) of Lemma~\ref{qbtcommlemma}, and a manipulation similar to the one in the first step, we obtain that
$$\qbtg {\Lambda^s \varphi_t} \leq C \, \sqtnorm{\Lambda^s  f}+C \, \sqtnorm{\Lambda^s  \varphi_t}
+C_t \, \sqtnorm{\Lambda^{s-1} \varphi_t}+ \eps \,\sqtnorm{\Lambda^s  u_t}$$ for some small $\eps >0$ independent of
$t,$ whereas by Lemma~\ref{Sob1global} along with part (i) of
Lemma~\ref{qbtcommlemma},
\begin{equation*}
\begin{split}
\sum_\nu \:{||\dze_\nu \potdg {\nu} \ze_\nu \Lambda^s \varphi_t^\nu||}^2_0 &\leq C \, (\: \sum_\nu \: \qbtgs{\dze_\nu
\Lambda^{-1} \potdg {\nu} \ze_\nu \Lambda^s \varphi_t^\nu}{G_t \dze_\nu \Lambda^{-1} \potdg {\nu} \ze_\nu \Lambda^s
\varphi_t^\nu}+||\varphi_t||^2_{s-1})\\&\leq \sum_\nu \: \tnorm{\dze_\nu
\Lambda^{-1} \potdg {\nu} \ze_\nu \Lambda^s \dbarb f}{G_t \dze_\nu \Lambda^{-1} \potdg {\nu} \ze_\nu \Lambda^s
\varphi_t^\nu}+C_t \, \sqtnorm{\Lambda^{s-1} \varphi_t}\\&\quad+ \eps \: \sum_\nu \:
\qbtg{G_t \dze_\nu \Lambda^{-1} \potdg {\nu} \ze_\nu \Lambda^s \varphi_t^\nu}.
\end{split}
\end{equation*}
$G_t$ is self-adjoint, so it can be moved to the other side of the first term of the right-hand side, which implies
\begin{equation*}
\begin{split}
\sum_\nu \:{||\dze_\nu \potdg {\nu} \ze_\nu \Lambda^s \varphi_t^\nu||}^2_0 &\leq
C_t \, \sqtnorm{\Lambda^s  f}+C_t \, \sqtnorm{\Lambda^{s-1} \varphi_t}\\&\quad+ \eps \: \sum_\nu \:
\qbtg{G_t \dze_\nu \Lambda^{-1} \potdg {\nu} \ze_\nu \Lambda^s \varphi_t^\nu}.
\end{split}
\end{equation*}
Since $G_t$ is not inverse $t$ dependent, the last term of the right-hand side is handled using just
Equation~\ref{qbtcommexpr} and not any of the simpler parts of Lemma~\ref{qbtcommlemma}. From the five terms thus obtained,
three are easily estimated as above and only two require a further explanation, namely
$$\tnorm{[\dbarb,G_t \dze_\nu \Lambda^{-1} \potdg {\nu} \ze_\nu \Lambda^s] \varphi_t^\nu}
{\dbarb G_t \dze_\nu \Lambda^{-1} \potdg {\nu} \ze_\nu \Lambda^s \varphi_t^\nu}$$ and
$$\tnorm{[\adt,G_t \dze_\nu \Lambda^{-1} \potdg {\nu} \ze_\nu \Lambda^s] \varphi_t^\nu}
{\adt G_t \dze_\nu \Lambda^{-1} \potdg {\nu} \ze_\nu \Lambda^s \varphi_t^\nu}.$$ Clearly, what needs to be done is
to commute $G_t$ outside of $\dbarb$ and $\adt$ respectively and then to throw it onto the other side where it is harmless.
Thus we obtain that $$\sum_\nu \:{||\dze_\nu \potdg {\nu} \ze_\nu \Lambda^s \varphi_t^\nu||}^2_0 \leq
C_t \, \sqtnorm{\Lambda^s  f}+C_t \, \sqtnorm{\Lambda^{s-1} \varphi_t}+ \eps' \, \sqtnorm{\Lambda^s \varphi_t}$$
for some small $\eps' >0$ independent of $t,$ so altogether,
\begin{equation*}
t \sqtnorm{\Lambda^s \varphi_t} \leq C_t \, \sqtnorm{\Lambda^s  f}+C \, \sqtnorm{\Lambda^s  \varphi_t}
+C_t \, \sqtnorm{\Lambda^{s-1} \varphi_t}+ \eps \,\sqtnorm{\Lambda^s  u_t},
\end{equation*}
where $\eps$ is the small, positive, $t$ independent constant from the estimate for $\qbtg {\Lambda^s \varphi_t}.$
Let $t$ be large enough that $C \, \sqtnorm{\Lambda^s  \varphi_t}$ can be absorbed on the left-hand side;
then, step two follows. Steps one and two show that $$\sqtnorm{\Lambda^s u_t} \leq C_t (\sqtnorm{\Lambda^s  f}+
\sqtnorm{\Lambda^{s-1} \varphi_t}.$$ Now, by Lemma~\ref{closedrange}, both $\dbarb$ and $\adt$ have closed range in $L^2$ and
any $H^s$ space for $s>0.$ It is an elementary application of the Closed Graph Theorem (see for example Section 4.1 of
\cite{ChenShaw}) that closeness of range in $H^s$ of
$\adt$ is equivalent to the following estimate for $s \geq 1$ $${\langle| \Lambda^{s-1} \varphi_t |\rangle}_t \leq C_t \,
{\langle| \Lambda^{s-1} u_t |\rangle}_t.$$ Therefore, $$\sqtnorm{\Lambda^s u_t} \leq C_t \, (\sqtnorm{\Lambda^s f}
+\sqtnorm{\Lambda^{s-1} u_t}).$$ We now pass to the usual Sobolev norms since ${\langle | \cdot | \rangle}_t$ is equivalent
to the $L^2$ norm and restate the estimate as $$||u_t||^2_s \leq C_t \, (||f||^2_s + ||u_t||^2_{s-1}).$$  Since the
orthogonal projection is a bounded operator in $L^2,$ $$||u_t||_0 \leq C_t \,|| f ||_0,$$ thus we proceed inductively
to finish the proof of the lemma. \qed

\medskip
\noindent Using these facts about the Szeg\"{o} projection and a Mittag-Leffler-type argument, one can find a smooth
solution to the $\dbarb$ problem when starting with smooth data. This proof originally appeared in \cite{Kohnmethods}.

\smallskip
\newtheorem{smoothsollemma}[proposition]{Lemma}
\begin{smoothsollemma}
\label{smoothsolution} Let $M$ be a compact, orientable, weakly pseudoconvex CR manifold of dimension
at least $5$, endowed with a strongly CR plurisubharmonic function
$\lambda$ as described in Section~\ref{CRpsh}. If $\alpha$ is a closed $(0,1)$ form such that $\alpha \in \smooth(M),$
then there exists a function $u \in \smooth(M)$ satisfying $\dbarb \, u \, = \, \alpha.$
\end{smoothsollemma}
\smallskip
\noindent {\bf Proof:} By Corollary~\ref{existreg}, for each $k \, = \, 1,2, \dots,$ there exists some $u_k \in H^k$
such that $\dbarb \, u_k \, = \, \alpha.$ We will modify each $u_k$ by an element of $\nul(\dbarb)$ in order to construct
a telescoping series that is in $H^k$ for each $k\geq 1.$ To do so, we need to show first that $H^s \cap \nul(\dbarb)$
is dense in $H^k \cap \nul(\dbarb)$ for each $s >k.$ Let $g$ be any element of $H^k \cap \nul(\dbarb).$
Smooth functions are dense in all $H^k (M),$ so there exists a sequence $\{g_i\}_i$ such that $g_i \in \smooth(M)$
and $g_i \rightarrow g$ in $H^k.$ $\dbarb \, g \, = \, 0$ implies that $$g - S_{b,t} \,g = \adt \, N_t \, \dbarb \, g = 0,$$
so $g \, = S_{b,t} \, g.$ Let $g'_i \, = \, S_{b,t}\, g_i.$ $g'_i \in H^s \cap \nul(\dbarb)$ since the Szeg\"{o} projection
is bounded as a map from $H^s$ to $H^s$ and for the same reason, $g'_i \rightarrow g$ in $H^k.$ Thus indeed,
$H^s \cap \nul(\dbarb)$ is dense in $H^k \cap \nul(\dbarb).$ Using this fact, we inductively construct a sequence
$\{ \tilde{u}_k \}_k$ as follows: $$\tilde{u}_1 = u_1,$$ $$\tilde{u}_2 = u_2 + v_2,$$ where $v_2 \in H^2 \cap \nul(\dbarb)$
is such that $$||\tilde{u}_2 - \tilde{u}_1||_1 \leq 2^{-1},$$ and in general, $$\tilde{u}_{k+1} = u_{k+1} + v_{k+1},$$
where $v_{k+1} \in H^{k+1} \cap \nul(\dbarb)$ is such that $$||\tilde{u}_{k+1} - \tilde{u}_k||_k \leq 2^{-k}.$$
Clearly, $\dbarb \, \tilde{u}_k \, = \, \alpha,$ so we set
$$u = \tilde{u}_J + \sum_{k=J}^\infty (\tilde{u}_{k+1} - \tilde{u}_k), \: \: J \in \N.$$ It follows that
$u \in H^k$ for each $k \in \N,$ hence that $u \in \smooth(M)$ and $\dbarb \, u \, = \, \alpha.$  \qed

\medskip
\noindent Following the method in \cite{Kohnglobalreg}, we shall prove next that as $t$ increases,
the harmonic spaces $\harmon$ get smoother
and smoother. For that, however, we need the following elementary result:

\smallskip
\newtheorem{Hscompactness}[proposition]{Lemma}
\begin{Hscompactness}
\label{Hscompactnesslemma} Let $\{ \varphi_k \}$ be a sequence of
$(0,1)$ forms such that $||\varphi_k||_s \leq C,$ for some $C$
independent of $k,$ then there exists a subsequence $\{
\varphi_{k_j} \}$ which converges to some $\varphi \in H^s$ in each $H^{s'}$ with $s'<s.$
\end{Hscompactness}
\smallskip
\noindent {\bf Proof:} Since $M$ is compact, by Rellich's theorem,
there exists a subsequence of $\{ \varphi_k \}$ which converges to
some $\varphi,$ where $\varphi \in H^{s'}$ for all $s' < s.$
$||\varphi_k||_s \leq C$ implies the subsequence which converges
to $\varphi$ in each $H^{s'}$ also converges weakly to some
$\varphi'$ in $H^s$ because $H^s$ is a separable Hilbert space. We
shall show that $\varphi \, = \, \varphi'$ almost everywhere.
According to \cite{RieszNagy} $\S 38,$ a theorem of Banach-Saks
guarantees that in every Hilbert space, a sequence that converges
weakly has a subsequence whose arithmetic means converge strongly
to the weak limit. Let $\{ \varphi_{k_j} \}$ be the subsequence of
$\{ \varphi_k \}$ whose arithmetic means converge strongly to
$\varphi'$ in $H^s$ hence also in $H^{s-1}$ Since $\{
\varphi_{k_j} \}$ converges to $\varphi$ in $H^{s-1},$ it follows
its arithmetic means also converge to $\varphi$ in $H^{s-1}.$ It
is then easy to see that indeed $\varphi \, = \, \varphi' \: a.e.$
which means $\varphi \in H^s.$ \qed

\smallskip
\noindent {\bf Remark:} It is not true in general that the subsequence
$\{\varphi_{k_j} \}$ converges to $\varphi$ in $H^s$ as well, else it would follow that $H^s$ were finite dimensional,
which is clearly not true.

\bigskip
\newtheorem{harmsoblemma}[proposition]{Lemma}
\begin{harmsoblemma}
\label{harmsob} Let $M$ be a compact, orientable, weakly pseudoconvex CR manifold of dimension
at least $5$, endowed with a strongly CR plurisubharmonic function
$\lambda$ as described in Section~\ref{CRpsh}. For $t$ large enough, $\:\harmon \subset H^s.$
\end{harmsoblemma}
\smallskip
\noindent {\bf Proof:} Let $\dim \harmon \, = \, S$ and $\theta_1, \dots, \theta_S$ be
a basis for $\harmon.$ If $S \, = \, 0,$ there is nothing to prove; otherwise, assume $\theta_0 \, = \, 0$
and $\theta_j \in H^s$ for all $0 \leq j \leq l < S.$ We shall construct some $\theta \in H^s \cap \harmon$ such that
${\langle| \theta |\rangle}_t \, = \, 1$ and ${\langle| \theta, \theta_j |\rangle}_t \, = \, 0$ for all
$j \leq l.$ Define a form $\qbtk{\, \cdot \,}{\, \cdot \,}$ as follows:
$$\qbtk{\varphi}{\phi} = \qbtgs{\varphi}{\phi} + \frac{1}{k} \, {\langle| \varphi, \phi |\rangle}_t$$ There exists
then some constant $C_k$ depending on $k$ such that
$${\langle| \varphi |\rangle}^2_t \leq C_k \, \qbtk{\varphi}{\varphi}.$$ Let $\alpha \in H^s$ be a
$(0,1)$ form such that $\alpha \perp \theta_j$ for all $j \leq l$ and $\alpha$ is not orthogonal to $\theta_{l+1},$
then by Lemma~\ref{Friedrichs} there exists a unique $\varphi_k \in Dom(\dbarb) \cap Dom(\ad)$ such
that $$\qbtk{\varphi_k}{\phi} ={\langle| \alpha, \phi |\rangle}_t,$$ for all $\phi \in Dom(\dbarb) \cap Dom(\ad).$
Moreover, the additional term $\frac{1}{k} \, {\langle| \varphi, \phi |\rangle}_t$ present in
$\qbtk{\, \cdot \,}{\, \cdot \,}$ does not prevent a proof similar to that of Proposition~\ref{estimate} to be
carried through in order to show that for $t$ large enough, $\varphi_k \in H^s,$ and
\begin{equation}
||\varphi_k||_s \leq C_t \, (||\alpha||_s + ||\varphi_k||_{\, 0}), \label{seqestimate}
\end{equation}
where $C_t$ is a $t$ dependent constant. We claim that the sequence
$\{||\varphi_k||_{\, 0}\}_k$ is unbounded. Assume that it is bounded, then by Equation~\ref{seqestimate} and
Lemma~\ref{Hscompactnesslemma} there exists a subsequence $\{\varphi_{k_i}\}$ which converges in each $H^{s'}$ for
$s'<s$ to some $\varphi \in H^s.$ Moreover, $\varphi$ satisfies
$$\qbtgs{\varphi}{\phi} ={\langle| \alpha, \phi |\rangle}_t,$$ for all $\phi \in Dom(\dbarb) \cap Dom(\ad).$ Let
$\phi \, = \, \theta_j$ in order to obtain a contradiction. It turns out the left-hand side is zero for all $j$, whereas
the right-hand side is non-zero for $j \, = \, l+1.$ The claim is established. It follows that there exists some
subsequence $\{\varphi_{k_i}\}$ such that $\lim_{i \rightarrow \infty}||\varphi_{k_i}||_{\, 0} \, = \, \infty.$ We set
$\gamma_i \, = \, \frac{\varphi_{k_i}}{||\varphi_{k_i}||_{\, 0}}$ and note that then
\begin{equation}
Q^{\, k_i}_{b,t} ({\gamma_i},{\phi}) = \frac{{\langle| \alpha, \phi |\rangle}_t}{||\varphi_{k_i}||_{\, 0}}.
\label{qbtkgamma}
\end{equation}
Once again, by Equation~\ref{seqestimate} and Lemma~\ref{Hscompactnesslemma} there exists a subsequence
of $\{\gamma_i\}_i$ which converges in each $H^{s'}$ for $s'<s$ to some $\theta \in H^s.$ $\theta$ is such that
${\langle| \theta |\rangle}_t \, = \, 1$ and $\qbtg{\theta} \, = \, 0,$ i.e. $\theta \in \harmon.$
Finally, Equation~\ref{qbtkgamma} with $\phi \, = \, \theta_j$ for $0 \leq j \leq l$ implies
${\langle| \theta, \theta_j |\rangle}_t \, = \, 0$ as needed. \qed

\medskip
\noindent Now, the last lemma of this section:

\smallskip
\newtheorem{isomorphylemma}[proposition]{Lemma}
\begin{isomorphylemma}
\label{isomorphy} Let $M$ be a compact, orientable, weakly pseudoconvex CR manifold of dimension
at least $5$, endowed with a strongly CR plurisubharmonic function
$\lambda$ as described in Section~\ref{CRpsh}. Then $$H_0^{0,1}(M,\dbarb) \cong H_s^{0,1}(M,\dbarb) \cong
H_\infty^{0,1}(M,\dbarb) \quad \forall \: s>0.$$
\end{isomorphylemma}
\smallskip
\noindent {\bf Proof:} We start by proving that $H_0^{0,1}(M,\dbarb) \cong H_s^{0,1}(M,\dbarb).$ Let $t$ be
large enough that by the previous lemma, $\harmon \subset H^s$ and such that $t \geq T_s,$ where $T_s$ is the
threshold value of $t$ given by Corollary~\ref{existreg} in order for the $\dbarb$ problem with $H^s$ data to
be solved in $H^s.$ We define $$\pi_{t,s} : H^s \mapsto \harmon.$$
The previous lemma then implies that $\pi_{t,s}$ is surjective. Now let $\alpha, \beta \in H^s$ be such that
$\pi_{t,s} (\alpha) \, = \, \pi_{t,s} (\beta).$ Since $\pi_{t,s}(\alpha), \: \pi_{t,s}(\beta) \in H^s,$ it follows
$\alpha - \beta \in H^s.$ By Corollary~\ref{existreg}, $\alpha \sim \beta,$ hence $\pi_{t,s}$ is also injective as
a map from $H_s^{0,1}(M,\dbarb)$ to $\harmon.$
By the same procedure, the map $$\pi_{t,0} : H_0^{0,1}(M,\dbarb) \mapsto \harmon$$ is bijective, so
$H_0^{0,1}(M,\dbarb) \cong H_s^{0,1}(M,\dbarb).$ Next, we shall prove that $H_0^{0,1}(M,\dbarb) \cong
H_\infty^{0,1}(M,\dbarb).$ Let us take another look at the map $$\pi_{t,0} : L^2 \mapsto \harmon.$$
$\pi_{t,0}$ is a linear operator which is bounded in $L^2$ and has norm $1.$ Also, $\smooth$ is dense in $L^2,$
namely $\forall \: \varphi \in L^2$ and $\forall \: \eps>0,$ there exists some $\phi \in \smooth(M),$ such that
$$||\varphi - \phi||_{\, 0} < \eps,$$ which implies
$$||\pi_{t,0}\varphi - \pi_{t,0}\phi||_{\, 0} = ||\pi_{t,0}(\varphi - \phi)||_{\, 0} \leq
||\varphi - \phi||_{\, 0} < \eps,$$ so $\pi_{t,0}(\smooth(M))$ is dense in $\harmon.$ By Lemma~\ref{harmsp},
$\harmon$ is a finite dimensional vector space, therefore the span of $\pi_{t,0}(\smooth(M))$ has to be the whole
of $\harmon.$ In other words, there exist some $\phi_1, \dots, \phi_S \in \smooth(M)$ such that $\pi_{t,0}(\phi_1),
\dots, \pi_{t,0}(\phi_S)$ is a basis for $\harmon.$ It follows $$\pi_{t,0} \Big|_{\smooth(M)} : L^2 \mapsto \harmon$$
is surjective. In order to prove that $\pi_{t,0}$ is injective as a map from $H_\infty^{0,1}(M,\dbarb)$ to $\harmon$
one has to employ a Mittag-Leffler type construction similar to the one used in the proof of Lemma~\ref{smoothsolution}. Let $\alpha \in \smooth(M).$ We shall construct some $\theta \in \smooth(M)$ such that $\alpha \sim \theta$ in each
$H^s,$ hence also in $\smooth.$ To start, we set $\theta_0 \, = \, \pi_{t,0}(\alpha)$ which means that
$\alpha \sim \theta_0$ in $L^2.$ Next, for each $s \in \N,$ $s \geq 1$ we can choose some $t_s$ such that
$\harmof{t_s} \subset H^s$ and the $\dbarb$ problem with $H^s$ data can be solved in $H^s$ by the previous
lemma and by Lemma~\ref{existreg}. We shall use the maps $$\pi_{t_s, s}: \smooth \mapsto \harmof{t_s}.$$
Clearly, $\pi_{t_1,1}(\alpha) - \pi_{t,0}(\alpha) \, = \, \dbarb \beta_0,$ for some $\beta_0 \in L^2.$
Since $\dbarb$ is defined as a closed operator, there exists a sequence $\{ \beta_{0,j} \}_j$ such that
$\beta_{0,j} \in \smooth,$ $\beta_{0,j} \rightarrow \beta_0$ in $L^2,$ and $\dbarb \beta_{0,j}
\rightarrow \dbarb \beta_0$ also in $L^2.$ Define $\theta_1 \, = \, \pi_{t_1, 1}(\alpha) -
\dbarb \beta_{0,j},$ where $j$ is chosen so that $||\theta_0 - \theta_1||_0 < 1.$ Similarly,
$\pi_{t_2,2}(\alpha)-\pi_{t_1,1}(\alpha) \, = \, \dbarb \beta_1,$ for some $\beta_1 \in H^1.$
Since the range of $\dbarb$ is closed in each $H^s$ space, it follows we can find a sequence
$\{\beta_{1,j}\}_j$ such that $\beta_{1,j} \in \smooth$ and $\dbarb \beta_{1,j} \rightarrow
\dbarb \beta_1$ in $H^1.$ Define $\theta_2 \, = \, \pi_{t_2,2}(\alpha) - \dbarb \beta_{1,j},$
where $j$ is chosen so that $||\theta_1 - \theta_2||_1 < \half.$ Inductively,
$\pi_{t_{s+1},{s+1}}(\alpha)-\pi_{t_s,s}(\alpha) \, = \, \dbarb \beta_s,$ for some
$\beta_s \in H^s,$ so there exists a sequence $\{\beta_{s,j}\}_j$ such that
$\beta_{s,j} \in \smooth$ and $\dbarb \beta_{s,j} \rightarrow \dbarb \beta_s$ in $H^s.$ Define
$\theta_{s+1} \, = \, \pi_{t_{s+1}, {s+1}}(\alpha) - \dbarb \beta_{s,j},$ where $j$ is chosen
so that $||\theta_s - \theta_{s+1}||_s < \frac{1}{2^s}.$ Set
\begin{equation*}
\begin{split}
\theta &= \theta_0 + \sum_{k=0}^\infty \:(\theta_{k+1} - \theta_k)\\&= \theta_m +
\sum_{k=m}^\infty \: (\theta_{k+1} - \theta_k).
\end{split}
\end{equation*}
By construction, $\theta \in H^s$ and $\theta \sim \alpha$ in $H^s$ for all $s\geq 0,$ so
$\theta \in \smooth$ and $\theta \sim \alpha$ in $\smooth.$ Now, to prove that $\pi_{t,0}$
is injective as a map from $H_\infty^{0,1}(M,\dbarb)$ to $\harmon,$ let $\alpha_1, \alpha_2 \in \smooth$
be such that $\pi_{t,0}(\alpha_1) \sim \pi_{t,0}(\alpha_2)$ in $L^2.$ By the construction above,
there exist $\theta', \theta'' \in \smooth$ such that $\theta' \sim \alpha_1$ in $\smooth$ and
$\theta'' \sim \alpha_2$ in $\smooth.$ It follows $\theta' \sim \alpha_1 \sim \pi_{t,0}(\alpha_1)$
in $L^2$ and $\theta'' \sim \alpha_2 \sim \pi_{t,0}(\alpha_2)$ in $L^2.$ This means $\theta' \sim \theta''$
in $L^2,$ but given the way $\theta'$ and $\theta''$ are defined, it follows $\theta' \sim \theta''$
in each $H^s$ space, hence also in $\smooth.$ Thus, $\alpha_1 \sim \alpha_2$ in $\smooth$ which
amounts to $\pi_{t,0}$ being injective, thus bijective, as a map from $H_\infty^{0,1}(M,\dbarb)$
to $\harmon.$ Since $\pi_{t,0}$ is also a bijection from $H_0^{0,1}(M,\dbarb)$ to $\harmon,$ the
conclusion of the lemma follows. \qed

\medskip
\medskip
\noindent {\bf Proof of the Main Theorem for $(0,1)$ forms:} Lemma~\ref{closedrange} proves part (i), then
Corollary~\ref{existreg} establishes (ii) and (iii). Lemma~\ref{smoothsolution} proves (iv), whereas (v) is a
consequence of Lemmas~\ref{harmsp} and ~\ref{isomorphy}. \qed

\bigskip
\section{Higher Level Forms}
\label{highforms}

\medskip
\medskip
We set out next to outline the proof of the Main Theorem for
$(p,q)$ forms with $1 \leq q \leq n-2.$ As noted in
Section~\ref{defnot}, the holomorphic part of the forms is
irrelevant for the argument, so in order to simplify the notation,
we shall consider only $(0,q)$ forms. The constructions in
Section~\ref{CRpsh} carry over as long as we impose that upper
bound $\eps_G$ for the absolute value of the errors $e_{ij}$ in
the bracket $[L_i, \lbar{j}]$ is such that $$q \, \eps_G \ll 1.$$
The construction of the global norm carries over as well, namely
if $\varphi \, = \, \sum_{|I|=q} \, \varphi_I \omi{I}$ is a
$(0,q)$ form, then $${\langle| \varphi |\rangle}_t^2 = \sum_\nu
\sum_{|I|=q} \: (|||\dze_\nu \pptg {\nu} \ze_\nu
\varphi_I^\nu|||_t^2 + ||\dze_\nu \potg {\nu} \ze_\nu
\varphi_I^\nu ||_0^2 + |||\dze_\nu \pmtg {\nu} \ze_\nu
\varphi_I^\nu|||_{-t}^2),$$ where $\varphi_I^\nu$ is the $I^{th}$
component of $\varphi$ expressed in the local coordinates on
$U_\nu.$ It follows that the computation of the adjoint of
$\dbarb$ with respect to this norm done for $(0,1)$ forms is
actually the same for $(0,q)$ forms, so we obtain the same
estimate for the energy form $\qbtg{\, \cdot \,}$ as in
Section~\ref{microlocsection}:
\begin{equation*}
\begin{split}
&K \qbtg {\varphi}+K_t \sum_\nu \sum_{|I|=q} \:{||\dze_\nu \potdg {\nu} \ze_\nu
\varphi_I^\nu||}^2_0+O({\langle |\varphi|\rangle}^2_t) +O_t
(||\varphi||^2_{-1})\\&\geq \sum_\nu \sum_{|I|=q} \: \qbtp{\dze_\nu \pptg {\nu}
\ze_\nu \varphi_I^\nu}+ \sum_\nu \sum_{|I|=q} \:\qbto{\dze_\nu \potg {\nu}
\ze_\nu  \varphi_I^\nu}\\&\quad + \sum_\nu \sum_{|I|=q} \: \qbtm{\dze_\nu \pmtg
{\nu} \ze_\nu \varphi_I^\nu}
\end{split}
\end{equation*}
The computations for $\qbtp{\, \cdot \,}$ and $\qbtm{\, \cdot \,}$ are different, however, in the case of $(0,q)$ forms
than their counterparts for $(0,1)$ forms. Let us outline first the one for $\qbtp{\, \cdot \,}.$ Just as in the case of
$(0,1)$ forms, we restrict the discussion to a small neighborhood $U',$ which acts as the model for each $U_\nu$ in the
covering of $M.$

\smallskip
\newtheorem{qbtpcomp0q}{Lemma}[section]
\begin{qbtpcomp0q}
Let $\varphi$ be a $(0,q)$ form \label{qbtp0q} supported in $U',$
$\varphi \in Dom(\dbarb) \cap Dom(\ad).$ If $[L_i, \lbar{j}]=
c_{ij}T + \sum_k a^k_{ij}\lbar{k} + \sum_k b^k_{ij} \lbar{k}^{*,
\, t} + g_{ij} +t \, e_{ij},$ where $a^k_{ij}$, $b^k_{ij}$,
$g_{ij},$ and $e_{ij}$ are \smooth functions independent of $t$
and $c_{ij}$ are the coefficients of the Levi form, then there
exists $1 \gg \eps' > 0$ such that
\begin{equation*}
\begin{split}
\qbtp{\varphi} &\geq (1-\eps')\, \sum_{|I|=q} \sum_i |||\lbar{i}
\varphi_I|||^2_{\, t} +  \sum_{|I|=q} \sum_{i \in I} \Re \{(c_{ii}
T \varphi_I, \varphi_I)_t \}\\&\quad  - \sum_{|I|=|J|=q} \sum_{i,j
\atop {i \neq j}} \eps^{iI}_{jJ} \, \Re \{(c_{ij} T \varphi_J,
\varphi_I)_t \}+ t \sum_{|I|=q} \sum_{i \in I} \Re \{ ( \lbar{i}
L_i (\lambda) \varphi_I, \varphi_I)_t \}\\&\quad- t
\sum_{|I|=|J|=q} \sum_{i,j \atop {i \neq j}} \eps^{iI}_{jJ}\, \Re
\{ (\lbar{j} L_i (\lambda) \varphi_J, \varphi_I)_t \}+ t
\sum_{|I|=q} \sum_{i \in I} \Re \{ ( e_{ii} \varphi_I,
\varphi_I)_t \}\\&\quad -t\sum_{|I|=|J|=q} \sum_{i,j \atop {i \neq
j}} \eps^{iI}_{jJ} \, \Re \{(e_{ij} \varphi_J, \varphi_I)_t \} +
O(|||\varphi|||^2_{\, t} ),
\end{split}
\end{equation*}
where $ \eps^{jJ}_{iI}$ is the sign of the permutation changing ${jJ}$ to ${iI}.$
\end{qbtpcomp0q}

\smallskip
\noindent {\bf Proof:} Since $\varphi \, = \, \sum_{|I|=q} \, \varphi_I \omi{I},$
$$\dbarb \varphi = \sum_{|I|=q} \sum_{i \notin I} \eps^{iI}_{\ord{iI}} \lbar{i} (\varphi_I) \, \omi{\ord{iI}} +
lower \: order \: terms,$$ $\ord{iI}$ denoting the multi-index $iI$ given canonically, i.e. in increasing order.
Therefore,
\begin{equation*}
\begin{split}
|||\dbarb \varphi|||^2_t &= \sum_{|I|=|J|=q} \sum_{i \notin I}
\sum_{j \notin J} \, \eps^{iI}_{jJ} \, {(\lbar{j}(\varphi_J),
\lbar{i} (\varphi_I))}_t + O \Big(|||\varphi|||_t^2 +(\sum_{|I|=q}
\sum_{i \notin I} |||\lbar{i}\varphi_I|||_t^2)^\half \,
|||\varphi|||_t \Big)\\&= \sum_{|I|=q} \sum_{i \notin I}
|||\lbar{i}\varphi_I|||_t^2+\sum_{|I|=|J|=q} \sum_{i \notin I}
\sum_{{j \notin J} \atop i \neq j} \, \eps^{iI}_{jJ} \,
{(\lbar{i}^{*,t}(\varphi_J), \lbar{j}^{*,t}
(\varphi_I))}_t\\&\quad+\sum_{|I|=|J|=q} \sum_{i \notin I}
\sum_{{j \notin J} \atop i \neq j} \, \eps^{iI}_{jJ} \,
{([\lbar{i}^{*,t},\lbar{j}]\varphi_J,  \varphi_I)}_t + O
\Big(|||\varphi|||_t^2 +(\sum_{|I|=q} \sum_{i \notin I}
|||\lbar{i}\varphi_I|||_t^2)^\half \, |||\varphi|||_t \Big),
\end{split}
\end{equation*}
by the same procedure of integrating by parts and commuting as in the proof of Lemma~\ref{qbtp}. The adjoint $\adp$ is
given by $$\adp \varphi = \sum_{|K'| = q} \sum_{{k \in K'} \atop {\ord{kK} = K'}} \eps^{kK}_{K'} \, \lbar{k}^{*,t}
(\varphi_{K'}) \, \omi{K} + lower \: order \: terms.$$ If we apply integration by parts to the cross terms, thus
converting $\lbar{i}^{*,t}$ into $\lbar{i},$ it follows that
\begin{equation*}
\begin{split}
|||\adp \varphi|||^2_t &= \sum_{|I|=|J|=q} \sum_{{j \in I} \atop
{\ord{j K}=I}} \sum_{{i \in J} \atop {\ord{i K} = J}} \,
\eps^{jK}_I \eps^{iK}_J \, {(\lbar{i}^{*,t} (\varphi_J),
\lbar{j}^{*,t} (\varphi_I))}_t\\&\quad + O \Big(|||\varphi|||_t^2
+(\sum_{|I|=q} \sum_i |||\lbar{i}\varphi_I|||_t^2)^\half \,
|||\varphi|||_t \Big)\\&= \sum_{|I|=q} \sum_{i \in I}
|||\lbar{i}\varphi_I|||_t^2- \sum_{|I|=|J|=q} \sum_{j \in I}
\sum_{{i \in J} \atop {i \neq j}} \, \eps^{iI}_{jJ}
{(\lbar{i}^{*,t} (\varphi_J), \lbar{j}^{*,t} (\varphi_I))}_t
\\&\quad+ \sum_{|I|=q} \sum_{i \in I} ([\lbar{i}, \lbar{i}^{*,t}] \varphi_I, \varphi_I)_t
+ O \Big(|||\varphi|||_t^2 +(\sum_{|I|=q} \sum_i
|||\lbar{i}\varphi_I|||_t^2)^\half \, |||\varphi|||_t \Big)
\end{split}
\end{equation*}
because $i \, = \, j$ forces $I \, = \, J,$ but if $i \neq j$ then
$\eps^{jK}_I \, = \, \eps^{ijK}_{iI}$ and $\eps^{iK}_J \, = \,
\eps^{jiK}_{jJ} \, = \, - \eps^{ijK}_{jJ},$ so $\eps^{jK}_I
\eps^{iK}_J \, = \, - \eps^{iI}_{jJ}.$ Now, we add up the
expressions for $|||\dbarb \varphi|||^2_t$ and $|||\adp
\varphi|||^2_t,$ unravel the commutators, and proceed as in
Lemma~\ref{qbtp} to finish this proof. \qed

\bigskip
\noindent Let us see what $\qbtm{\, \cdot \,}$ works out to be in this case.

\smallskip
\newtheorem{qbtmcomp0q}[qbtpcomp0q]{Lemma}
\begin{qbtmcomp0q}
Let $\varphi$ be a $(0,q)$ form \label{qbtm0q} supported in $U',$
$\varphi \in Dom(\dbarb) \cap Dom(\ad).$ If $[L_i, \lbar{j}]=
c_{ij}T + \sum_k a^k_{ij}\lbar{k} + \sum_k b^k_{ij} \lbar{k}^{*,
\, -t} + g_{ij} -t \, e_{ij},$ where $a^k_{ij}$, $b^k_{ij}$,
$g_{ij},$ and $e_{ij}$ are \smooth functions independent of $t$
and $c_{ij}$ are the coefficients of the Levi form, then there
exists $1 \gg \eps" > 0$ such that
\begin{equation*}
\begin{split}
\qbtm{\varphi} &\geq (1-\eps")\, \sum_{|I|=q} \sum_i
|||\lbar{i}^{*,-t} \varphi_I|||^2_{-t} - \sum_{|I|=q} \sum_{i
\notin I} \Re \{(c_{ii} T \varphi_I, \varphi_I)_{-t} \}\\&\quad  -
\sum_{|I|=|J|=q} \sum_{i,j \atop {i \neq j}} \eps^{iI}_{jJ} \, \Re
\{(c_{ij} T \varphi_J, \varphi_I)_{-t} \}+ t \sum_{|I|=q} \sum_{i
\notin I} \Re \{ ( \lbar{i} L_i (\lambda) \varphi_I,
\varphi_I)_{-t} \}\\&\quad+ t \sum_{|I|=|J|=q} \sum_{i,j \atop {i
\neq j}} \eps^{iI}_{jJ}\, \Re \{ (\lbar{j} L_i (\lambda)
\varphi_J, \varphi_I)_{-t} \}+ t \sum_{|I|=q} \sum_{i \notin I}
\Re \{ ( e_{ii} \varphi_I, \varphi_I)_{-t} \}\\&\quad
+t\sum_{|I|=|J|=q} \sum_{i,j \atop {i \neq j}} \eps^{iI}_{jJ}\,
\Re \{(e_{ij} \varphi_J, \varphi_I)_{-t} \} +
O(|||\varphi|||^2_{-t} ).
\end{split}
\end{equation*}
\end{qbtmcomp0q}

\smallskip
\noindent {\bf Proof:} By using the same proof as in the previous lemma with $t$ replaced by $-t$ and then by converting
all the $\lbar{i}$ derivatives into $\lbar{i}^{*,-t}$ derivatives, we obtain:
\begin{equation*}
\begin{split}
\qbtm{\varphi} &= \sum_{|I|=q} \sum_i |||\lbar{i}^{*,-t}
\varphi_I|||^2_{-t} + \sum_{|I|=q} \sum_{i \notin I}
([\lbar{i}^{*,-t}, \lbar{i}] \varphi_I, \varphi_I)_{-t} \\&\quad +
\sum_{|I|=|J|=q} \sum_{i,j \atop {i \neq j}} \eps^{iI}_{jJ}
([\lbar{i}^{*,-t}, \lbar{j}] \varphi_J, \varphi_I)_{-t}
\\&\quad+ O \Big(|||\varphi|||^2_{-t}+\big(\sum_{|I|=q} \sum_i |||\lbar{i}^{*,-t}
\varphi_I|||^2_{-t} \big)^\half \, |||\varphi|||_{-t}  \Big).
\end{split}
\end{equation*}
We then unravel the commutators and proceed as above. \qed

\medskip
\noindent These expressions for $\qbtp{\, \cdot \,}$ and $\qbtm{\, \cdot \,}$ for $(0,q)$ forms may seem different from
their counterparts for $(0,1)$ forms, but the difference only rests in two Hermitian matrices that behave the same
way for any level of forms. The following two lemmas will elucidate this observation:

\smallskip
\newtheorem{Hermpluslemma}[qbtpcomp0q]{Lemma}
\begin{Hermpluslemma}
Let $(g_{ij})$ be a Hermitian matrix and let $1 \leq q \leq n-2.$ \label{Hermplus} Then the ${n-1 \choose q}$ by
${n-1 \choose q}$ matrix $(G^{\, q}_{IJ})$ given by
\begin{equation*}
\begin{aligned}
&G^{\, q}_{II}= \sum_{i \in I} g_{ii} \\ &G^{\, q}_{IJ} = - \sum_{i,j \atop {i \neq j}} \eps^{iI}_{jJ} \, g_{ji} \quad
if \: I \neq J,
\end{aligned}
\end{equation*}
where $I$ and $J$ are multi-indices, $|I| \, = \, |J| \, = \, q,$ is also Hermitian. Moreover, the eigenvalues of
$(G^{\, q}_{IJ})$ are sums of the eigenvalues of $(g_{ij})$ taken $q$ at a time, so $(G^{\, q}_{IJ})$
is positive definite if $(g_{ij})$ is positive definite and positive semi-definite if $(g_{ij})$ is positive
semi-definite.
\end{Hermpluslemma}

\smallskip
\noindent {\bf Proof:} Given its definition, $(G^{\, q}_{IJ})$ is clearly Hermitian. Also, notice that
$(G^{\, 1}_{IJ}) \, = \, (g_{ij}).$ Now, since $(g_{ij})$ is Hermitian, we can diagonalize it at each point, and it is
easy to see that this diagonalizes $(G^{\, q}_{IJ})$ as well, and that the eigenvalues of $(G^{\, q}_{IJ})$ are sums
of the eigenvalues of $(g_{ij})$ taken $q$ at a time. \qed

\medskip
\newtheorem{Hermminuslemma}[qbtpcomp0q]{Lemma}
\begin{Hermminuslemma}
Let $(h_{ij})$ be a Hermitian matrix and let $1 \leq q \leq n-2.$ \label{Hermminus} Then the ${n-1 \choose q}$ by
${n-1 \choose q}$ matrix $(H^{\, q}_{IJ})$ given by
\begin{equation*}
\begin{aligned}
&H^{\, q}_{II}= \sum_{i \notin I} h_{ii} \\ &H^{\, q}_{IJ} =  \sum_{i,j \atop {i \neq j}} \eps^{iI}_{jJ} \, h_{ji} \quad
if \: I \neq J,
\end{aligned}
\end{equation*}
where $I$ and $J$ are multi-indices, $|I| \, = \, |J| \, = \, q,$ is also Hermitian. Moreover, the eigenvalues of
$(H^{\, q}_{IJ})$ are sums of the eigenvalues of $(h_{ij})$ taken $n-1-q$ at a time, so $(H^{\, q}_{IJ})$
is positive definite if $(h_{ij})$ is positive definite and $n-1-q >0;$ $(H^{\, q}_{IJ})$
is positive semi-definite if $(h_{ij})$ is positive semi-definite, regardless of the value of $n.$
\end{Hermminuslemma}

\smallskip
\noindent {\bf Proof:} It is proven the same way as the previous lemma. We note that
$$(H^{\, 1}_{IJ}) \, = \, \Big(\delta_{ij} \sum_i h_{ii} - h_{ij} \Big),$$ so Lemma~\ref{linalg} is
a special case of this one.   \qed

\bigskip
Armed with these two lemmas and having noticed that the proof of Lemma~\ref{Gardcp} is unaffected
by replacing the Levi form with some other positive semi-definite Hermitian form, we can easily see that
the Lemmas~\ref{Gardcp} and ~\ref{Gardcm} generalize to $(0,q)$ forms with minor modifications of notation.
It follows that the same holds for Propositions~\ref{localplus} and ~\ref{localminus}, since the covering of $M$
was chosen so that $q \, \eps_G \ll 1$ as noted at the beginning of this section. All these results then imply
the analogue of Proposition~\ref{global}, namely the same main estimate for $(0,q)$ forms as for $(0,1)$ forms up
to minor modifications of notation:
\medskip
\newtheorem{mainestimate0q}[qbtpcomp0q]{Proposition}
\begin{mainestimate0q}
Let $1 \leq q \leq n-2$ and $\varphi$ be a $(0,q)$ form supported on $M$, a compact,
orientable, weakly pseudoconvex \label{global0q} CR manifold of
dimension at least $5$, and let $M$ be endowed with a strongly CR
plurisubharmonic function $\lambda$ as described in
Section~\ref{CRpsh}. If $\varphi$ satisfies $\varphi \in Dom(\dbarb)
\cap Dom(\ad),$ then there exist constants $K,$ $K_t,$ and
$K'_t,$ and a positive number $T_0$ such that for any $t \geq
T_0$,
\begin{equation*}
K \qbtg {\varphi}+K_t \sum_\nu \sum_{|I|=q} \:{||\dze_\nu
\potdg {\nu} \ze_\nu \varphi_I^\nu||}^2_0 +K'_t \,||\varphi||^2_{-1}
\geq  t {\langle |\varphi|\rangle}^2_t.
\end{equation*}
\end{mainestimate0q}

\bigskip
The different expression for $\qbto{\, \cdot \,}$ that one obtains for $(0,q)$ forms by plugging $t=0$ into the
expression for $\qbtp{\, \cdot \,}$ from Lemma~\ref{qbtp0q} does not affect the proof of Lemma~\ref{Sob1loc}, hence
we obtain the same result for $(0,q)$ forms; therefore, the same is true for Lemma~\ref{Sob1global} which is based on it.
This version of Lemma~\ref{Sob1global} then allows us to extend to $(0,q)$ forms all the results in Section~\ref{regdbar}
pertaining to $(0,1)$ forms, thus completing the proof of the Main Theorem.

\bigskip
\section{The G{\aa}rding Inequality}
\label{Gardsect}

\medskip
\medskip
The following version of the sharp G\aa rding inequality can be found either in \cite{LaxNirenberg} or derived
from the scalar version in \cite{Hormander}:

\smallskip
\newtheorem{1Garding}{Theorem}[section]
\begin{1Garding}
If $P$ is a matrix first order pseudodifferential operator with positive semi-definite symbol, then there
exists a constant $C$ such that $\sum_{ij} Re \{(P_{ij} u_i, u_j)
\} \geq -C ||u||^2$ for any vector of smooth functions $u$.
\end{1Garding}

\medskip
\smallskip
\newtheorem{2Garding}[1Garding]{Lemma}
\begin{2Garding}
Let $R$ be a first order pseudodifferential operator such that $\sigma(R) \geq
\kappa$ \label{Gardsecond} where $\kappa$ is some positive
constant and let $(h_{ij})$ be a positive semi-definite matrix.
Then there exists a constant $C$ such that $\sum_{ij} Re \{(h_{ij}
R u_i, u_j)\} \geq \kappa \sum_{ij} (h_{ij}u_i, u_j)-C ||u||$ for
any vector of smooth functions $u$.
\end{2Garding}
\noindent {\bf{Proof:}} Apply the previous theorem with $P$, where
$P_{ij} \, = \, h_{ij}(R - \kappa)$ which is positive
semi-definite, given the hypotheses. Therefore, there must exist a
constant $C$ so that $$\sum_{ij} Re \{(P_{ij} u_i, u_j) \} \geq -C
||u||$$  $$\Updownarrow$$ $$ \sum_{ij} Re \{((h_{ij} R - \kappa
h_{ij}) u_i, u_j) \} \, = \,\sum_{ij} Re \{(h_{ij} R u_i, u_j) \}
- \kappa \sum_{ij} (h_{ij} u_i, u_j) \geq -C ||u||,$$ because
$(h_{ij})$ is positive semi-definite, so $\sum_{ij} (h_{ij} u_i,
u_j)$ is real-valued. Rearranging this inequality we get precisely
the conclusion of the lemma. \qed

\bigskip
\section{Main Estimate Computation Details}
\label{Compdetails}

\medskip
\medskip
 The first lemma shows that compositions of pseudodifferential
operators whose symbols have disjoint support are smoothing
operators, hence easily controlled by the error of order zero.
\smallskip
\newtheorem{disjointlemma}{Lemma}[section]
\begin{disjointlemma}
Let $P_1$ and $P_2$ be two pseudodifferential operators with
symbols $\sigma(P_1) \, = \, p_1(x,\xi)$ and $\sigma(P_2) \, = \,
p_2(x,\xi)$ respectively, \label{disjlemma} which satisfy $supp \,
(p_1 (x,\xi)) \cap supp \,(p_2 (x,\xi)) \, = \, \emptyset$. Then
the composition $P_1 P_2$ is a smoothing operator.
\end{disjointlemma}
\noindent {\bf{Proof:}} This is an elementary fact which follows
from the pseudodifferential operator theory in \cite{Treves}. \qed

\medskip\noindent Relying upon this result, we now give the details of the computation of the adjoint $\adt,$
Lemma~\ref{adclaim}, and that of the energy form,
Lemma~\ref{globalenergy}, which were omitted in
Section~\ref{microlocsection}.

\medskip \noindent {\bf Proof of Lemma~\ref{adclaim}:} The proof
proceeds in two steps. First, we exploit $$\tnorm{f}{\adt \,
\varphi}=\tnorm{\dbarb \, f}{\varphi}$$ in order to obtain an
expression for $\adt.$ Second, we compare this expression with the
one claimed for $\adt$ and compute the error $E_t$ to show it has
the desired form.

\smallskip \noindent {\bf Step 1:}\begin{equation*}
\begin{split}
\tnorm{f}{\adt \varphi} &= \tnorm {\dbarb f}{\varphi}\\ &=
\sum_\nu \: ({(\dze_\nu \pptg {\nu} \ze_\nu {\dbarb f^\nu},
\dze_\nu \pptg {\nu} \ze_\nu {\varphi^\nu})}_t+ {(\dze_\nu \potg
{\nu}\ze_\nu {\dbarb f^\nu}, \dze_\nu \potg {\nu} \ze_\nu {\varphi^\nu})}_0 \\
& \quad+ {(\dze_\nu \pmtg {\nu} \ze_\nu {\dbarb f^\nu}, \dze_\nu
\pmtg {\nu} \ze_\nu {\varphi^\nu})}_{-t})
\end{split}
\end{equation*}
Let us look at ${(\dze_\nu \pptg {\nu} \ze_\nu {\dbarb f^\nu},
\dze_\nu \pptg {\nu} \ze_\nu {\varphi^\nu})}_t$ since the other
two terms will behave in a similar manner:
\begin{equation*}
\begin{split}
{(\dze_\nu \pptg {\nu} \ze_\nu {\dbarb f^\nu}, \dze_\nu \pptg
{\nu} \ze_\nu {\varphi^\nu})}_t &={(\dbarb \, \dze_\nu \pptg {\nu}
\ze_\nu { f^\nu}, \dze_\nu \pptg {\nu} \ze_\nu {\varphi^\nu})}_t +
{([\dze_\nu \pptg {\nu} \ze_\nu, \dbarb] \dze_\nu f^\nu, \dze_\nu
\pptg {\nu} \ze_\nu {\varphi^\nu})}_t
\\&= {(\dze_\nu \pptg {\nu} \ze_\nu { f^\nu},(\ad -t [\ad, \lambda])
\dze_\nu \pptg {\nu} \ze_\nu {\varphi^\nu})}_t  \\ &\quad
+{([\dze_\nu \pptg {\nu} \ze_\nu, \dbarb] \dze_\nu f^\nu, \dze_\nu
\pptg {\nu} \ze_\nu {\varphi^\nu})}_t \\ &={(\dze_\nu \pptg {\nu}
\ze_\nu { f^\nu}, \dze_\nu \pptg {\nu} \ze_\nu (\ad -t [\ad,
\lambda]){\varphi^\nu})}_t \\ &\quad + {( { f^\nu}, \ze_\nu (\pptg
{\nu})^* \dze_\nu [\ad -t [\ad, \lambda], \dze_\nu \pptg {\nu}
\ze_\nu] \dze_\nu {\varphi^\nu})}_t \\ &\quad +
{([e^{-t\lambda},\pptg {\nu}] \ze_\nu { f^\nu},\dze_\nu [\ad -t
[\ad, \lambda], \dze_\nu \pptg {\nu} \ze_\nu] \dze_\nu
{\varphi^\nu})}_0\\&\quad + {(f^\nu, \dze_\nu [\dze_\nu \pptg
{\nu} \ze_\nu, \dbarb]^* \dze_\nu \pptg {\nu} \ze_\nu
{\varphi^\nu})}_t  \\&\quad + {([e^{-t \lambda},[\dze_\nu \pptg
{\nu} \ze_\nu, \dbarb]] \dze_\nu f^\nu, \dze_\nu \pptg {\nu}
\ze_\nu {\varphi^\nu})}_0
\end{split}
\end{equation*}
where $[\, \cdot \, , \, \cdot \,]^*$ is the adjoint of the
commutator $[\, \cdot \, , \, \cdot \,].$ Then,
\begin{equation*}
\begin{split}
{(\dze_\nu \pptg {\nu} \ze_\nu {\dbarb f^\nu}, \dze_\nu \pptg
{\nu} \ze_\nu {\varphi^\nu})}_t &={(\dze_\nu \pptg {\nu} \ze_\nu {
f^\nu}, \dze_\nu \pptg {\nu} \ze_\nu (\ad -t [\ad,
\lambda]){\varphi^\nu})}_t \\ &\quad + {( { f^\nu}, \ze_\nu (\pptg
{\nu})^* \dze_\nu [\ad -t [\ad, \lambda], \dze_\nu \pptg {\nu}
\ze_\nu] \dze_\nu {\varphi^\nu})}_t
\\&\quad+{({f^\nu},\ze_\nu[e^{-t\lambda},\pptg {\nu}]^*\dze_\nu
[\ad -t [\ad, \lambda], \dze_\nu \pptg {\nu} \ze_\nu] \dze_\nu
{\varphi^\nu})}_0\\&\quad + {(f^\nu, \dze_\nu [\dze_\nu \pptg
{\nu} \ze_\nu, \dbarb]^* \dze_\nu \pptg {\nu} \ze_\nu
{\varphi^\nu})}_t
\\&\quad + {( f^\nu, \dze_\nu [e^{-t \lambda},[\dze_\nu \pptg {\nu} \ze_\nu,
\dbarb]]^*\dze_\nu \pptg {\nu} \ze_\nu {\varphi^\nu})}_0
\end{split}
\end{equation*}

\noindent Let  $${( { f^\nu}, \ze_\nu (\pptg {\nu})^* \dze_\nu
[\ad -t [\ad, \lambda], \dze_\nu \pptg {\nu} \ze_\nu] \dze_\nu
{\varphi^\nu})}_t =I.$$ Since $\sum_\mu \ze^2_\mu \, = \, 1$ and
$(\pptg {\mu})^* \pptg{\mu}  + (\potg {\mu})^* \potg {\mu} +
(\pmtg {\mu})^* \pmtg {\mu} = Id$ we have,
\begin{equation*}
\begin{split}
I &= \sum_\mu {( \pptg{\mu} \ze_\mu { f^\mu},\pptg{\mu} \ze_\mu
\ze_\nu (\pptg {\nu})^* \dze_\nu [\ad -t [\ad, \lambda], \dze_\nu
\pptg {\nu} \ze_\nu] \dze_\nu {\varphi^\mu})}_t\\ &\quad +
\sum_\mu {( \potg{\mu} \ze_\mu { f^\mu},\potg{\mu} \ze_\mu \ze_\nu
(\pptg {\nu})^* \dze_\nu [\ad -t [\ad, \lambda], \dze_\nu \pptg
{\nu} \ze_\nu] \dze_\nu {\varphi^\mu})}_t\\ &\quad + \sum_\mu {(
\pmtg{\mu} \ze_\mu { f^\mu},\pmtg{\mu} \ze_\mu \ze_\nu (\pptg
{\nu})^* \dze_\nu [\ad -t [\ad, \lambda], \dze_\nu \pptg {\nu}
\ze_\nu] \dze_\nu {\varphi^\mu})}_t
\\&= \sum_\mu {(\dze_\mu \pptg{\mu} \ze_\mu { f^\mu},\dze_\mu \pptg{\mu} \ze_\mu
\ze_\nu (\pptg {\nu})^* \dze_\nu [\ad -t [\ad, \lambda], \dze_\nu
\pptg {\nu} \ze_\nu] \dze_\nu {\varphi^\mu})}_t\\ &\quad +
\sum_\mu {(   { f^\mu},\ze_\mu(\potg{\mu})^* \dze_\mu e^{-t
\lambda}\dze_\mu \potg{\mu} \ze_\mu\ze_\nu (\pptg {\nu})^*
\dze_\nu [\ad -t [\ad, \lambda], \dze_\nu \pptg {\nu} \ze_\nu]
\dze_\nu {\varphi^\mu})}_0\\ &\quad + lower \: order \: terms
\end{split}
\end{equation*}
because $\dze_\mu$ dominates $\ze_\mu$ and the principal symbols
of $\pmtg{\mu}$ and $\pptg{\nu}$ have disjoint supports in $\xi$
space after transferring one to the other via some diffeomorphism
taking $U_\nu$ to $U_\mu$. In other words, all the other terms are
of order $-1$ or lower.

\smallskip\noindent Similarly, let $${(f^\nu, \dze_\nu [\dze_\nu \pptg {\nu} \ze_\nu,
\dbarb]^* \dze_\nu \pptg {\nu} \ze_\nu {\varphi^\nu})}_t = II,$$
then
\begin{equation*}
\begin{split}
II &=\sum_\mu {(\dze_\mu \pptg{\mu} \ze_\mu { f^\mu},\dze_\mu
\pptg{\mu} \ze_\mu \dze_\nu [\dze_\nu \pptg {\nu} \ze_\nu,
\dbarb]^* \dze_\nu \pptg {\nu} \ze_\nu {\varphi^\mu})}_t\\ &\quad
+ \sum_\mu {(  { f^\mu},\ze_\mu(\potg{\mu})^* \dze_\mu  e^{-t
\lambda}\dze_\mu \potg{\mu} \ze_\mu \dze_\nu [\dze_\nu \pptg {\nu}
\ze_\nu, \dbarb]^* \dze_\nu \pptg {\nu} \ze_\nu
{\varphi^\mu})}_0\\ &\quad +lower \: order \: terms.
\end{split}
\end{equation*}
Let ${(\dze_\nu \pptg {\nu} \ze_\nu {\dbarb f^\nu}, \dze_\nu \pptg
{\nu} \ze_\nu {\varphi^\nu})}_t \, = \, III.$ The above means
\begin{equation*}
\begin{split}
III &= {(\dze_\nu \pptg {\nu} \ze_\nu { f^\nu}, \dze_\nu \pptg
{\nu} \ze_\nu (\ad -t [\ad, \lambda]){\varphi^\nu})}_t\\ &\quad
+\sum_\mu \:\Big(\, {(\dze_\mu \pptg{\mu} \ze_\mu {
f^\mu},\dze_\mu \pptg{\mu} \ze_\mu \dze_\nu [\dze_\nu \pptg {\nu}
\ze_\nu, \dbarb]^* \dze_\nu \pptg {\nu} \ze_\nu {\varphi^\mu})}_t
\\ &\quad + {(\dze_\mu \pptg{\mu} \ze_\mu { f^\mu},\dze_\mu
\pptg{\mu} \ze_\mu \ze_\nu (\pptg {\nu})^* \dze_\nu [\ad -t [\ad,
\lambda], \dze_\nu \pptg {\nu} \ze_\nu] \dze_\nu
{\varphi^\mu})}_t\\&\quad+{(   { f^\mu},\ze_\mu(\potg{\mu})^*
\dze_\mu e^{-t \lambda}\dze_\mu \potg{\mu} \ze_\mu\ze_\nu (\pptg
{\nu})^* \dze_\nu [\ad -t [\ad, \lambda], \dze_\nu \pptg {\nu}
\ze_\nu] \dze_\nu {\varphi^\mu})}_0\\&\quad+{(  {
f^\mu},\ze_\mu(\potg{\mu})^* \dze_\mu  e^{-t \lambda}\dze_\mu
\potg{\mu} \ze_\mu \dze_\nu [\dze_\nu \pptg {\nu} \ze_\nu,
\dbarb]^* \dze_\nu \pptg {\nu}
\ze_\nu {\varphi^\mu})}_0\,\Big)\\
&\quad+{({f^\nu},\ze_\nu[e^{-t\lambda},\pptg {\nu}]^*\dze_\nu [\ad
-t [\ad, \lambda], \dze_\nu \pptg {\nu} \ze_\nu] \dze_\nu
{\varphi^\nu})}_0
\\&\quad + {( f^\nu, \dze_\nu [e^{-t \lambda},[\dze_\nu \pptg {\nu} \ze_\nu,
\dbarb]]^*\dze_\nu \pptg {\nu} \ze_\nu {\varphi^\nu})}_0\\ &\quad
+ lower \: order \: terms
\end{split}
\end{equation*}
Similarly, let $ {(\dze_\nu \pmtg {\nu} \ze_\nu {\dbarb f^\nu},
\dze_\nu \pmtg {\nu} \ze_\nu {\varphi^\nu})}_{-t} \, = \, IV,$
then
\begin{equation*}
\begin{split}
IV &= {(\dze_\nu \pmtg {\nu} \ze_\nu { f^\nu}, \dze_\nu \pmtg
{\nu} \ze_\nu (\ad +t [\ad,
\lambda]){\varphi^\nu})}_{-t}\\
&\quad +\sum_\mu \:\Big(\,{(\dze_\mu \pmtg{\mu} \ze_\mu {
f^\mu},\dze_\mu \pmtg{\mu} \ze_\mu \dze_\nu [\dze_\nu \pmtg {\nu}
\ze_\nu,
\dbarb]^* \dze_\nu \pmtg {\nu} \ze_\nu {\varphi^\mu})}_{-t} \\
&\quad + {(\dze_\mu \pmtg{\mu} \ze_\mu { f^\mu},\dze_\mu
\pmtg{\mu} \ze_\mu \ze_\nu (\pmtg {\nu})^* \dze_\nu [\ad +t [\ad,
\lambda], \dze_\nu \pmtg {\nu} \ze_\nu] \dze_\nu
{\varphi^\mu})}_{-t}\\
&\quad +  {( { f^\mu}, \ze_\mu (\potg{\mu})^*\dze_\mu e^{t
\lambda}\dze_\mu \potg{\mu} \ze_\mu  \dze_\nu [\dze_\nu \pmtg
{\nu} \ze_\nu, \dbarb]^* \dze_\nu \pmtg {\nu} \ze_\nu
{\varphi^\mu})}_0\\ &\quad +  {( { f^\mu},\ze_\mu (\potg{\mu})^*
\dze_\mu e^{t \lambda}\dze_\mu \potg{\mu} \ze_\mu \ze_\nu (\pmtg
{\nu})^* \dze_\nu [\ad +t [\ad, \lambda],
\dze_\nu \pmtg {\nu} \ze_\nu] \dze_\nu {\varphi^\mu})}_0\,\Big)\\
&\quad+{({f^\nu},\ze_\nu[e^{t\lambda},\pmtg {\nu}]^*\dze_\nu [\ad
+t [\ad, \lambda], \dze_\nu \pmtg {\nu} \ze_\nu] \dze_\nu
{\varphi^\nu})}_0
\\&\quad + {( f^\nu, \dze_\nu [e^{t \lambda},[\dze_\nu \pmtg {\nu} \ze_\nu,
\dbarb]]^*\dze_\nu \pmtg {\nu} \ze_\nu {\varphi^\nu})}_0\\
&\quad + lower \: order \: terms.
\end{split}
\end{equation*}
For the zero norm, the following suffices:
\begin{equation*}
\begin{split}
{(\dze_\nu \potg {\nu} \ze_\nu {\dbarb f^\nu}, \dze_\nu \potg
{\nu} \ze_\nu {\varphi^\nu})}_0 &= {(\dze_\nu \potg {\nu} \ze_\nu
{ f^\nu}, \dze_\nu \potg {\nu} \ze_\nu \ad {\varphi^\nu})}_0  +
{(\dze_\nu \potg {\nu} \ze_\nu { f^\nu},[\ad, \dze_\nu \potg {\nu}
\ze_\nu] \dze_\nu {\varphi^\nu})}_0
\\&\quad + {(f^\nu, \dze_\nu [\dze_\nu \potg {\nu} \ze_\nu,
\dbarb]^* \dze_\nu \potg {\nu} \ze_\nu {\varphi^\nu})}_0 \\&=
{(\dze_\nu \potg {\nu} \ze_\nu { f^\nu}, \dze_\nu \potg {\nu}
\ze_\nu \ad {\varphi^\nu})}_0  + {(  f^\nu, \ze_\nu(\potg {\nu})^*
\dze_\nu[\ad, \dze_\nu \potg {\nu} \ze_\nu] \dze_\nu
{\varphi^\nu})}_0 \\&\quad + {(f^\nu, \dze_\nu [\dze_\nu \potg
{\nu} \ze_\nu, \dbarb]^* \dze_\nu \potg {\nu} \ze_\nu
{\varphi^\nu})}_0
\end{split}
\end{equation*}
Finally we get,
\begin{equation*}
\begin{split}
\tnorm{f}{\adt \varphi} &=\sum_\nu \: \Big(\,{(\dze_\nu \pptg{\nu}
\ze_\nu f,\dze_\nu \pptg {\nu} \ze_\nu (\ad -t [\ad ,\lambda])
\varphi_\nu)}_t+ {(\dze_\nu \potg{\nu} \ze_\nu f,\dze_\nu \potg
{\nu} \ze_\nu \ad \varphi_\nu)}_0\\ & \quad+ {(\dze_\nu \pmtg{\nu}
\ze_\nu f,\dze_\nu \pmtg {\nu} \ze_\nu (\ad + t [\ad ,\lambda])
\varphi_\nu)}_{-t}\,\Big)\\ &\quad +\sum_\nu \sum_{\mu \atop U_\nu
\cap U_\mu \neq \emptyset} \: \Big(\, {(\dze_\mu \pptg{\mu}
\ze_\mu { f^\mu},\dze_\mu \pptg{\mu} \ze_\mu \dze_\nu [\dze_\nu
\pptg {\nu} \ze_\nu, \dbarb]^* \dze_\nu \pptg {\nu} \ze_\nu
{\varphi^\mu})}_t\\ &\quad + {(\dze_\mu \pptg{\mu} \ze_\mu {
f^\mu},\dze_\mu \pptg{\mu} \ze_\mu \ze_\nu (\pptg {\nu})^*
\dze_\nu [\ad -t [\ad, \lambda], \dze_\nu \pptg {\nu} \ze_\nu]
\dze_\nu {\varphi^\mu})}_t\\ &\quad + {(\dze_\mu \pmtg{\mu}
\ze_\mu { f^\mu},\dze_\mu \pmtg{\mu} \ze_\mu \ze_\nu (\pmtg
{\nu})^* \dze_\nu [\ad +t [\ad, \lambda], \dze_\nu \pmtg {\nu}
\ze_\nu] \dze_\nu {\varphi^\mu})}_{-t}\\&\quad + {(\dze_\mu
\pmtg{\mu} \ze_\mu { f^\mu},\dze_\mu \pmtg{\mu} \ze_\mu \dze_\nu
[\dze_\nu \pmtg {\nu} \ze_\nu, \dbarb]^* \dze_\nu \pmtg {\nu}
\ze_\nu {\varphi^\mu})}_{-t}\,\Big)\\&\quad + \sum_\nu \: \Big(\,
{(  f^\nu, \ze_\nu(\potg {\nu})^* \dze_\nu[\ad, \dze_\nu \potg
{\nu} \ze_\nu] \dze_\nu {\varphi^\nu})}_0  + {(f^\nu, \dze_\nu
[\dze_\nu \potg {\nu} \ze_\nu, \dbarb]^* \dze_\nu \potg {\nu}
\ze_\nu
{\varphi^\nu})}_0\,\Big)\\
&\quad+\sum_\nu\:\Big(\,{({f^\nu},\ze_\nu[e^{-t\lambda},\pptg
{\nu}]^*\dze_\nu [\ad -t [\ad, \lambda], \dze_\nu \pptg {\nu}
\ze_\nu] \dze_\nu {\varphi^\nu})}_0 \\&\quad+ {( f^\nu, \dze_\nu
[e^{-t \lambda},[\dze_\nu \pptg {\nu} \ze_\nu, \dbarb]]^*\dze_\nu
\pptg
{\nu} \ze_\nu {\varphi^\nu})}_0\\
&\quad+{({f^\nu},\ze_\nu[e^{t\lambda},\pmtg {\nu}]^*\dze_\nu [\ad
+t [\ad, \lambda], \dze_\nu \pmtg {\nu} \ze_\nu] \dze_\nu
{\varphi^\nu})}_0 \\&\quad+ {( f^\nu, \dze_\nu [e^{t
\lambda},[\dze_\nu \pmtg {\nu} \ze_\nu, \dbarb]]^*\dze_\nu \pmtg
{\nu} \ze_\nu {\varphi^\nu})}_0\,\Big)\\ &\quad +\sum_\nu
\sum_{\mu \atop U_\nu \cap U_\mu \neq \emptyset} \: \Big(\,{(  {
f^\mu},\ze_\mu(\potg{\mu})^* \dze_\mu e^{-t \lambda}\dze_\mu
\potg{\mu} \ze_\mu \dze_\nu [\dze_\nu \pptg {\nu} \ze_\nu,
\dbarb]^* \dze_\nu \pptg {\nu} \ze_\nu
{\varphi^\mu})}_0\\&\quad+{(   { f^\mu},\ze_\mu(\potg{\mu})^*
\dze_\mu e^{-t \lambda}\dze_\mu \potg{\mu} \ze_\mu\ze_\nu (\pptg
{\nu})^* \dze_\nu [\ad -t [\ad, \lambda], \dze_\nu \pptg {\nu}
\ze_\nu] \dze_\nu
{\varphi^\mu})}_0 \\
&\quad +  {( { f^\mu}, \ze_\mu (\potg{\mu})^*\dze_\mu e^{t
\lambda}\dze_\mu \potg{\mu} \ze_\mu  \dze_\nu [\dze_\nu \pmtg
{\nu} \ze_\nu, \dbarb]^* \dze_\nu \pmtg {\nu} \ze_\nu
{\varphi^\mu})}_0\\ &\quad +  {( { f^\mu},\ze_\mu (\potg{\mu})^*
\dze_\mu e^{t \lambda}\dze_\mu \potg{\mu} \ze_\mu \ze_\nu (\pmtg
{\nu})^* \dze_\nu [\ad +t [\ad, \lambda], \dze_\nu \pmtg {\nu}
\ze_\nu] \dze_\nu {\varphi^\mu})}_0\,\Big)\\&\quad + (f, S_{t,+}
\, \varphi)_0 + (f, S_{t,-} \, \varphi)_0
\end{split}
\end{equation*}
since the lower order errors for the $t$ and $-t$ norms can be
rewritten so that they can be represented by two $t$ dependent
pseudodifferential operators of order $-1$, $S_{t,+}$ and
$S_{t,-}.$ On each $U_\rho$, we define operators $\pptdg{\rho}$
and $\pmtdg{\rho}$ fulfilling Property $\dag$ on page
\pageref{propdagger}.

\smallskip \noindent {\bf Step 2:} The expression for $\adt$ is so long that it cannot all
be handled at the same time, so we split it as follows:

\smallskip \noindent Let $$A = \ad -t\sum_\rho \ze^2_\rho \pptdg{\rho} [\ad , \lambda] +
t \sum_\rho \ze^2_\rho \pmtdg{\rho} [\ad, \lambda], $$ $$B =
\sum_\rho \: (\dze_\rho [\dze_\rho \pptg {\rho} \ze_\rho,
\dbarb]^* \dze_\rho \pptg {\rho} \ze_\rho +\ze_\rho (\pptg
{\rho})^* \dze_\rho [\ad -t [\ad, \lambda], \dze_\rho \pptg {\rho}
\ze_\rho] \dze_\rho),$$ and $$C = \sum_\rho \: (\ze_\rho (\pmtg
{\rho})^* \dze_\rho [\ad +t [\ad, \lambda], \dze_\rho \pmtg {\rho}
\ze_\rho] \dze_\rho +\dze_\rho [\dze_\rho \pmtg {\rho} \ze_\rho,
\dbarb]^* \dze_\rho \pmtg {\rho} \ze_\rho).$$ First we compute
$\tnorm{f}{A \, \varphi}$:
\begin{equation*}
\begin{split}
\tnorm{f}{A \, \varphi} &=  \sum_\nu \: \bigg({\Big(\dze_\nu
\pptg{\nu} \ze_\nu f^\nu,\dze_\nu \pptg {\nu} \ze_\nu (\ad
-t\sum_\rho \ze^2_\rho \pptdg{\rho} [\ad , \lambda] + t\sum_\rho
 \ze^2_\rho \pmtdg{\rho} [\ad, \lambda]) \, \varphi^\nu \Big)}_t\\
& \quad+ {\Big(\dze_\nu \potg{\nu} \ze_\nu f^\nu,\dze_\nu \potg
{\nu} \ze_\nu (\ad -t \sum_\rho\ze^2_\rho \pptdg{\rho} [\ad ,
\lambda] +t\sum_\rho\ze^2_\rho \pmtdg{\rho} [\ad, \lambda]) \varphi^\nu \Big)}_0\\
& \quad+ {\Big(\dze_\nu \pmtg{\nu} \ze_\nu f^\nu,\dze_\nu \pmtg
{\nu} \ze_\nu (\ad -t \sum_\rho\ze^2_\rho \pptdg{\rho} [\ad ,
\lambda] + t\sum_\rho
 \ze^2_\rho \pmtdg{\rho} [\ad, \lambda]) \varphi^\nu
\Big)}_{-t}\,\bigg)
\end{split}
\end{equation*}
Rearranging the terms and commuting $\pptdg{\rho}$ and
$\pmtdg{\rho}$ outside of $\ze_\nu \ze^2_\rho,$ we obtain
\begin{equation*}
\begin{split}
\tnorm{f}{A \, \varphi} &= \tnorm{f}{\ad \varphi} +  \sum_{\nu, \,
\rho} \: \Big(-t{(\dze_\nu \pptg{\nu} \ze_\nu f^\nu,\dze_\nu \pptg
{\nu} \pptdg{\rho} \ze_\nu
\ze^2_\rho [\ad,\lambda] \varphi^\nu)}_t \\
& \quad-t{(\dze_\nu \pptg{\nu} \ze_\nu f^\nu,\dze_\nu \pptg {\nu}
[\ze_\nu \ze^2_\rho, \pptdg{\rho}]  [\ad,\lambda]
\varphi^\nu)}_t\\&\quad+t{(\dze_\nu \pptg{\nu} \ze_\nu
f^\nu,\dze_\nu \pptg {\nu} \pmtdg{\rho} \ze_\nu
\ze^2_\rho  [\ad,\lambda] \varphi^\nu)}_t \\
& \quad+t{(\dze_\nu \pptg{\nu} \ze_\nu f^\nu,\dze_\nu \pptg {\nu}
[\ze_\nu \ze^2_\rho, \pmtdg{\rho}] [\ad,\lambda] \varphi^\nu)}_t
\\ & \quad+ {(\dze_\nu \potg{\nu} \ze_\nu f^\nu,\dze_\nu \potg {\nu}
\ze_\nu (-t\ze^2_\rho \pptdg{\rho}  [\ad,\lambda] + t\ze^2_\rho
\pmtdg{\rho}  [\ad,\lambda]) \varphi^\nu)}_0\\
& \quad-t{(\dze_\nu \pmtg{\nu} \ze_\nu f^\nu,\dze_\nu \pmtg {\nu}
\pptdg{\rho} \ze_\nu \ze^2_\rho  [\ad,\lambda]
\varphi^\nu)}_{-t}\\ & \quad-t {(\dze_\nu \pmtg{\nu} \ze_\nu
f^\nu,\dze_\nu \pmtg {\nu} [\ze_\nu \ze^2_\rho, \pptdg{\rho}]
[\ad,\lambda] \varphi^\nu)}_{-t}\\&\quad+t {(\dze_\nu \pmtg{\nu}
\ze_\nu f^\nu,\dze_\nu \pmtg {\nu} \pmtdg{\rho} \ze_\nu \ze^2_\rho
[\ad,\lambda] \varphi^\nu)}_{-t}\\ & \quad+t {(\dze_\nu \pmtg{\nu}
\ze_\nu f^\nu,\dze_\nu \pmtg {\nu} [\ze_\nu \ze^2_\rho,
\pmtdg{\rho}] [\ad,\lambda] \varphi^\nu)}_{-t}\,\Big)
\end{split}
\end{equation*}
Since the symbols of $\pptdg{\rho}$ and $\pmtdg{\rho}$ are
identically $1$ on $\cp_\rho$ and $\cm_\rho$ respectively, the
symbols of the commutators $[\ze_\nu \ze^2_\rho, \pptdg{\rho}]$
and $[\ze_\nu \ze^2_\rho, \pmtdg{\rho}]$ are supported in
$\co_\rho.$ Moreover, by part (i) of Property $\dag$ if $U_\nu
\cap U_\rho \neq \emptyset$ or by Lemma~\ref{disjlemma} otherwise,
$\pptg{\nu} \pmtdg{\rho},$ $\pmtg{\nu} \pptdg{\rho},$
$\pptg{\nu}[\ze_\nu \ze^2_\rho, \pmtdg{\rho}],$ and
$\pmtg{\nu}[\ze_\nu \ze^2_\rho, \pptdg{\rho}]$ are smoothing
operators. Also, by part (ii) of Property $\dag$, $\pptg {\nu}
\pptdg{\rho} \, = \, \pptg {\nu}$ up to lower order terms for all
$\rho$ such that $U_\nu \cap U_\rho \neq \emptyset$ and $\pptg
{\nu} \pptdg{\rho}$ is a smoothing operator by
Lemma~\ref{disjlemma} otherwise. Similarly, $\pmtg {\nu}
\pmtdg{\rho} \, = \, \pmtg {\nu}$ up to lower order terms for all
$\rho$ such that $U_\nu \cap U_\rho \neq \emptyset$ and $\pmtg
{\nu} \pmtdg{\rho}$ is a smoothing operator by
Lemma~\ref{disjlemma} otherwise. This means
\begin{equation*}
\begin{split}
\tnorm{f}{A \, \varphi} &= \tnorm{f}{\ad \varphi} +  \sum_\nu
\sum_{\rho \atop {U_\nu \cap U_\rho \neq \emptyset}} \:
\Big(-t{(\dze_\nu \pptg{\nu} \ze_\nu f^\nu,\dze_\nu \pptg {\nu}
\ze_\nu
\ze^2_\rho [\ad,\lambda] \varphi^\nu)}_t \\
& \quad+t {(\dze_\nu \pmtg{\nu} \ze_\nu f^\nu,\dze_\nu \pmtg
{\nu}\ze_\nu \ze^2_\rho [\ad,\lambda] \varphi^\nu)}_{-t}\\&\quad
-t{(\dze_\nu \pptg{\nu} \ze_\nu f^\nu,\dze_\nu \pptg {\nu}
[\ze_\nu \ze^2_\rho, \pptdg{\rho}]  [\ad,\lambda] \varphi^\nu)}_t
\\ & \quad+ {(\dze_\nu \potg{\nu} \ze_\nu f^\nu,\dze_\nu \potg {\nu}
\ze_\nu (-t\ze^2_\rho \pptdg{\rho}  [\ad,\lambda] + t\ze^2_\rho
\pmtdg{\rho} [\ad,\lambda]) \varphi^\nu)}_0\\ & \quad+t {(\dze_\nu
\pmtg{\nu} \ze_\nu f^\nu,\dze_\nu \pmtg {\nu} [\ze_\nu \ze^2_\rho,
\pmtdg{\rho}] [\ad,\lambda] \varphi^\nu)}_{-t}\,\Big)+ {(f, A_{-1}
\, \varphi)}_0,
\end{split}
\end{equation*}
where the lower order errors are rewritten so that they can be
represented by a pseudodifferential operator of order $-1$ denoted
by $A_{-1}.$ For the first two terms of the sum, pull $\sum_\rho$
inside in front of $\ze^2_\rho$ and notice that $\sum_\rho
\ze^2_\rho \, = \, 1,$ so $\tnorm{f}{A \, \varphi}$ equals the
following:
\begin{equation*}
\begin{split}
\tnorm{f}{A \, \varphi} &= \tnorm{f}{\ad \varphi} +  \sum_\nu
\:\Big(-t{(\dze_\nu \pptg{\nu} \ze_\nu f^\nu,\dze_\nu \pptg {\nu}
\ze_\nu [\ad,\lambda] \varphi^\nu)}_t\\&\quad +t {(\dze_\nu
\pmtg{\nu} \ze_\nu f^\nu,\dze_\nu \pmtg {\nu}\ze_\nu [\ad,\lambda]
\varphi^\nu)}_{-t}\,\Big)\\&\quad +\sum_\nu \sum_{\rho \atop
{U_\nu \cap U_\rho
\neq\emptyset}}\:\Big(-t{(f^\nu,\ze_\nu(\pptg{\nu})^*\dze^2_\nu
\pptg {\nu} [\ze_\nu \ze^2_\rho, \pptdg{\rho}] [\ad,\lambda]
\varphi^\nu)}_t\\ & \quad+ {(f^\nu,\ze_\nu(\potg{\nu})^*\dze^2_\nu
\potg {\nu} \ze_\nu (-t\ze^2_\rho \pptdg{\rho} [\ad,\lambda] +
t\ze^2_\rho \pmtdg{\rho} [\ad,\lambda]) \varphi^\nu)}_0\\ &
\quad+t {(f^\nu,\ze_\nu(\pmtg{\nu})^*\dze^2_\nu \pmtg {\nu}
[\ze_\nu \ze^2_\rho, \pmtdg{\rho}] [\ad,\lambda]
\varphi^\nu)}_{-t}\,\Big)+ {(f, A_{-1} \, \varphi)}_0
\end{split}
\end{equation*}
Before we compute $\tnorm{f}{B \, \varphi}$ notice that each term
in $B$ is a composition of pseudodifferential operators from which
at least one has a symbol supported in $\cp_\rho$. This implies
$\sum_\nu {(\dze_\nu \pmtg{\nu} \ze_\nu f^\nu, \dze_\nu \pmtg{\nu}
\ze_\nu (B \, \varphi)^\nu)}_{-t}$ is a smooth error term. Next,
notice that if we write $B \, = \, \sum_\rho B_\rho,$ for each
$\rho,$ $B_\rho$ is supported in $U_\rho,$ so all terms in the
norm corresponding to neighborhoods $U_\nu$ such that $U_\nu \cap
U_\rho = \emptyset$ are likewise smooth. Thus,
\begin{equation*}
\begin{split}
\tnorm{f}{B \, \varphi} &= \sum_\nu \sum_{\rho \atop {U_\nu \cap
U_\rho \neq \emptyset}} \: \Big({(\dze_\nu \pptg{\nu} \ze_\nu
f^\nu, \dze_\nu \pptg{\nu} \ze_\nu \dze_\rho [\dze_\rho \pptg
{\rho} \ze_\rho, \dbarb]^* \dze_\rho \pptg {\rho} \ze_\rho
\varphi^\nu)}_t\\&\quad+ {(\dze_\nu \potg{\nu} \ze_\nu f^\nu,
\dze_\nu \potg{\nu} \ze_\nu \dze_\rho [\dze_\rho \pptg {\rho}
\ze_\rho, \dbarb]^* \dze_\rho \pptg {\rho} \ze_\rho
\varphi^\nu)}_0\\ &\quad+{(\dze_\nu \pptg {\nu} \ze_\nu
f^\nu,\dze_\nu \pptg{\nu} \ze_\nu \ze_\rho (\pptg {\rho})^*
\dze_\rho [\ad -t [\ad, \lambda], \dze_\rho \pptg {\rho} \ze_\rho]
\dze_\rho \varphi^\nu)}_t\\ &\quad+{(\dze_\nu \potg{\nu} \ze_\nu
f^\nu, \dze_\nu \potg{\nu} \ze_\nu \ze_\rho (\pptg {\rho})^*
\dze_\rho [\ad -t [\ad, \lambda], \dze_\rho \pptg {\rho} \ze_\rho]
\dze_\rho \varphi^\nu)}_0\,\Big)\\&\quad+ smooth \: errors.
\end{split}
\end{equation*}
Again, before we compute $\tnorm{f}{C \, \varphi}$ we note that
each term in $C$ is a composition of pseudodifferential operators
from which at least one has a symbol supported in $\cm_\rho$. By
the same logic as above, then
\begin{equation*}
\begin{split}
\tnorm{f}{C \, \varphi} &= \sum_\nu \sum_{\rho \atop {U_\nu \cap
U_\rho \neq \emptyset}} \: \Big({(\dze_\nu \potg{\nu} \ze_\nu
f^\nu, \dze_\nu \potg{\nu} \ze_\nu \ze_\rho (\pmtg {\rho})^*
\dze_\rho [\ad +t [\ad, \lambda], \dze_\rho \pmtg {\rho} \ze_\rho]
\dze_\rho \varphi^\nu)}_0\\&\quad+ {(\dze_\nu \pmtg{\nu} \ze_\nu
f^\nu, \dze_\nu \pmtg{\nu} \ze_\nu \ze_\rho (\pmtg {\rho})^*
\dze_\rho [\ad +t [\ad, \lambda], \dze_\rho \pmtg {\rho} \ze_\rho]
\dze_\rho \varphi^\nu)}_{-t}\\&\quad+{(\dze_\nu \potg{\nu} \ze_\nu
f^\nu, \dze_\nu \potg{\nu} \ze_\nu \dze_\rho [\dze_\rho \pmtg
{\rho} \ze_\rho, \dbarb]^* \dze_\rho \pmtg {\rho} \ze_\rho
\varphi^\nu)}_0\\&\quad+{(\dze_\nu \pmtg{\nu} \ze_\nu f^\nu,
\dze_\nu \pmtg{\nu} \ze_\nu \dze_\rho [\dze_\rho \pmtg {\rho}
\ze_\rho, \dbarb]^* \dze_\rho \pmtg {\rho} \ze_\rho
\varphi^\nu)}_{-t}\Big)
\\&\quad+ smooth \: errors.
\end{split}
\end{equation*}
Next we add $\tnorm{f}{A \, \varphi},$ and $\tnorm{f}{B \,
\varphi},$ and $\tnorm{f}{C \, \varphi}$ and by comparing this sum
with the expression for $\tnorm {f}{\adt \varphi}$ obtained in
Step $1$, we see that $\tnorm {f}{E_t \, \varphi}$ equals the
following:
\begin{equation*}
\begin{split}
\tnorm{f}{E_t \, \varphi} &= \sum_\nu \: \Big( {(  f^\nu,
\ze_\nu(\potg {\nu})^* \dze_\nu[\ad, \dze_\nu \potg {\nu} \ze_\nu]
\dze_\nu {\varphi^\nu})}_0  + {(f^\nu, \dze_\nu [\dze_\nu \potg
{\nu} \ze_\nu, \dbarb]^* \dze_\nu \potg {\nu} \ze_\nu
{\varphi^\nu})}_0\Big)\\&\quad +\sum_\nu \sum_{\rho \atop {U_\nu
\cap U_\rho
\neq\emptyset}}\:\Big(t{(f^\nu,\ze_\nu(\pptg{\nu})^*\dze^2_\nu
\pptg {\nu} [\ze_\nu \ze^2_\rho, \pptdg{\rho}] [\ad,\lambda]
\varphi^\nu)}_t\\ & \quad- {(f^\nu,\ze_\nu(\potg{\nu})^*\dze^2_\nu
\potg {\nu} \ze_\nu (-t\ze^2_\rho \pptdg{\rho} [\ad,\lambda] +
t\ze^2_\rho \pmtdg{\rho} [\ad,\lambda]) \varphi^\nu)}_0\\ &
\quad-t {(f^\nu,\ze_\nu(\pmtg{\nu})^*\dze^2_\nu \pmtg {\nu}
[\ze_\nu \ze^2_\rho, \pmtdg{\rho}] [\ad,\lambda]
\varphi^\nu)}_{-t}\\&\quad+{(  { f^\nu},\ze_\nu(\potg{\nu})^*
\dze_\nu e^{-t \lambda}\dze_\nu \potg{\nu} \ze_\nu \dze_\rho
[\dze_\rho \pptg {\rho} \ze_\rho, \dbarb]^* \dze_\rho \pptg {\rho}
\ze_\rho {\varphi^\nu})}_0\\&\quad+{(   {
f^\nu},\ze_\nu(\potg{\nu})^* \dze_\nu e^{-t \lambda}\dze_\nu
\potg{\nu} \ze_\nu\ze_\rho (\pptg {\rho})^* \dze_\rho [\ad -t
[\ad, \lambda], \dze_\rho \pptg {\rho} \ze_\rho] \dze_\rho
{\varphi^\nu})}_0 \\
&\quad +  {( { f^\nu}, \ze_\nu (\potg{\nu})^*\dze_\nu e^{t
\lambda}\dze_\nu \potg{\nu} \ze_\nu  \dze_\rho [\dze_\rho \pmtg
{\rho} \ze_\rho, \dbarb]^* \dze_\rho \pmtg {\rho} \ze_\rho
{\varphi^\nu})}_0\\ &\quad +  {( { f^\nu},\ze_\nu (\potg{\nu})^*
\dze_\nu e^{t \lambda}\dze_\nu \potg{\nu} \ze_\nu \ze_\rho (\pmtg
{\rho})^* \dze_\rho [\ad +t [\ad, \lambda], \dze_\rho \pmtg {\rho}
\ze_\rho] \dze_\rho {\varphi^\nu})}_0\\&\quad- {( f^\nu, \ze_\nu
(\potg{\nu})^* \dze^2_\nu \potg{\nu} \ze_\nu \dze_\rho [\dze_\rho
\pptg {\rho} \ze_\rho, \dbarb]^* \dze_\rho \pptg {\rho} \ze_\rho
\varphi^\nu)}_0\\ &\quad-{( f^\nu,\ze_\nu (\potg{\nu})^*
\dze^2_\nu \potg{\nu} \ze_\nu \ze_\rho (\pptg {\rho})^* \dze_\rho
[\ad -t [\ad, \lambda], \dze_\rho \pptg {\rho} \ze_\rho] \dze_\rho
\varphi^\nu)}_0\\&\quad-{( f^\nu, \ze_\nu (\potg{\nu})^*
\dze^2_\nu \potg{\nu} \ze_\nu \ze_\rho (\pmtg {\rho})^* \dze_\rho
[\ad +t [\ad, \lambda], \dze_\rho \pmtg {\rho} \ze_\rho] \dze_\rho
\varphi^\nu)}_0\\&\quad-{( f^\nu, \ze_\nu (\potg{\nu})^*
\dze^2_\nu \potg{\nu} \ze_\nu \dze_\rho [\dze_\rho \pmtg {\rho}
\ze_\rho, \dbarb]^* \dze_\rho \pmtg {\rho} \ze_\rho
\varphi^\nu)}_0\Big)\\
&\quad+\sum_\nu\:\Big({({f^\nu},\ze_\nu[e^{-t\lambda},\pptg
{\nu}]^*\dze_\nu [\ad -t [\ad, \lambda], \dze_\nu \pptg {\nu}
\ze_\nu] \dze_\nu {\varphi^\nu})}_0 \\&\quad+ {( f^\nu, \dze_\nu
[e^{-t \lambda},[\dze_\nu \pptg {\nu} \ze_\nu, \dbarb]]^*\dze_\nu
\pptg
{\nu} \ze_\nu {\varphi^\nu})}_0\\
&\quad+{({f^\nu},\ze_\nu[e^{t\lambda},\pmtg {\nu}]^*\dze_\nu [\ad
+t [\ad, \lambda], \dze_\nu \pmtg {\nu} \ze_\nu] \dze_\nu
{\varphi^\nu})}_0 \\&\quad+ {( f^\nu, \dze_\nu [e^{t
\lambda},[\dze_\nu \pmtg {\nu} \ze_\nu, \dbarb]]^*\dze_\nu \pmtg
{\nu} \ze_\nu {\varphi^\nu})}_0\Big)\\&\quad + (f, S_t \,
\varphi)_0,
\end{split}
\end{equation*}
where $S_t$ includes $S_{t,+}$ and $S_{t,-}$, the
pseudodifferential operators of order $-1$ from the expression for
$\tnorm {f}{\adt \varphi}$, as well as all the other lower order
errors arising from the computations of $\tnorm{f}{A \, \varphi},$
and $\tnorm{f}{B \, \varphi},$ and $\tnorm{f}{C \, \varphi}$.
Using the self-adjoint $t$ dependent pseudodifferential operator
of order zero $G_t$ from Corollary~\ref{equivop} and remembering
that the weight functions for the $t$ norm and the $-t$ norm are
$e^{-t \lambda}$ and $e^{t \lambda}$ respectively,
\begin{equation*}
\begin{split}
E_t &= G_t \bigg( \sum_\nu \: \big(\,\ze_\nu(\potg {\nu})^*
\dze_\nu[\ad, \dze_\nu \potg {\nu} \ze_\nu] \dze_\nu +\dze_\nu
[\dze_\nu \potg {\nu} \ze_\nu, \dbarb]^* \dze_\nu \potg {\nu}
\ze_\nu\,\big)\\&\quad +\sum_\nu \sum_{\rho \atop {U_\nu \cap
U_\rho \neq\emptyset}}\:\big(\,t e^{-t
\lambda}\ze_\nu(\pptg{\nu})^*\dze^2_\nu \pptg {\nu} [\ze_\nu
\ze^2_\rho, \pptdg{\rho}] [\ad,\lambda]\\&\quad
-\ze_\nu(\potg{\nu})^*\dze^2_\nu \potg {\nu} \ze_\nu (-t\ze^2_\rho
\pptdg{\rho} [\ad,\lambda] + t\ze^2_\rho \pmtdg{\rho}
[\ad,\lambda])\\&\quad -t e^{t
\lambda}\ze_\nu(\pmtg{\nu})^*\dze^2_\nu \pmtg {\nu} [\ze_\nu
\ze^2_\rho, \pmtdg{\rho}]
[\ad,\lambda]\\&\quad+\ze_\nu(\potg{\nu})^* \dze_\nu e^{-t
\lambda}\dze_\nu \potg{\nu} \ze_\nu \dze_\rho [\dze_\rho \pptg
{\rho} \ze_\rho, \dbarb]^* \dze_\rho \pptg {\rho} \ze_\rho
\\&\quad+\ze_\nu(\potg{\nu})^* \dze_\nu e^{-t \lambda}\dze_\nu
\potg{\nu} \ze_\nu\ze_\rho (\pptg {\rho})^* \dze_\rho [\ad -t
[\ad, \lambda], \dze_\rho \pptg {\rho} \ze_\rho]
\dze_\rho\\&\quad+\ze_\nu (\potg{\nu})^*\dze_\nu e^{t
\lambda}\dze_\nu \potg{\nu} \ze_\nu  \dze_\rho [\dze_\rho \pmtg
{\rho} \ze_\rho, \dbarb]^* \dze_\rho \pmtg {\rho}
\ze_\rho\\&\quad+\ze_\nu (\potg{\nu})^* \dze_\nu e^{t
\lambda}\dze_\nu \potg{\nu} \ze_\nu \ze_\rho (\pmtg {\rho})^*
\dze_\rho [\ad +t [\ad, \lambda], \dze_\rho \pmtg {\rho} \ze_\rho]
\dze_\rho  \\&\quad-\ze_\nu(\potg{\nu})^* \dze^2_\nu \potg{\nu}
\ze_\nu \dze_\rho [\dze_\rho \pptg {\rho} \ze_\rho, \dbarb]^*
\dze_\rho \pptg {\rho} \ze_\rho
\\&\quad-\ze_\nu(\potg{\nu})^* \dze^2_\nu
\potg{\nu} \ze_\nu\ze_\rho (\pptg {\rho})^* \dze_\rho [\ad -t
[\ad, \lambda], \dze_\rho \pptg {\rho} \ze_\rho]
\dze_\rho\\&\quad-\ze_\nu (\potg{\nu})^* \dze^2_\nu \potg{\nu}
\ze_\nu  \dze_\rho [\dze_\rho \pmtg {\rho} \ze_\rho, \dbarb]^*
\dze_\rho \pmtg {\rho} \ze_\rho\\&\quad-\ze_\nu (\potg{\nu})^*
\dze^2_\nu \potg{\nu} \ze_\nu \ze_\rho (\pmtg {\rho})^* \dze_\rho
[\ad +t [\ad, \lambda], \dze_\rho \pmtg {\rho} \ze_\rho]
\dze_\rho\,
 \big)\\&\quad +\sum_\nu \:\big(\,\ze_\nu[e^{-t\lambda},\pptg
{\nu}]^*\dze_\nu [\ad -t [\ad, \lambda], \dze_\nu \pptg {\nu}
\ze_\nu] \dze_\nu\\&\quad+\ze_\nu[e^{t\lambda},\pmtg
{\nu}]^*\dze_\nu [\ad +t [\ad, \lambda], \dze_\nu \pmtg {\nu}
\ze_\nu] \dze_\nu\\&\quad+\dze_\nu [e^{-t \lambda},[\dze_\nu \pptg
{\nu} \ze_\nu, \dbarb]]^*\dze_\nu \pptg {\nu} \ze_\nu+\dze_\nu
[e^{t \lambda},[\dze_\nu \pmtg {\nu} \ze_\nu, \dbarb]]^*\dze_\nu
\pmtg {\nu} \ze_\nu\,\big)+ S_t \bigg)
\end{split}
\end{equation*}
Notice that $E_t$ is a pseudodifferential operator of order zero
and all of its top order terms are compositions involving
$\potg{\nu},$ its adjoint,
 or one of the commutators $[\ze_\nu \ze^2_\rho,
\pptdg{\rho}]$ or $[\ze_\nu \ze^2_\rho, \pmtdg{\rho}]$, all four
of which have symbols supported in $\co_\nu$ for the first two and
$\co_\rho$ for the other two, which finishes the proof of
Lemma~\ref{adclaim}. \qed

\bigskip \noindent Here is the second proof postponed from
Section~\ref{microlocsection}.

\medskip \noindent {\bf Proof of Lemma~\ref{globalenergy}:} We just need to compute the square of the norm of
$\adt \varphi$. Since $\adt \, = \, A+B+C+E_t,$
\begin{equation*}
\begin{split}
\sqtnorm{\adt \varphi}&= \sqtnorm{A \, \varphi} + \sqtnorm {B \,
\varphi} + \sqtnorm {C \, \varphi} + \sqtnorm{ E_t \,
\varphi}\\&\quad + 2 \, \Re \{ \tnorm{A \, \varphi}{B \, \varphi}
\}+ 2 \, \Re \{ \tnorm{A \, \varphi}{C \, \varphi} \}+ 2 \, \Re \{
\tnorm{A \, \varphi}{E_t \, \varphi} \}\\&\quad+ 2 \, \Re \{
\tnorm{B \, \varphi}{C \, \varphi} \}+ 2 \, \Re \{ \tnorm{B \,
\varphi}{E_t \, \varphi} \}+ 2 \, \Re \{ \tnorm{C \, \varphi}{E_t
\, \varphi} \}
\end{split}
\end{equation*}
By applying Lemma~\ref{sclc} to the cross terms with a very small
$\eps$ for the first three and $\eps = 1$ for the other ones, the
previous expression can be rewritten as follows:
\begin{equation}
\begin{split}
\sqtnorm{\adt \varphi}+\left(1+ \frac{1}{\eps}\right) \sqtnorm {B
\, \varphi}+\left(1+ \frac{1}{\eps}\right)\sqtnorm {C \, \varphi}
+\left(1+ \frac{1}{\eps}\right) \sqtnorm{ E_t \, \varphi}&\geq
(1-3 \, \eps) \sqtnorm{A \, \varphi} \label{adexpr}
\end{split}
\end{equation}
We thus need to compute $\sqtnorm{A \, \varphi},$ $\sqtnorm{B \,
\varphi},$ $\sqtnorm{C \, \varphi},$ and $\sqtnorm{E_t \,
\varphi}.$

\begin{equation*}
\begin{split}
\sqtnorm{A \, \varphi} &=  \sum_\nu \:
\bigg({\bigg|\bigg|\bigg|\dze_\nu \pptg {\nu} \ze_\nu (\ad
-t\sum_\rho \ze^2_\rho \pptdg{\rho} [\ad , \lambda]
+t\sum_\rho\ze^2_\rho \pmtdg{\rho} [\ad , \lambda]) \,
\varphi^\nu\bigg|\bigg|\bigg|}^2_{\,t}\\ & \quad+
{\bigg|\bigg|\dze_\nu \potg {\nu} \ze_\nu (\ad -t\sum_\rho
\ze^2_\rho \pptdg{\rho} [\ad , \lambda] +t\sum_\rho\ze^2_\rho
\pmtdg{\rho} [\ad , \lambda]) \varphi^\nu\bigg|\bigg|}^2_0\\&
\quad+ {\bigg|\bigg|\bigg|\dze_\nu \pmtg {\nu} \ze_\nu (\ad -t
\sum_\rho\ze^2_\rho \pptdg{\rho} [\ad , \lambda]
+t\sum_\rho\ze^2_\rho \pmtdg{\rho} [\ad , \lambda])
\varphi^\nu\bigg|\bigg|\bigg|}^2_{-t}\bigg)
\end{split}
\end{equation*}
Let us analyze first the $t$ norm term. The symbols of
$\pptg{\nu}$ and $\pmtdg{\rho}$ have disjoint supports by
construction, so the pseudodifferential operator $\dze_\nu \pptg
{\nu} \ze_\nu \ze^2_\rho \pmtdg{\rho} [\ad , \lambda]$ is
smoothing. By Lemma~\ref{sclc}, we obtain:
\begin{equation*}
\begin{split}
&{\bigg|\bigg|\bigg|\dze_\nu \pptg {\nu} \ze_\nu (\ad -t\sum_\rho
\ze^2_\rho \pptdg{\rho} [\ad , \lambda] +t\sum_\rho\ze^2_\rho
\pmtdg{\rho} [\ad , \lambda]) \,
\varphi^\nu\bigg|\bigg|\bigg|}^2_{\,t}\\&=
{\bigg|\bigg|\bigg|\dze_\nu \pptg {\nu} \ze_\nu (\ad -t\sum_\rho
\ze^2_\rho \pptdg{\rho} [\ad , \lambda]) \,
\varphi^\nu\bigg|\bigg|\bigg|}^2_{\,t} +t^2 \,
{\bigg|\bigg|\bigg|\sum_\rho \:\dze_\nu \pptg {\nu} \ze_\nu
\ze^2_\rho \pmtdg{\rho} [\ad , \lambda] \,
\varphi^\nu\bigg|\bigg|\bigg|}^2_{\,t}+ cross \: \: term
\\&\geq (1-\eps)\,{\bigg|\bigg|\bigg|\dze_\nu \pptg {\nu} \ze_\nu (\ad
-t\sum_\rho \ze^2_\rho \pptdg{\rho} [\ad , \lambda]) \,
\varphi^\nu\bigg|\bigg|\bigg|}^2_{\,t} +O_t (||\varphi||^2_{-1})
\end{split}
\end{equation*}
where $\eps$ is a small positive number. Next,
\begin{equation*}
\begin{split}
&{\bigg|\bigg|\bigg|\dze_\nu \pptg {\nu} \ze_\nu (\ad -t\sum_\rho
\ze^2_\rho \pptdg{\rho} [\ad , \lambda]) \,
\varphi^\nu\bigg|\bigg|\bigg|}^2_{\,t}\\&={\bigg|\bigg|\bigg|\dze_\nu
\pptg {\nu}  (\ad \ze_\nu -t\sum_\rho \pptdg{\rho} \ze_\nu
\ze^2_\rho [\ad , \lambda]-t\sum_\rho  [\ze_\nu \ze^2_\rho,
\pptdg{\rho}][\ad , \lambda]+[\ze_\nu, \ad]) \,
\varphi^\nu\bigg|\bigg|\bigg|}^2_{\,t}\\&\geq
(1-\eps)\,{\bigg|\bigg|\bigg|\dze_\nu \pptg {\nu}  (\ad \ze_\nu
-t\sum_\rho \pptdg{\rho} \ze_\nu \ze^2_\rho [\ad , \lambda]) \,
\varphi^\nu\bigg|\bigg|\bigg|}^2_{\,t}\\&\quad-(\frac{1}{\eps}-1)\,{\bigg|\bigg|\bigg|\dze_\nu
\pptg {\nu}  ([\ze_\nu, \ad]-t\sum_\rho  [\ze_\nu \ze^2_\rho,
\pptdg{\rho}][\ad , \lambda]) \,
\varphi^\nu\bigg|\bigg|\bigg|}^2_{\,t}
\end{split}
\end{equation*}
Since $\pptdg{\rho}$ is inverse zero order $t$ dependent, $-t
[\ze_\nu \ze^2_\rho, \pptdg{\rho}]$ can be bounded independently
of $t$. $[\ze_\nu, \ad]$ is of order zero and independent of $t$,
so
$${\bigg|\bigg|\bigg|\dze_\nu \pptg {\nu}  ([\ze_\nu, \ad]-t\sum_\rho
[\ze_\nu \ze^2_\rho, \pptdg{\rho}][\ad , \lambda]) \,
\varphi^\nu\bigg|\bigg|\bigg|}^2_{\,t} = O(|||\varphi|||^2_{\,
t}).$$ As in the proof of Lemma~\ref{adclaim}, by part (ii) of
Property $\dag$ on page~\pageref{propdagger}, $\pptg {\nu}
\pptdg{\rho} \, = \, \pptg {\nu}$ up to lower order terms for all
$\rho$ such that $U_\nu \cap U_\rho \neq \emptyset$ and $\pptg
{\nu} \pptdg{\rho}$ is a smoothing operator by
Lemma~\ref{disjlemma} otherwise. The terms containing these lower
order errors can be thrown into $O_t (||\varphi||^2_{-1}).$ Thus
we can split ${||\dze_\nu \pptg {\nu} (\ad \ze_\nu -t\sum_\rho
\pptdg{\rho} \ze_\nu \ze^2_\rho [\ad , \lambda]) \,
\varphi^\nu||}^2_{\,t}$ into the relevant terms and the lower
order ones and take care of the cross term using Lemma~\ref{sclc}
to obtain
\begin{equation*}
\begin{split}
&{\bigg|\bigg|\bigg|\dze_\nu \pptg {\nu} \ze_\nu (\ad -t\sum_\rho
\ze^2_\rho \pptdg{\rho} [\ad , \lambda]) \,
\varphi^\nu\bigg|\bigg|\bigg|}^2_{\,t}\\&\geq
(1-2\,\eps)\,{\bigg|\bigg|\bigg|\dze_\nu \pptg {\nu}  (\ad \ze_\nu
-t\sum_\rho \ze_\nu \ze^2_\rho [\ad , \lambda]) \,
\varphi^\nu\bigg|\bigg|\bigg|}^2_{\,t}+O(|||\varphi|||^2_{\,
t})+O_t (||\varphi||^2_{-1}).
\end{split}
\end{equation*}
Finally, $\sum_\rho \ze^2_\rho \, = \,1,$ so
\begin{equation*}
\begin{split}
&{\bigg|\bigg|\bigg|\dze_\nu \pptg {\nu} \ze_\nu (\ad -t\sum_\rho
\ze^2_\rho \pptdg{\rho} [\ad , \lambda]) \,
\varphi^\nu\bigg|\bigg|\bigg|}^2_{\,t}\\&\geq
(1-2\,\eps)\,{|||\dze_\nu \pptg {\nu}  (\ad  -t [\ad , \lambda])
\,\ze_\nu \varphi^\nu|||}^2_{\,t}+O(|||\varphi|||^2_{\, t})+O_t
(||\varphi||^2_{-1}).
\end{split}
\end{equation*}
So then we sum over $\nu$ to obtain
\begin{equation}
\begin{split}
&\sum_\nu{\bigg|\bigg|\bigg|\dze_\nu \pptg {\nu} \ze_\nu (\ad
-t\sum_\rho \ze^2_\rho \pptdg{\rho} [\ad , \lambda]
+t\sum_\rho\ze^2_\rho \pmtdg{\rho} [\ad , \lambda]) \,
\varphi^\nu\bigg|\bigg|\bigg|}^2_{\,t}\\&\geq (1-\eps')\,\sum_\nu
\: {|||\dze_\nu \pptg {\nu}  (\ad  -t [\ad , \lambda]) \,\ze_\nu
\varphi^\nu|||}^2_{\,t}+O(|||\varphi|||^2_{\, t}) +O_t
(||\varphi||^2_{-1})\label{Aposnorm}
\end{split}
\end{equation}
Similarly,
\begin{equation}
\begin{split}
&\sum_\nu{\bigg|\bigg|\bigg|\dze_\nu \pmtg {\nu} \ze_\nu (\ad
-t\sum_\rho \ze^2_\rho \pptdg{\rho} [\ad , \lambda]
+t\sum_\rho\ze^2_\rho \pmtdg{\rho} [\ad , \lambda]) \,
\varphi^\nu\bigg|\bigg|\bigg|}^2_{-t}\\&\geq (1-\eps')\,\sum_\nu
\: {|||\dze_\nu \pmtg {\nu}  (\ad  +t [\ad , \lambda]) \,\ze_\nu
\varphi^\nu|||}^2_{-t}+O(|||\varphi|||^2_{-t}) +O_t
(||\varphi||^2_{-1})\label{Anegnorm}
\end{split}
\end{equation}
Now, let us now deal with the zero norm term in the expression for
$\sqtnorm{A \, \varphi}$ :
\begin{equation*}
\begin{split}
&{\bigg|\bigg|\dze_\nu \potg {\nu} \ze_\nu (\ad -t\sum_\rho
\ze^2_\rho \pptdg{\rho} [\ad , \lambda] +t\sum_\rho\ze^2_\rho
\pmtdg{\rho} [\ad , \lambda]) \, \varphi^\nu\bigg|\bigg|}^2_0 =
{||\dze_\nu \potg {\nu} \ze_\nu \ad \varphi^\nu||}^2_0
\\ & \quad+t^2{\bigg|\bigg|\sum_\rho\dze_\nu \potg {\nu} \ze_\nu
\ze^2_\rho \pptdg{\rho} [\ad , \lambda]
\varphi^\nu\bigg|\bigg|}^2_0+t^2{\bigg|\bigg|\sum_\rho\dze_\nu
\potg {\nu} \ze_\nu \ze^2_\rho \pmtdg{\rho} [\ad , \lambda]
\varphi^\nu\bigg|\bigg|}^2_0 \\ & \quad- 2 t\sum_\rho \Re \{
{(\dze_\nu \potg {\nu} \ze_\nu \ad \varphi^\nu,\dze_\nu \potg
{\nu} \ze_\nu \ze^2_\rho \pptdg{\rho} [\ad,\lambda]
\varphi^\nu)}_0 \}\\ & \quad+ 2 t\sum_\rho \Re \{ {(\dze_\nu \potg
{\nu} \ze_\nu \ad \varphi^\nu,\dze_\nu \potg {\nu} \ze_\nu
\ze^2_\rho \pmtdg{\rho} [\ad,\lambda]
\varphi^\nu)}_0 \}\\
&  \quad- 2\, t^2\sum_{\rho,\, \mu} Re \{ {(\dze_\nu \potg {\nu}
\ze_\nu \ze^2_\rho \pptdg{\rho} [\ad,\lambda] \varphi^\nu,\dze_\nu
\potg {\nu} \ze_\nu \ze^2_\mu \pmtdg{\mu} [\ad,\lambda]
\varphi^\nu)}_0 \}
\end{split}
\end{equation*}
Now using Lemma~\ref{sclc} on page~\pageref{sclc} and summing over
$\nu,$ we see that for some small $\eps>0$
\begin{equation}
\begin{split}
&\sum_\nu \: {\bigg|\bigg|\dze_\nu \potg {\nu} \ze_\nu (\ad
-t\sum_\rho \ze^2_\rho \pptdg{\rho} [\ad , \lambda]
+t\sum_\rho\ze^2_\rho \pmtdg{\rho} [\ad , \lambda]) \,
\varphi^\nu\bigg|\bigg|}^2_0 \\&\geq (1-2\,\eps)\,\sum_\nu
{||\dze_\nu \potg {\nu} \ze_\nu \ad \varphi||}^2_0 -\frac
{1}{\epsilon} \, t^2\sum_\nu{\bigg|\bigg|\sum_\rho\dze_\nu \potg
{\nu} \ze_\nu \ze^2_\rho \pptdg{\rho} [\ad , \lambda]
\varphi^\nu\bigg|\bigg|}^2_0\\ & \quad- \frac {1}{\epsilon} \,
t^2\sum_\nu{\bigg|\bigg|\sum_\rho\dze_\nu \potg {\nu} \ze_\nu
\ze^2_\rho \pmtdg{\rho} [\ad , \lambda]
\varphi^\nu\bigg|\bigg|}^2_0. \label{Azeronorm}
\end{split}
\end{equation}
Put together \ref{Aposnorm}, \ref{Anegnorm}, and \ref{Azeronorm}
to obtain
\begin{equation*}
\begin{split}
\sqtnorm{A\, \varphi} &\geq (1-\eps')\,\bigg[\sum_\nu \:
{|||\dze_\nu \pptg {\nu}  (\ad  -t [\ad , \lambda]) \,\ze_\nu
\varphi^\nu|||}^2_{\,t}+\sum_\nu \: {||\dze_\nu \potg {\nu}
\ze_\nu \ad \varphi^\nu||}^2_0 \\ & \quad+ \sum_\nu \:
{|||\dze_\nu \pmtg {\nu}  (\ad  +t [\ad , \lambda]) \,\ze_\nu
\varphi^\nu|||}^2_{-t}\bigg] -\frac {1}{\epsilon} \,
t^2\sum_\nu{\bigg|\bigg|\sum_\rho\dze_\nu \potg {\nu} \ze_\nu
\ze^2_\rho \pptdg{\rho} [\ad , \lambda]
\varphi^\nu\bigg|\bigg|}^2_0\\ & \quad- \frac {1}{\epsilon} \,
t^2\sum_\nu{\bigg|\bigg|\sum_\rho\dze_\nu \potg {\nu} \ze_\nu
\ze^2_\rho \pmtdg{\rho} [\ad , \lambda]
\varphi^\nu\bigg|\bigg|}^2_0+O({\langle |\varphi|\rangle}^2_t)
+O_t (||\varphi||^2_{-1})
\end{split}
\end{equation*}
As in Section~\ref{defnot}, let $\ad  -t [\ad , \lambda] \, = \,
\adp$, $\ad  +t [\ad , \lambda] \, = \, \adm$, and $\ad \, = \,
\ado$ so then
\begin{equation}
\begin{split}
\sqtnorm{A\, \varphi} &\geq (1-\eps')\,\sum_\nu \:\bigg[
{|||\dze_\nu \pptg {\nu}  \adp\ze_\nu \varphi^\nu|||}^2_{\,t}+
{||\dze_\nu \potg {\nu} \ze_\nu \ado \varphi^\nu||}^2_0+
{|||\dze_\nu \pmtg {\nu} \adm\ze_\nu \varphi^\nu|||}^2_{-t}\bigg]
\\ & \quad-\frac {1}{\epsilon} \,
t^2\sum_\nu{\bigg|\bigg|\sum_\rho\dze_\nu \potg {\nu} \ze_\nu
\ze^2_\rho \pptdg{\rho} [\ad , \lambda]
\varphi^\nu\bigg|\bigg|}^2_0- \frac {1}{\epsilon} \,
t^2\sum_\nu{\bigg|\bigg|\sum_\rho\dze_\nu \potg {\nu} \ze_\nu
\ze^2_\rho \pmtdg{\rho} [\ad , \lambda]
\varphi^\nu\bigg|\bigg|}^2_0\\ & \quad+O({\langle
|\varphi|\rangle}^2_t) +O_t (||\varphi||^2_{-1})\label{Atotalexpr}
\end{split}
\end{equation}
Since $$B =\sum_\rho \: (\dze_\rho [\dze_\rho \pptg {\rho}
\ze_\rho, \dbarb]^* \dze_\rho \pptg {\rho} \ze_\rho +\ze_\rho
(\pptg {\rho})^* \dze_\rho [\ad -t [\ad, \lambda], \dze_\rho \pptg
{\rho} \ze_\rho] \dze_\rho)$$ and $\pptg{\rho}$ is inverse zero
order $t$ dependent, $[\ad -t [\ad, \lambda], \dze_\rho \pptg
{\rho} \ze_\rho]$ can be bounded independently of $t$ and
$$\sqtnorm{B \, \varphi} =O({\langle |\varphi|\rangle}^2_t).$$
Similarly, since
$$C = \sum_\rho \: (\ze_\rho (\pmtg {\rho})^* \dze_\rho [\ad +t
[\ad, \lambda], \dze_\rho \pmtg {\rho} \ze_\rho] \dze_\rho
+\dze_\rho [\dze_\rho \pmtg {\rho} \ze_\rho, \dbarb]^* \dze_\rho
\pmtg {\rho} \ze_\rho)$$ and $\pmtg{\rho}$ is inverse zero order
$t$ dependent, $[\ad +t [\ad, \lambda], \dze_\rho \pmtg {\rho}
\ze_\rho]$ can be bounded independently of $t$ and $$\sqtnorm{C \,
\varphi} =O({\langle |\varphi|\rangle}^2_t).$$ Lemma~\ref{adclaim}
established that $E_t$ was of order zero and its principal symbol
was supported in $\co_\nu$. The lower order terms are controlled
by $O_t (||\varphi||^2_{-1})$, as for the rest, let $\potdg{\nu}$
dominate $\potg{\nu}$ and be also supported in $\co_\nu$, then
there exists a $t$ dependent constant $K_t$ such that $$K_t
\sum_\nu \:{||\dze_\nu \potdg {\nu} \ze_\nu \varphi^\nu||}^2_0+O_t
(||\varphi||^2_{-1}) \geq \sqtnorm{E_t \, \varphi}.$$ Notice that
$K_t$ grows exponentially with $t$ since the sum it multiplies has
to dominate squares of $t$ and $-t$ norms of terms of $E_t$ and
these norms have weight functions which are exponential in $t$.
For $K_t$ large enough, $K_t \sum_\nu \:{||\dze_\nu \potdg {\nu}
\ze_\nu \varphi^\nu||}^2_0$ also controls $$\frac {1}{\epsilon} \,
t^2\sum_\nu{\bigg|\bigg|\sum_\rho\dze_\nu \potg {\nu} \ze_\nu
\ze^2_\rho \pptdg{\rho} [\ad , \lambda]
\varphi^\nu\bigg|\bigg|}^2_0$$ and $$\frac {1}{\epsilon} \,
t^2\sum_\nu{\bigg|\bigg|\sum_\rho\dze_\nu \potg {\nu} \ze_\nu
\ze^2_\rho \pmtdg{\rho} [\ad , \lambda]
\varphi^\nu\bigg|\bigg|}^2_0,$$ two of the terms in the
expression~\ref{Atotalexpr}. We conclude that $\sqtnorm{\adt
\varphi}$ must satisfy
\begin{equation}
\begin{split}
&K \sqtnorm{\adt \varphi}+K_t \sum_\nu \:{||\dze_\nu \potdg {\nu}
\ze_\nu \varphi^\nu||}^2_0+O({\langle |\varphi|\rangle}^2_t) +O_t
(||\varphi||^2_{-1})\\&\geq \sum_\nu \: {|||\dze_\nu \pptg {\nu}
\adp\ze_\nu \varphi^\nu|||}^2_{\,t}+ \sum_\nu \:{||\dze_\nu \potg
{\nu} \ze_\nu \ado \varphi^\nu||}^2_0+ \sum_\nu \:{|||\dze_\nu
\pmtg {\nu} \adm\ze_\nu \varphi^\nu|||}^2_{-t},
\end{split}
\end{equation}
where $K$ is independent of $t$ $\left(K \, = \, \frac {1}{(1-3 \,
\eps)(1- \eps')}\right)$ and $K_t$ has been slightly increased.
The lemma follows. \qed

\end{document}